\documentclass[final,3p,times,authoryear]{elsarticle}
\usepackage{amssymb,amsmath}
\usepackage{hyperref}
\usepackage{adjustbox}
\usepackage{tabularx}
\usepackage{rotating}
\usepackage{float}
\usepackage[super]{nth}

\usepackage{comment}
\journal{...}

\begin{document}

\begin{frontmatter}

%% Title, authors and addresses

%% use the tnoteref command within \title for footnotes;
%% use the tnotetext command for theassociated footnote;
%% use the fnref command within \author or \affiliation for footnotes;
%% use the fntext command for theassociated footnote;
%% use the corref command within \author for corresponding author footnotes;
%% use the cortext command for theassociated footnote;
%% use the ead command for the email address,
%% and the form \ead[url] for the home page:
%% \title{Title\tnoteref{label1}}
%% \tnotetext[label1]{}
%% \author{Name\corref{cor1}\fnref{label2}}
%% \ead{email address}
%% \ead[url]{home page}
%% \fntext[label2]{}
%% \cortext[cor1]{}
%% \affiliation{organization={},
%%            addressline={}, 
%%            city={},
%%            postcode={}, 
%%            state={},
%%            country={}}
%% \fntext[label3]{}

\title{Two-echelon Electric Vehicle Routing Problem in Parcel Delivery: A Literature Review}

%% use optional labels to link authors explicitly to addresses:
%% \author[label1,label2]{}
%% \affiliation[label1]{organization={},
%%             addressline={},
%%             city={},
%%             postcode={},
%%             state={},
%%             country={}}
%%
%% \affiliation[label2]{organization={},
%%             addressline={},
%%             city={},
%%             postcode={},
%%             state={},
%%             country={}}

\author[1,a]{Nima Moradi\footnote{Corresponding author: 251 Rue St-Marguerite, H4C2W7, Montreal, Canada, nima.moradi@mail.concordia.ca}}
\author[2,a]{Niloufar Mirzavand Boroujeni}
\author[3]{Navid Aftabi}
\author[4]{Amin Aslani}

\affiliation[1]{organization={Concordia Institute for Information Systems Engineering},%Department and Organization
            %addressline={}, 
            city={Montreal},
            postcode={H3G1M8}, 
            state={QC},
            country={Canada, nima.moradi@mail.concordia.ca}}
\affiliation[2]{organization={University of Minnesota, Department of Industrial and Systems Engineering},%Department and Organization
            %addressline={}, 
            city={Minneapolis},
            postcode={55455}, 
            state={MN},
            country={United States, mirza047@umn.edu}}
\affiliation[3]{organization={University of Wisconsin-Madison, Industrial and Systems Engineering Department},%Department and Organization
            %addressline={}, 
            city={Madison},
            postcode={53707}, 
            state={WI},
            country={United States, aftabi@wisc.edu}}
\affiliation[4]{organization={University of Wisconsin-Whitewater, College of Business and Economics},%Department and Organization
            %addressline={}, 
            city={Whitewater},
            postcode={53190}, 
            state={WI},
            country={United States, aslania@uww.edu}}
\affiliation[a]{country={These authors equally contributed to the present work and shared the first authorship.}}

\begin{abstract}
Multi-echelon parcel delivery systems using electric vehicles (EVs) are crucial for managing urban logistics complexity and promoting sustainability. In multi-echelon systems, particularly within two-stage systems, larger vehicles transport parcels from a central depot to satellite hubs, where smaller EVs pick up the parcels and carry out last-mile deliveries. This system could increase efficiency, reduce emissions, and improve service reliability. The two-echelon electric vehicle routing problem (2E-EVRP), an extension of the traditional two-echelon vehicle routing problem (2E-VRP), addresses EV-specific challenges such as battery constraints and recharging stations to tackle environmental impacts, urban congestion, and e-commerce demands. While effectively reducing costs, energy use, and emissions, the 2E-EVRP faces modeling challenges due to multi-echelon structures, EV limitations, and recharging station selection. This paper systematically reviews 2E-EVRP literature, analyzing key studies. It proposes a classification scheme to categorize the papers based on the problem variants, objectives, constraints, and solution methods. It identifies gaps such as delivery tardiness, environmental trade-offs, multi-objective optimization, multiple depots, split deliveries, and time-dependent travel conditions. Future research directions include aligning models with urban policies, integrating parcel lockers, enabling same-day delivery, and incorporating advanced technologies like autonomous vehicles. Methodological advancements suggest using machine learning, reinforcement learning, and simulation-based approaches to enhance dynamic routing and real-time decision-making. These directions aim to expand the 2E-EVRP applicability, addressing theoretical and practical challenges in sustainable urban logistics for future works.

\end{abstract}

%%Graphical abstract
%\begin{graphicalabstract}
%\includegraphics{grabs}
%\end{graphicalabstract}

%%Research highlights
\begin{highlights}
\item The first review work on two-echelon electric vehicle routing problems.
\item A classification scheme: problem variants, objectives, constraints, solution methods.  
\item Research gaps are identified based on collected papers.  
\item Future research directions: problem extensions, real-life applications, methodology.  
\end{highlights}

\begin{keyword}
%% keywords here, in the form: keyword \sep keyword
Two-echelon routing \sep electric vehicles \sep parcel delivery \sep review
%% PACS codes here, in the form: \PACS code \sep code

%% MSC codes here, in the form: \MSC code \sep code
%% or \MSC[2008] code \sep code (2000 is the default)

\end{keyword}

\end{frontmatter}

%% \linenumbers

%% main text
\section{Introduction}
Multi-echelon parcel delivery systems in urban settings have emerged as a critical logistics model to handle the increasing complexity of urban transportation \citep{kayvanfar2018analysis}. These systems utilize a hierarchical structure, where goods are transported in two stages to intermediate facilities, such as satellites or urban consolidation centers, and then to final destinations. This approach allows for better route optimization, reduced vehicle usage in urban cores, and improved service reliability \citep{nielsen2024systematic}. Incorporating electric vehicles (EVs) into these systems offers additional benefits, such as reduced emissions and noise pollution, making them an increasingly important solution for sustainable urban logistics. According to the International Energy Agency (IEA), EVs emit 50\% less greenhouse gases on average over their lifecycle compared to internal combustion engine vehicles, even when accounting for electricity generation, making them a critical tool in reducing global CO2 emissions and combating climate change \citep{iea}.

The Two-Echelon Electric Vehicle Routing Problem (2E-EVRP), initially introduced by \cite{breunig2017two}, has garnered increasing attention in recent years due to its relevance in optimizing sustainable logistics, especially parcel delivery systems. Integrating EVs into two-tier delivery networks addresses key challenges such as reducing environmental impact, managing urban congestion, and meeting the growing demands of e-commerce-driven last-mile delivery \citep{breunig2019electric,siragusa2022electric}. Building on multi-echelon delivery systems, the 2E-EVRP extends the traditional two-echelon vehicle routing problem (2E-VRP), initially introduced by \cite{crainic2009models}, by introducing EVs and recharging stations in the second echelon. In the first echelon, larger vehicles transport goods from a central depot to satellites, while in the second echelon, smaller EVs handle last-mile deliveries to customers. This two-echelon framework supports operational efficiency and addresses sustainability challenges by leveraging EVs for urban delivery, which are better suited for short-distance trips with strict environmental constraints. Ryder highlights that adopting sustainable practices in last-mile delivery enhances brand reputation and fosters customer loyalty, indicating a competitive advantage for companies embracing eco-friendly initiatives \citep{ryder}. 

Applications of the 2E-EVRP are particularly relevant in urban logistics, where multi-echelon frameworks enhance operational efficiency while adhering to urban planning regulations and environmental targets. Real-world data demonstrates the effectiveness of these systems in reducing delivery costs, energy consumption, and greenhouse gas emissions \citep{benhassine2023optimization}. Companies in the logistics and e-commerce sectors are increasingly adopting such models to meet sustainability goals while maintaining customer satisfaction. DHL reports that implementing green logistics strategies, including multi-echelon systems, addresses environmental concerns and aligns with customer preferences for sustainable practices \citep{dhl}. Despite the advantages, significant challenges remain in modeling and solving 2E-EVRPs due to the complexities introduced by multi-echelon structures, EV-specific constraints like battery capacities and charging station (CS) locations, and urban delivery restrictions. In the literature, no review work on 2E-EVRP and its variants discusses the type of problems, their properties, proposed solutions methodologies, scalability, and a generic classification framework from an operations research (OR) perspective. The lack of a consolidated review limits the ability to synthesize knowledge across studies and identify future research directions. To fill this gap, this paper offers a thorough review of the 2E-EVRP in the parcel delivery context. It aims to analyze the current state of research, identify emerging trends, and propose directions for future work to advance the field.

This review employs a systematic methodology, categorizing existing literature based on problem variants, objective functions, operational constraints, and solution techniques related to the 2E-EVRPs from the OR viewpoint. It also examines mathematical models, datasets, solved instances, and practical implementations to provide a holistic understanding of the 2E-EVRP landscape. This work differs from the other review papers on related problems. For example, \cite{gonzalez2011two,prodhon2014survey,cuda2015survey,sluijk2023two,nielsen2024systematic} have reviewed the 2E-VRP or two-echelon location-routing problem without elaboration on EV features or recharging stations and related challenges and limitations. They mainly consider the problems in which conventional vehicles are used in both echelons. Also, \cite{li2021ground,yu2024collaborative} only focused on conventional ground vehicles in the first echelon and drones or robots in the second echelon without using the recharging stations in the problems or models. The contributions of this paper include a detailed literature analysis, a review and analysis of the problems, mathematical models and solution methods, the identification of challenges and trends, and guidance for future studies on the 2E-EVRPs. The findings aim to bridge gaps in the literature and inspire novel problem extensions, models, and efficient solutions for sustainable parcel delivery in two-tier last-mile logistics.
 
The remainder of the paper is structured as follows: Section \ref{review-method} provides the review methodology applied in this survey paper. Section \ref{descriptive-class} presents the descriptive analysis and classification scheme based on the collected works. Section \ref{problem-type-variant} explains the problem types and variants introduced in the literature of 2-EVRPs besides the related mathematical models, discussing each problem's assumptions, features, and constraints. Section \ref{obj-functions-section} discusses the objective functions presented in the literature on 2E-EVRPs. Section \ref{solution-methodology} describes the solution methodology developed in the related works for 2E-EVRPs, including their scalability and size of solved instances. Section \ref{gaps-future} provides the study gaps and future research directions. Finally, Section \ref{conclusion} presents the paper's conclusion.    

\section{Review Methodology}\label{review-method}
\begin{figure}[ht!]
    \centering
    \includegraphics[width=0.9\linewidth]{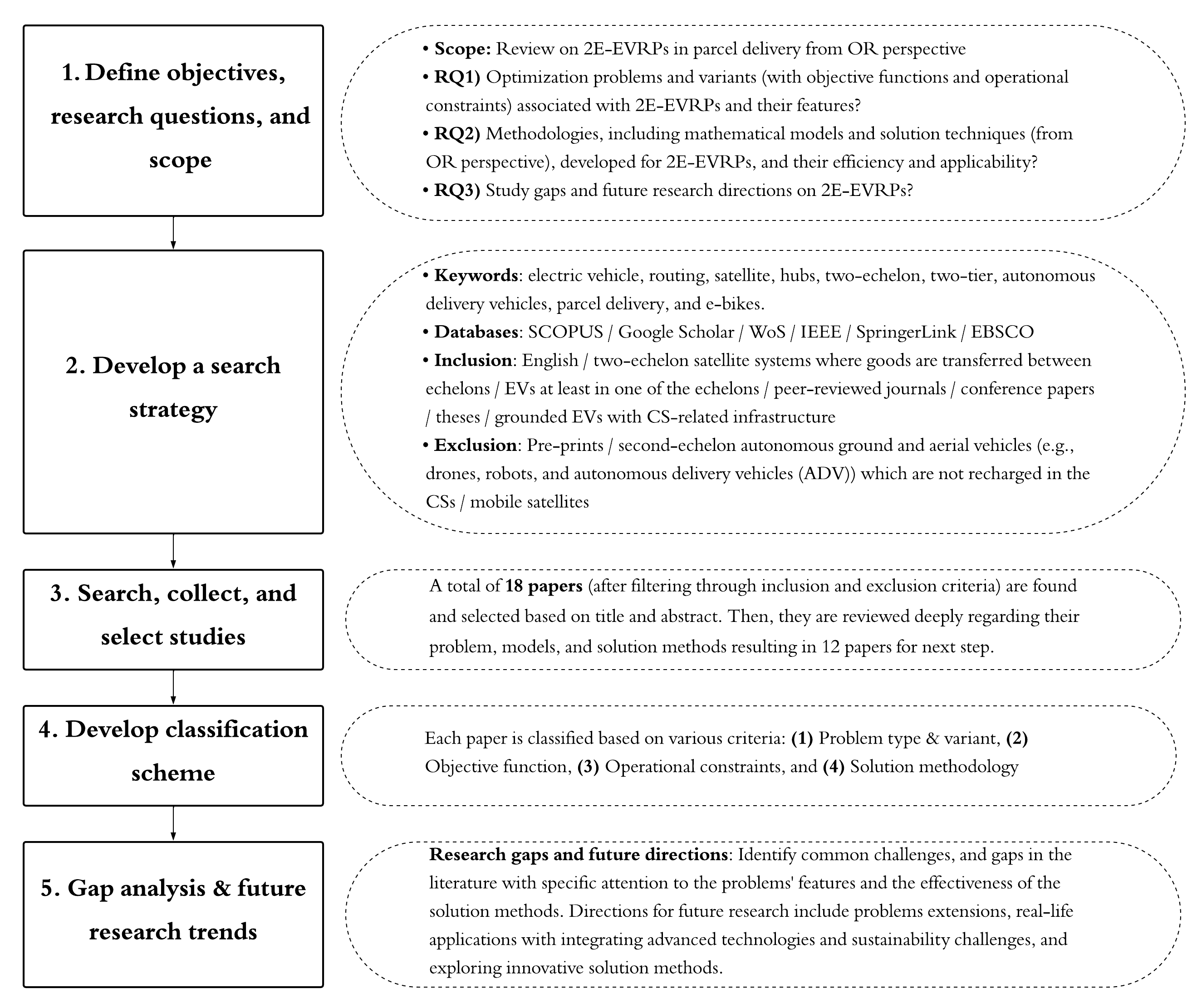}
    \caption{A description of the review methodology steps in the present paper}
    \label{steps}
\end{figure}

In this work, we apply a review methodology similar to SCRP \cite{durach2017new}, based on Systematic Literature Review (SLR). SLR is a structured and systematic approach to reviewing and synthesizing research on a specific topic, often used in academic and professional studies to provide a comprehensive overview of existing knowledge, identify gaps, and propose future directions. It consists of five main steps, as shown in Fig.~\ref{steps}, described in more detail in the following.

\begin{itemize}
    \item Define objectives and scope: We began by defining the review's objectives (research questions) and outlining its scope, delineating specific research areas and topics.
    \item Develop a search strategy: We established keywords, selected databases, and set inclusion/exclusion criteria to guide our search and selection process.
    \item Search, collect, and select studies: Relevant studies were systematically searched using the previous step and selected based on screening processes.
    \item Develop a classification scheme: A framework was developed to categorize the studies, aiding in the structured analysis of findings.
    \item Gap analysis \& future research trends: Models, variants, and methods were compared, and gaps were identified. Also, key findings were summarized, and recommendations for future research were outlined.
\end{itemize}

\subsection{Define objectives and scope}
This review aims to comprehensively analyze the 2E-EVRP literature, focusing on problem classifications, solution methodologies, and emerging trends. The scope encompasses studies involving optimization of the two-echelon delivery systems, in which the first-echelon vehicles transport goods to intermediate depots (satellites or hubs), and multiple EVs in the second echelon complete the deliveries while recharging at CSs if needed. In other words, specific attention is given to two-tier delivery systems using EVs in the second echelon, where the mathematical modeling and OR-based solution methodologies are presented. In addition, it is important to define the research questions at the start of a review study, as they provide the foundation for structuring the study and help shape the anticipated outcomes. In this current review, we have expressed the research questions in the following manner:

\begin{itemize}
\item \textbf{RQ1)} What are the main optimization problems and variants (with objective functions and operational constraints) associated with implementing 2E-EVRPs in urban logistics, and what are their features?
\item \textbf{RQ2)} What methodologies, including mathematical models and solution techniques (from OR perspective), have been developed to address 2E-EVRPs, and how do they compare efficiency and applicability?
\item \textbf{RQ3)} What research gaps and future directions exist for advancing the study of 2E-EVRPs, particularly in problem extensions, methodologies, and real-life applications integrating sustainability and emerging technologies?
\end{itemize}  

\subsection{Develop a search strategy}
The following keywords were carefully selected to encompass the domain and its variations. They guided the search process: electric vehicle, routing, satellite, hubs, two-echelon, two-tier, autonomous delivery vehicles, parcel delivery, and e-bikes. Scientific databases, including SCOPUS, IEEE Xplore, Google Scholar, Web of Science, EBSCO, and SpringerLink, were chosen for their comprehensive subject coverage. The following inclusion and exclusion criteria were also selected to narrow down our search process and paper selection. 

\begin{itemize}
\item \textbf{Inclusion:} English studies focusing on two-echelon satellite systems where goods are transferred between echelons, use of EVs at least in one of the echelons, peer-reviewed journals, conference papers, and theses involving grounded EV with CS-related infrastructure.
\item \textbf{Exclusion:} Pre-prints, studies involving second-echelon autonomous ground and aerial vehicles (e.g., drones, robots, and autonomous delivery vehicles (ADV)) which are not recharged in the CSs during the routes, and mobile satellites.
\end{itemize}

\subsection{Search, collect, and select studies}
After searching and finding the papers according to the previous step, the initial selection was based on titles and abstracts (resulting in 18 papers). The papers were then reviewed thoroughly to ensure relevance to the study's scope, and finally, 12 papers were sent to the next step for deep analysis and providing the classification scheme.

\subsection{Develop classification scheme}
Studies were categorized based on themes such as problem types (variants), objective functions at each problem, constraints, and solution methodologies. This classification facilitated structured comparison and analysis (please see Section \ref{descriptive-class} for more details on the proposed classification framework).

\subsection{Gap analysis \& future research trends}
The collected studies were compared to identify patterns, common challenges, and gaps in the literature. Specific attention was paid to the problems' features and the effectiveness of the solution methods employed. The review's main findings were summarized, highlighting gaps and proposing directions for future research. These include integrating advanced technologies, addressing sustainability challenges, and exploring innovative solutions. The study gaps and future research suggestions are detailed in Section \ref{gaps-future}. This structured methodology and its steps ensured a rigorous and comprehensive review of the literature on 2E-EVRPs, contributing to both academic knowledge and practical applications.

\section{Descriptive analysis \& classification scheme}\label{descriptive-class}
This section first conducts a descriptive analysis of the twelve reviewed articles. This analysis examines several variables, including publication years, journal sources, and solution methodologies. By exploring these factors, we aim to highlight areas of significant research focus, emerging trends, and potential gaps in the literature on 2E-EVRP. Table~\ref{tab:publications-journals-variants} summarizes key studies on 2E-EVRP. The authors' names are listed alongside the corresponding years of publication, offering a chronological view of the progression in this research area. The table also identifies the journals, highlighting the diversity of academic platforms contributing to disseminating 2E-EVRP-related findings.

\begin{table}[ht!]
\centering
\caption{Twelve published articles in the domain of 2E-EVRP}
\resizebox{1\textwidth}{!}{
\begin{tabular}{|l|l|l|}
\hline
\textbf{Authors and Year} & \textbf{Journal} & \textbf{Publisher}\\ \hline
\cite{breunig2017two}       & XLIX Simpósio Brasileiro de Pesquisa Operacional & N/A \\ \hline
\cite{agardi2019two}        & Transport and Telecommunication & Transport and Telecommunication Institute, Lomonosova\\ \hline
\cite{breunig2019electric}  & Computers and Operations Research & Elsevier\\ \hline
\cite{jie2019two}           & European Journal of Operational Research & Elsevier\\ \hline
%\cite{poeting2019comprehensive} & 2019 Winter Simulation Conference (WSC)\\ \hline
%\cite{poeting2019simulation}    & Advances in Production, Logistics and Traffic \\ \hline
\cite{wang2019two}          & Journal of Physics: Conference Series & IOP Publishing \\ \hline
%\cite{lii2020two}           & Computers \& Industrial Engineering \\ \hline
\cite{affi2020general}      & 2020 International Multi-Conference on OCTA & IEEE \\ \hline
%\cite{liu2020two}           & IEEE Access & IEEE\\ \hline
%\cite{bakach2021two}        & Networks\\ \hline
\cite{caggiani2021green}    & Transportation Research Procedia & Elsevier\\ \hline
%\cite{liu2021hybrid}        & Transportation Research Part E: Logistics and Transportation Review\\ \hline
\cite{wang2021two}          & Applied Sciences & MDPI\\ \hline
\cite{zijlstra2021integrating} & Erasmus School of Economics: Econometrie Series & Erasmus University Rotterdam\\ \hline
%\cite{alewijnse2021minimising} & IFIP International Conference on Advances in Production Management Systems  & Springer\\ \hline
\cite{akbay2022variable}    & Applied Sciences & MDPI\\ \hline
%\cite{alfandari2022tailored} & European Journal of Operational Research\\ \hline
%\cite{ghobadi2022fuzzy}& Sustainability \\ \hline
%\cite{liu2022physical}      & International Journal of Production Economics\\ \hline
\cite{akbay2023application} & European Conference on Evolutionary Computation in Combinatorial Optimization & Springer \\ \hline
\cite{wu2023branch}         & Complex \& Intelligent Systems & Springer\\ \hline
%\cite{yu2024collaborative}  & Annals of Operations Research\\ \hline
\end{tabular}}
\label{tab:publications-journals-variants}
\end{table}

\begin{figure}[ht!]
 \centering \includegraphics[width=0.7\textwidth]{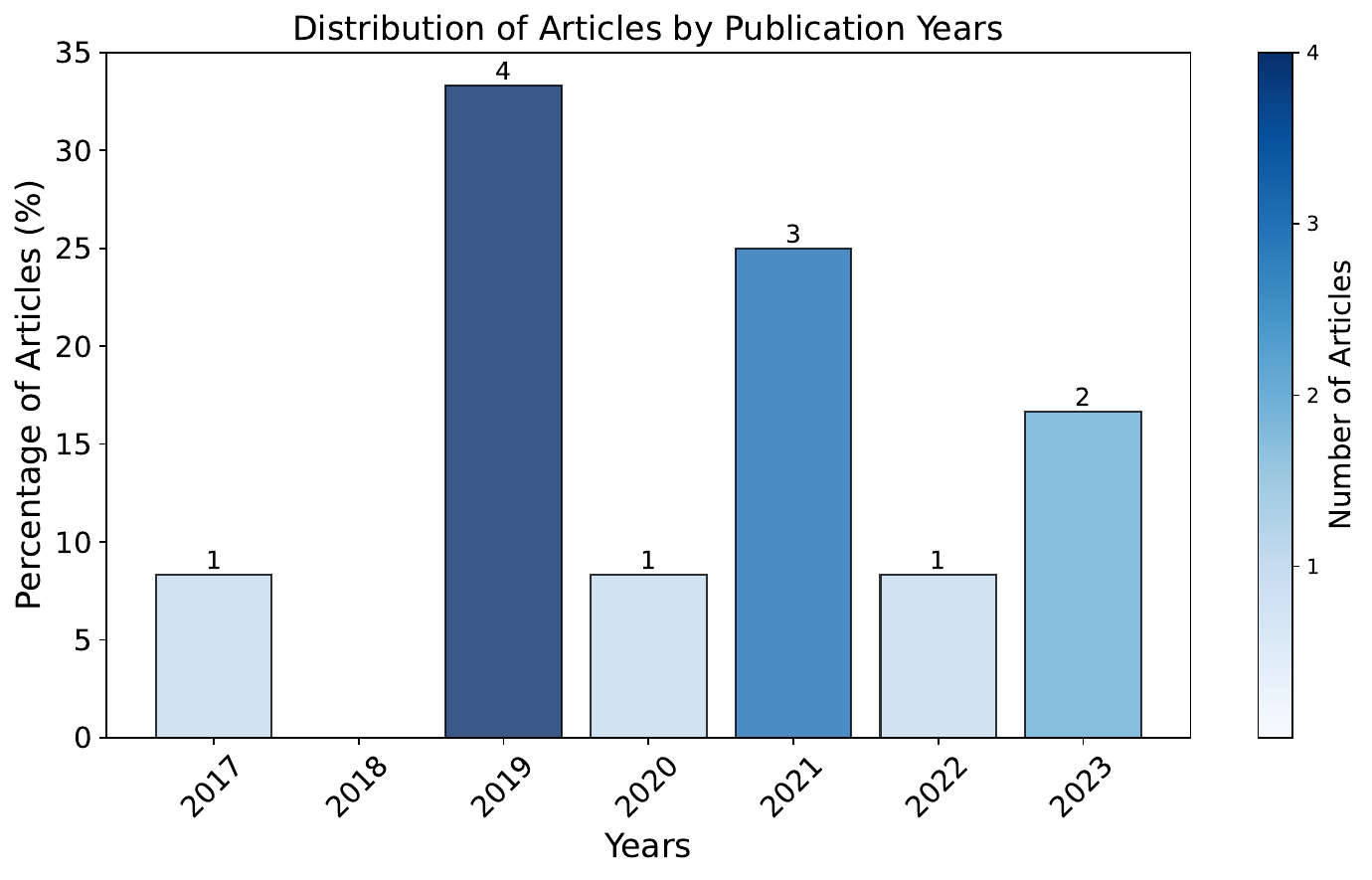}
 \caption{Distribution of articles by publication years and formulation types}
\label{fig:publication-year}
\end{figure}

As Fig. \ref{fig:publication-year} illustrates, research on the 2E-EVRP began in 2017 with the publication of the first article in this emerging field, i.e., the work by \cite{breunig2017two}. A two-year hiatus followed this until 2019, when there was a significant rise in activity, with four articles published, accounting for 33\% of the total studies. The momentum continued into 2020 with one additional article (8\%), while 2021 saw another surge with three publications (25\%), reflecting sustained interest in the topic. However, research output tapered off in subsequent years, with one article published in 2022 (8\%) and two in 2023 (17\%). 
All the studied articles     employed mixed-integer linear programming (MILP) formulations. MILP effectively captures the problem's complexity while offering computational efficiency and scalability, making it well-suited for addressing the challenges of 2E-EVRPs in large-scale logistics.

\begin{figure}[ht!]
 \centering \includegraphics[width=0.7\textwidth]{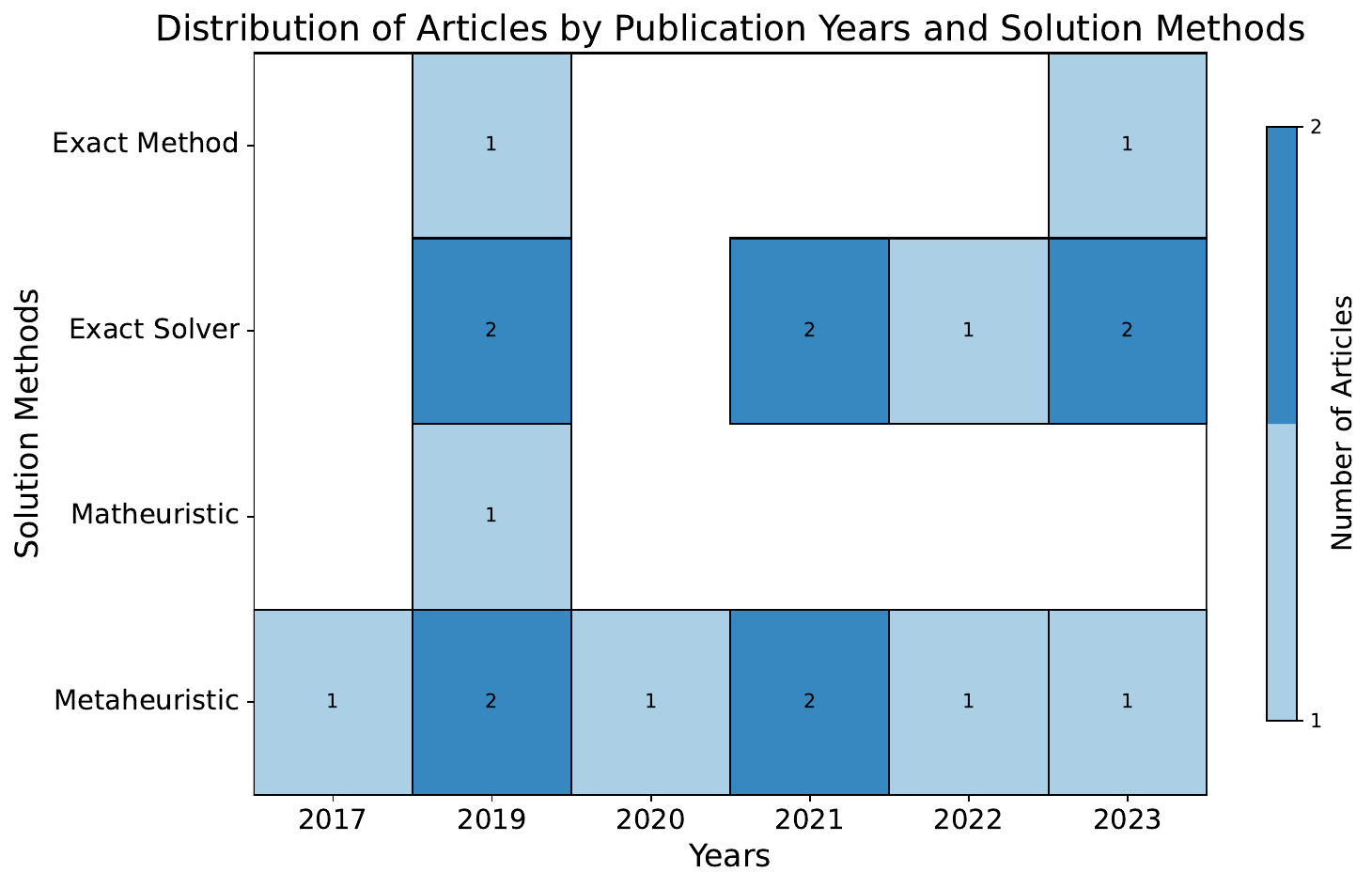}
 \caption{Distribution of articles by publication years and solution methods}
\label{fig:solution-method}
\end{figure}

Furthermore, Fig. \ref{fig:solution-method} displays the distribution of solution methods for addressing 2E-EVRPs published in recent years. Exact methods, such as branch-and-bound, emerged in 2019 but were used sparingly up to 2023. In contrast, exact solvers like Gurobi and CPLEX have gained moderate popularity using pre-built algorithms, particularly in alternate years since 2019. Matheuristics were only seen in 2019 and have not been widely used since then, while metaheuristics have consistently grown in prominence since 2017. Metaheuristics and exact solvers are closely aligned in usage, primarily due to a common strategy: using exact solvers for smaller instances and integrating them with metaheuristics to effectively address larger, more complex problems. All these solution methods are explained in detail in Section \ref{solution-methodology}.

Moreover, we analyze the key elements of each reviewed study, focusing on objective functions and constraints to reveal how these components reflect the underlying problem contexts and the research objectives pursued. This analysis enhances our understanding of how these models contribute to the field of 2E-EVRP. To facilitate content analysis, we categorize the objective functions and constraints into distinct groups and present them in Tables \ref{objective-functions-table}-\ref{constraints-table} with their definitions. Each category is assigned a notation to provide a structured framework for statistical interpretation. These objective functions and constraints will be employed to propose our classification framework (for more details on the objective functions addressed in each paper, please see Section \ref{obj-functions-section}).

\begin{table}[ht!]
\centering
\caption{Objective functions, descriptions, and notations}\label{objective-functions-table}
\resizebox{1\textwidth}{!}{
\begin{tabular}{|l|l|l|}
\hline
\textbf{Objective functions (minimization)} & \textbf{Description} & \textbf{Symbol} \\ \hline
Transportation (shipping) costs& Costs based on the distance traveled across \nth{1} and \nth{2} echelons.& 
$\mathcal{O}_1$ \\ \hline
Traveled distances& 
Total travel distance across \nth{1} and \nth{2} echelons. & 
$\mathcal{O}_2$ \\ \hline
Fixed vehicle usage costs&  Costs associated with \nth{1} and \nth{2} echelon vehicles, including initial investments.&  $\mathcal{O}_3$ \\ \hline
%Fixed depot utilization costs&  Costs associated with depots, including installment or locating costs.&  $\mathcal{O}_4$ \\ \hline
Fixed satellite utilization costs&  Costs associated with satellites, including installment or establishing costs.&  $\mathcal{O}_4$ \\ \hline
Handling costs& 
Costs associated with parcel processing, including load-dependent costs and satellite operations. & 
$\mathcal{O}_5$ \\ \hline
Energy/Battery consumption costs& 
Energy-related costs include battery swapping, charging, electricity usage, and operational energy expenses. & 
$\mathcal{O}_6$ \\ \hline
%Environmental &  Emissions from the vehicles and environmental factors from depot and satellite operations. &  $\mathcal{O}_8$ \\ \hline
%Delivery tardiness&  Weighted or maximum delivery tardiness, waiting times, or late deliveries across \nth{1} and \nth{2} echelons.&  $\mathcal{O}_7$ \\ \hline
Drivers' wages & Salaries paid to drivers working in the first or second echelons.&  $\mathcal{O}_{7}$ \\ \hline
\end{tabular}}
\end{table}

\begin{table}[ht!]
\centering
\caption{Constraints, descriptions, and notations}\label{constraints-table}
\resizebox{1\textwidth}{!}{
\begin{tabular}{|l|l|l|l|}
\hline
\textbf{Constraints} & \textbf{Subcategory} & \textbf{Description} & \textbf{Symbol} \\
\hline
%\textbf{Depot} & Capacity constraints & Limit the number of parcels stored at each depot and the number of open depots. & \(\mathcal{C}_1\) \\
%\hline
\textbf{\nth{1} Echelon Vehicles}& Load capacity constraints & Limit the number of parcels (cargo) that \nth{1}-echelon vehicle could carry & \(\mathcal{C}_1\) \\
\cline{2-4}
&Energy/Battery constraints & Limit the battery usage (traveled distance) by the \nth{1}-echelon vehicles  & \(\mathcal{C}_2\) \\
\cline{2-4}
&Route length constraints & Limit the traveling time of a \nth{1}-echelon vehicles on its route & \(\mathcal{C}_{3}\) \\
\cline{2-4}
&Fleet size constraints & Limit the number of used \nth{1}-echelon vehicles & \(\mathcal{C}_4\) \\
\hline
\textbf{Satellite} 
 &Vehicle storage constraints & Limit the number of \nth{2}-echelon vehicles located the satellite & \(\mathcal{C}_{5}\) \\
 \cline{2-4}
 &Capacity constraints & Limit the number of customers (demands) satisfied by the satellite& \(\mathcal{C}_{6}\) \\
\cline{2-4}
&Time window constraints & The satellite must be delivered within a specific period & \(\mathcal{C}_{7}\) \\
\hline
\textbf{\nth{2} Echelon Vehicles} &Load capacity constraints & Limit the number of parcels (cargo) that \nth{2}-echelon vehicle could carry & \(\mathcal{C}_{8}\) \\
\cline{2-4}
&Energy/Battery constraints &  Limit the battery usage (traveled distance) by the \nth{2}-echelon vehicles  & \(\mathcal{C}_{9}\) \\
\cline{2-4}
&Route length constraints &Limit the traveling time of a \nth{2}-echelon vehicles on its route  & \(\mathcal{C}_{10}\) \\
\cline{2-4}
&Fleet size constraints & Limit the number of used \nth{2}-echelon vehicles & \(\mathcal{C}_{11}\) \\
\hline
\textbf{Customer} &Time window constraints & The customer must be served within a specific period & \(\mathcal{C}_{12}\) \\
\hline
\end{tabular}}
\end{table}

\begin{figure}[ht!]
 \centering \includegraphics[width=0.7\textwidth]{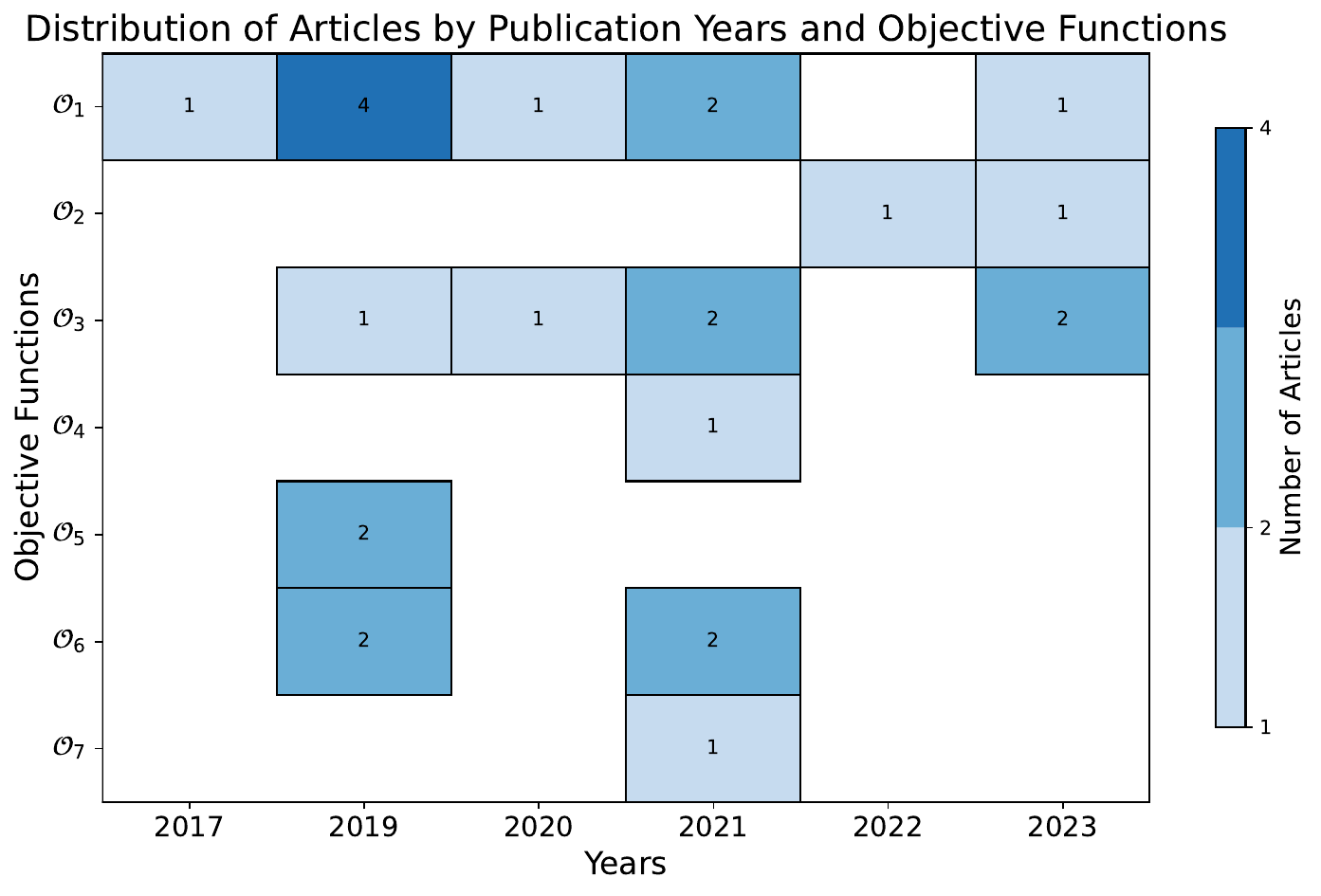}
 \caption{Distribution of articles by publication years and objective functions}
\label{fig:obj-funcs}
\end{figure}

Additionally, Fig.~\ref{fig:obj-funcs} provides a comprehensive view of the evolving priorities of objective functions in 2E-EVRP, revealing key shifts in focus over time. The dominance of transportation (shipping) costs (\(\mathcal{O}_1\)) across studies highlights the consistent emphasis on cost efficiency, a fundamental driver in logistics optimization. However, the rising attention to objectives like fixed vehicle usage costs (\(\mathcal{O}_3\)) and traveled distances (\(\mathcal{O}_2\)) suggests a growing awareness of operational inefficiencies and fleet management challenges as critical components of cost reduction strategies. Interestingly, sustainability-oriented objectives, such as energy and battery consumption costs (\(\mathcal{O}_6\)), have gained prominence in specific years, reflecting an alignment with global trends toward greener logistics practices. The concentration of studies in 2019 and 2021 suggests a period of heightened diversification, where research broadened to incorporate a broader range of objectives beyond traditional cost considerations. Less frequently studied objectives, such as fixed satellite costs (\(\mathcal{O}_4\)) and drivers’ wages (\(\mathcal{O}_7\)), highlight potential research gaps in addressing infrastructure and labor-related challenges. These areas may present opportunities for future exploration, particularly as the logistics industry evolves to meet increasing demands for resilience and adaptability. Overall, the literature demonstrates a dynamic progression, transitioning from a narrow focus on core cost objectives to a more nuanced understanding of trade-offs between economic, operational, and environmental factors. This trend indicates a maturing field that increasingly recognizes the complexity and interconnectedness of modern logistics systems.

\begin{figure}[ht!]
 \centering \includegraphics[width=0.7\textwidth]{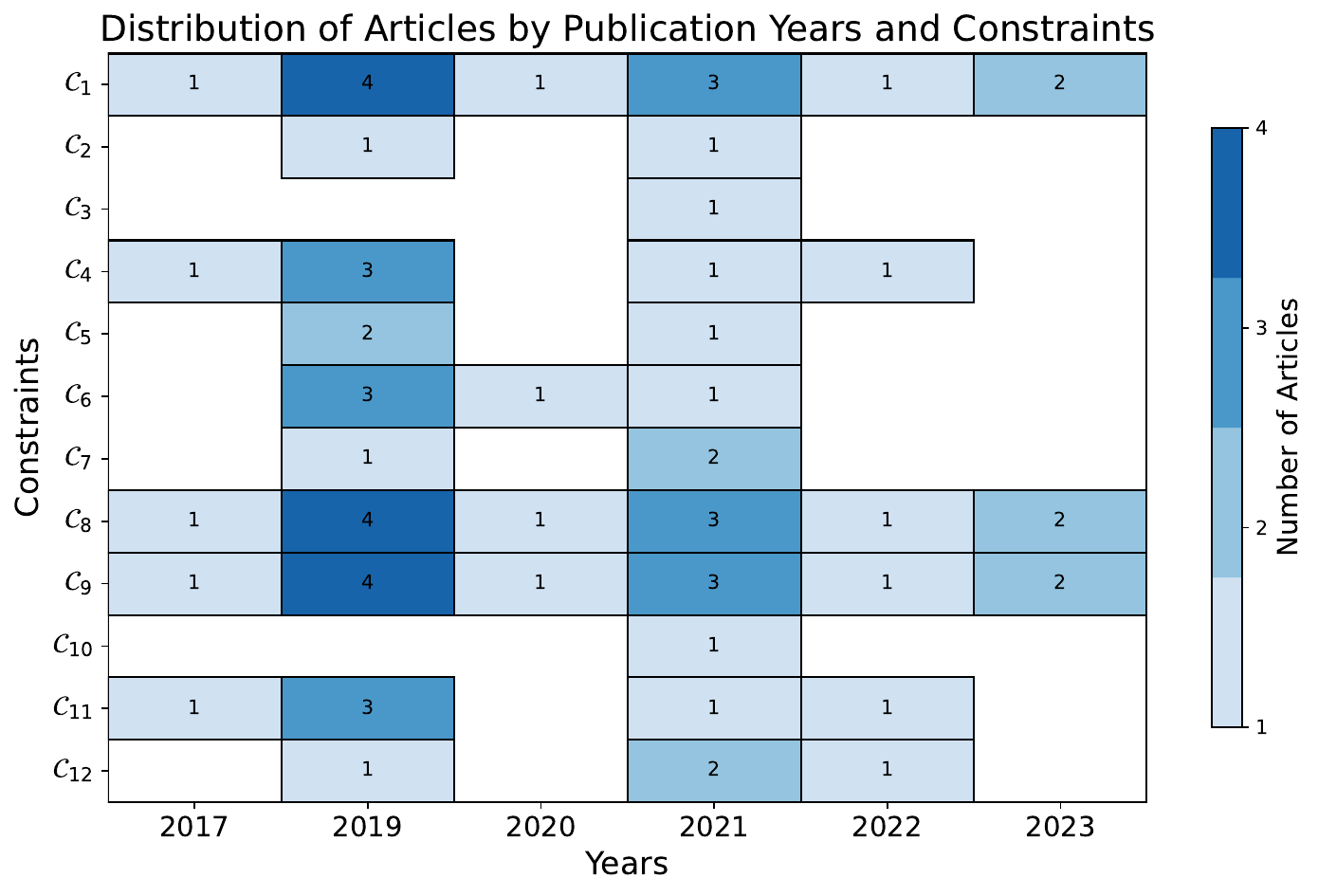}
\caption{Distribution of articles by publication years and constraints}
\label{fig:cons}
\end{figure}

Finally, Fig~\ref{fig:cons} illustrates the evolution of constraints in the 2E-EVRP problem from 2017 to 2023, reflecting increasing complexity and refinement in the model over time. In 2017 and 2022, the focus was primarily on operational constraints for the \nth{1} and \nth{2} echelon vehicles, with both years emphasizing load capacity (\(\mathcal{C}_1\), \(\mathcal{C}_8\)) and fleet size (\(\mathcal{C}_4\), \(\mathcal{C}_{11}\)) constraints. The \nth{2} echelon vehicles also had energy/battery consumption constraints (\(\mathcal{C}_9\)). The key difference in 2022 was the addition of customer time window constraints (\(\mathcal{C}_{12}\)), indicating a slight shift towards incorporating customer-specific needs while maintaining a similar foundational focus as in 2017. In 2019 and 2021, the model became more comprehensive. Both years saw the introduction of additional constraints, including energy/battery consumption (\(\mathcal{C}_2\)) for the \nth{1} echelon vehicles, as well as satellite-specific constraints like vehicle storage (\(\mathcal{C}_5\)), capacity (\(\mathcal{C}_6\)), and time windows (\(\mathcal{C}_7\)). A notable difference in 2021 was the inclusion of route length constraints (\(\mathcal{C}_3\), \(\mathcal{C}_{10}\)) for both echelons, reflecting an increased focus on optimizing operational efficiency and route planning. In contrast, 2020 and 2023 saw a return to a more streamlined set of constraints. Both years focused primarily on load capacity (\(\mathcal{C}_1\), \(\mathcal{C}_8\)) of both echelon vehicles and energy/battery consumption (\(\mathcal{C}_9\)) for the \nth{2} echelon vehicles, with 2020 introducing satellite capacity constraints (\(\mathcal{C}_6\)). Overall, the evolution of constraints reflects a trend of gradual complexity, starting with foundational operational constraints in 2017, expanding to include various operational and customer-specific factors in 2019 and 2021, and then narrowing back to core constraints in 2020 and 2023 for refined optimization.

\begin{figure}[ht!]
    \centering
\includegraphics[width=0.9\linewidth]{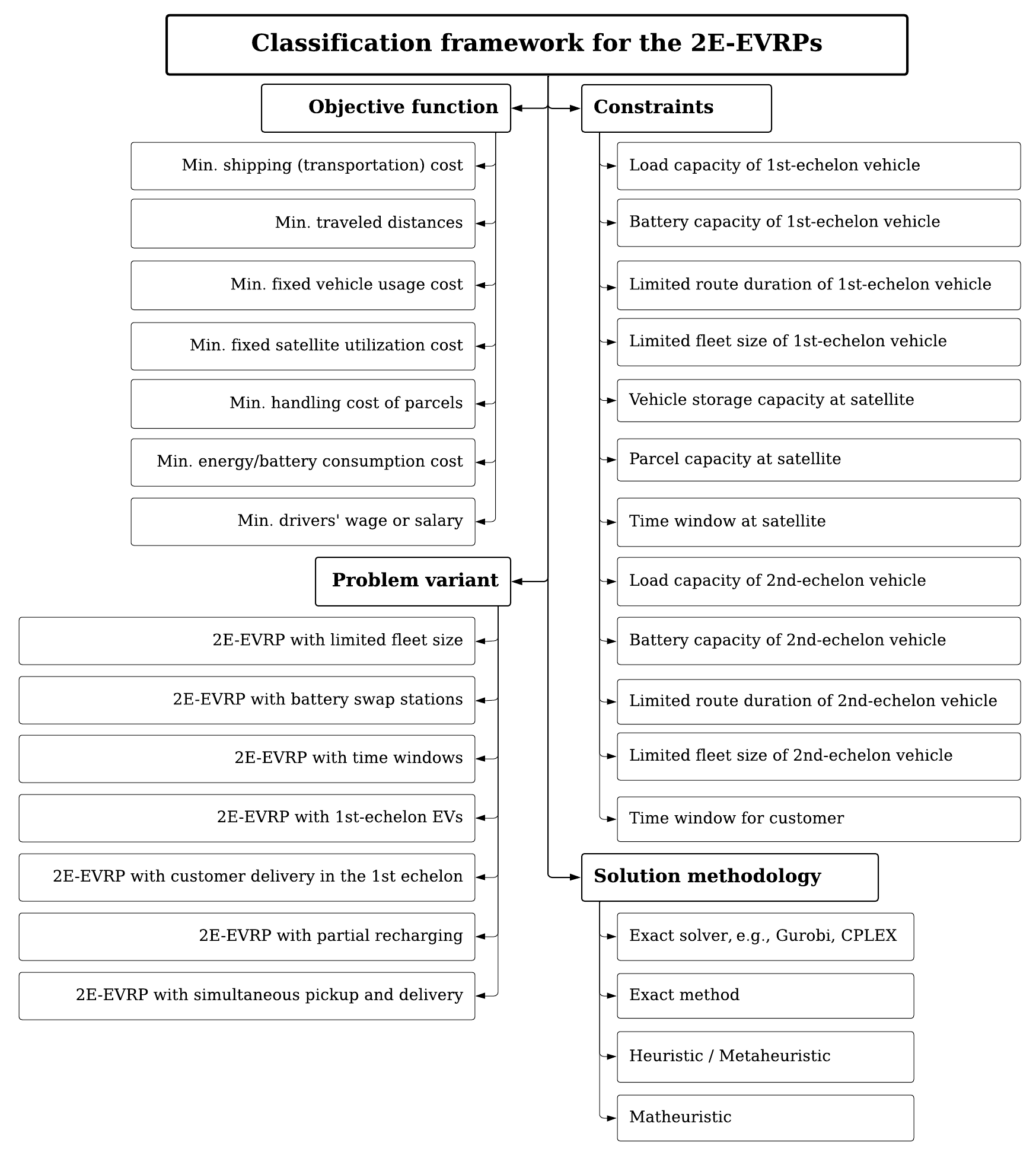}
    \caption{Proposed classification scheme for 2E-EVRPs in the literature}
    \label{fig:classifcation-scheme}
\end{figure}

In the following, we present the classification scheme for 2E-EVRP-related papers, as given by Fig. \ref{fig:classifcation-scheme}. This figure shows that each paper could be categorized into four main items such as problem variants, objective functions, constraints, and solution methodology. Each main item is also classified into several sub-items, which were explained previously. Notably, the problem variants are described in Section \ref{problem-type-variant}. Given this classification, Table \ref{tab:content-analysis} summarizes the twelve papers on 2E-EVRPs. As immediately observed, most papers have studied the 2E-EVRP with a limited fleet size of vehicles, minimizing the vehicle's transportation (shipping cost), limited load and battery capacity for vehicles while developing mathematical models (solved by exact solvers), and metaheuristic-based algorithms. These papers are explained in more detail in the following Sections \ref{problem-type-variant}-\ref{solution-methodology}.  

\begin{table}[ht!]
\centering
\footnotesize
\caption{Comparison of 2E-EVRPs in the literature according to the proposed classification framework}
\label{tab:content-analysis}
\begin{adjustbox}{max width=\textwidth}
\begin{tabular}{l l l l l }
\hline
\textbf{Reference} & \textbf{Problem variant} & \textbf{Objective function} & \textbf{Constraints} & \textbf{Solution methodology} \\ 
\hline
\cite{breunig2017two} & 2E-EVRP with limited fleet size & $\mathcal{O}_1$  &$\mathcal{C}_1\mathcal{C}_4\mathcal{C}_8\mathcal{C}_{9}\mathcal{C}_{11}$ & Metaheuristic  \\
\cite{agardi2019two} & 2E-EVRP with limited fleet size &$\mathcal{O}_1$  & $\mathcal{C}_1\mathcal{C}_4\mathcal{C}_6\mathcal{C}_8\mathcal{C}_{9}\mathcal{C}_{11}$  & Metaheuristic \\
\cite{breunig2019electric} & 2E-EVRP with limited fleet size &$\mathcal{O}_1$ & $\mathcal{C}_1\mathcal{C}_4\mathcal{C}_5\mathcal{C}_{6}\mathcal{C}_{8}\mathcal{C}_{9}\mathcal{C}_{11}$ & Exact method / Metaheuristic  \\
\cite{jie2019two} & 2E-EVRP-BSS with limited fleet size and \nth{1}-echelon EVs &$\mathcal{O}_1 \mathcal{O}_5\mathcal{O}_6$  & $\mathcal{C}_1\mathcal{C}_2\mathcal{C}_4\mathcal{C}_{5}\mathcal{C}_{6}\mathcal{C}_{8}\mathcal{C}_{9}\mathcal{C}_{11}$ & Exact solver / Matheuristic \\
\cite{wang2019two} & 2E-EVRPTW-BSS & $\mathcal{O}_1\mathcal{O}_3\mathcal{O}_5\mathcal{O}_6$ &$\mathcal{C}_1\mathcal{C}_7 \mathcal{C}_8\mathcal{C}_9\mathcal{C}_{12}$& Exact solver \\
\cite{affi2020general} & 2E-EVRP  &$\mathcal{O}_1\mathcal{O}_3$  &$\mathcal{C}_1\mathcal{C}_6 \mathcal{C}_8\mathcal{C}_9$  & Metaheuristic  \\
\cite{caggiani2021green} & 2E-EVRPTW-PR with \nth{1}-echelon EVs and customer delivery &$\mathcal{O}_3\mathcal{O}_4\mathcal{O}_6\mathcal{O}_7$ &$\mathcal{C}_1\mathcal{C}_2 \mathcal{C}_3\mathcal{C}_{7}\mathcal{C}_{8}\mathcal{C}_{9}\mathcal{C}_{10}\mathcal{C}_{12}$ & Exact solver \\
\cite{wang2021two} & 2E-EVRPTW-BSS   &$\mathcal{O}_1\mathcal{O}_3\mathcal{O}_6$ &$\mathcal{C}_1\mathcal{C}_7 \mathcal{C}_8 \mathcal{C}_9\mathcal{C}_{12}$& Exact solver / Metaheuristic  \\
\cite{zijlstra2021integrating} & 2E-EVRP with limited fleet size   & $\mathcal{O}_1$ & $\mathcal{C}_1\mathcal{C}_4\mathcal{C}_5\mathcal{C}_{6}\mathcal{C}_{8}\mathcal{C}_{9}\mathcal{C}_{11}$ & Metaheuristic \\
\cite{akbay2022variable} & 2E-EVRP-TW with limited fleet size & $\mathcal{O}_2$ & $\mathcal{C}_1\mathcal{C}_4\mathcal{C}_8\mathcal{C}_{9}\mathcal{C}_{11}\mathcal{C}_{12}$& Exact solver / Metaheuristic \\
\cite{akbay2023application} & 2E-EVRPSPD &$\mathcal{O}_2\mathcal{O}_3$  &  $\mathcal{C}_1\mathcal{C}_8\mathcal{C}_9$ & Exact solver / Metaheuristic \\
\cite{wu2023branch} & 2E-EVRP  & $\mathcal{O}_1\mathcal{O}_3$  &$\mathcal{C}_1\mathcal{C}_8\mathcal{C}_{9}$  & Exact solver / Exact method \\
\hline 
\end{tabular}
\end{adjustbox}
\end{table}

\section{Problem types \& variants}\label{problem-type-variant}
In this section, the various types of 2E-EVRP and its variants introduced in the literature are investigated in more detail. First, Section \ref{probelm-statement} describes the basic 2E-EVRP, including its assumptions and a generic mathematical model (in Section \ref{math-formulation}) inspired by the related works in the literature. Second, Section \ref{variants} thoroughly explains the 2E-EVRP variant and related papers while presenting their features and associated mathematical formulations.

\subsection{Problem statement}\label{probelm-statement}
\begin{figure}[ht!]
    \centering
    \fbox{\includegraphics[width=0.7\linewidth]{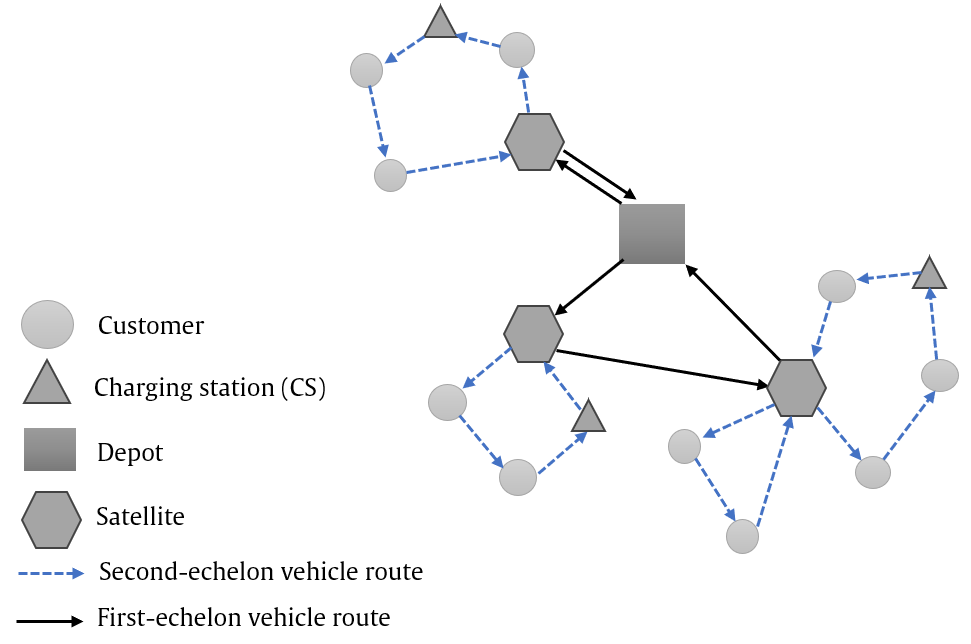}}
    \caption{Schematic example for a basic 2E-EVRP}
    \label{fig:basic-2e-evrp}
\end{figure}

Basic 2E-EVRP \citep{breunig2017two,breunig2019electric,affi2020general,zijlstra2021integrating,wu2023branch} has two levels: (i) between the central depot and satellite facilities and (ii) between satellites and customers, as represented by Fig. \ref{fig:basic-2e-evrp}. There is one main depot where the \nth{1}-echelon vehicles start their trips to deliver the parcels to the selected satellites. A satellite could be visited by multiple \nth{1}-echelon vehicles. When a satellite is selected, and the parcels are unloaded by \nth{1}-echelon vehicles, \nth{2}-echelon vehicles (e.g., EVs) are loaded with these parcels and leave the depot to deliver the goods to the end customers. The \nth{2}-echelon vehicles may be recharged at CSs if needed. However, precisely one \nth{2}-echelon vehicle must visit a customer place. Also, there is no opening cost for satellites, but a handling cost per parcel is assumed for each satellite. The goal is to deliver the parcels from the main depot to the satellite (in the first echelon), where the \nth{2}-echelon vehicles are loaded with these parcels and delivered to the end customers (in the second echelon, or last mile).    

\subsubsection{Mathematical formulation}\label{math-formulation}
The basic 2E-EVRP could be represented as a graph $G=(V,A)$ where $V$ is the set of all nodes; $V=D\cup C\cup F\cup S$. $D=\{0\}$ is the central depot, where a fleet of \nth{1}-echelon vehicles denoted by $K^1$ are located and ready to start their trips. Each \nth{1}-echelon vehicle has a maximum load capacity $Q^1$ and fixed usage cost $\mathcal{F}^1$. $C$ is the set of customer nodes, each associated with a delivery demand $q_i$, and they must be satisfied by a \nth{2}-echelon vehicle. $F$ represents the set of CS (charging stations), where the \nth{2}-echelon vehicles may visit to recharge their battery to a full capacity $B^2$. $S$ is the set of available satellite locations, and a dummy node for each satellite $s\in S$ is defined as $s'\in S'$, where $S'$ is the set of dummy satellites. each satellite $s\in S$ could launch at most $|K^2_s|$ \nth{2}-echelon vehicles, each associated with a maximum load capacity $Q^2$, battery capacity $B^2$, battery consumption rate per unit distance $h^2$, and a fixed usage cost $\mathcal{F}^2$. Also, $A$ refers to the set of directed arcs in the graph, where $A=A^1\cup A^2$. The set of arcs in the first echelon $A^1$ is defined as the directed edges between the central depot and satellites; $A^1=\{(i,j)|i,j\in D\cup S,i\neq j\}$. Similarly, $A^2$ is defined as the arcs in the second echelon connecting the satellites to customer and CS nodes; $A^2=\{(i,j)|i,j\in C\cup S\cup F\cup S',i\neq j\}$. In addition, a handling cost per parcel $h_s$ is assumed for each satellite $s\in S$, and the total demands that a satellite $s$ could satisfy is limited by a capacity parameter $Q'_s$.

Eight decision variables are defined for the basic 2E-EVRP. $x_{ijk}$ is one if the \nth{1}-echelon vehicle travels on the arc $(i,j)\in A^1$; zero otherwise. A similar variable is defined for \nth{2}-echelon vehicles: $y_{ijks}$ is one if the \nth{2}-echelon vehicle $k\in K^2$ launching from satellite $s\in S$ travels on the arc $(i,j)\in A^2$; zero otherwise. The selected satellite $s\in S$ is indicated by $\lambda_s$, determining the total number of parcels (demands) delivered to that satellite. The total number of the parcels delivered to satellite $s\in S$ via \nth{1}-echelon vehicle $k\in K^1$ is indicated by variable $w_{sk}$. Two decision variables are defined to track the cargo on the vehicles: $f_{ijk}$ denotes the remaining cargo load on the \nth{1}-echelon vehicle $k\in K^1$ while traveling on arc $(i,j)\in A^1$; similarly, $f'_{ijks}$ represents the remaining cargo load on the \nth{2}-echelon vehicle $k\in K^2$ launched from satellite $s\in S$ when traveling on arc $(i,j)\in A^2$. In addition, two variables are defined to monitor the vehicles' battery levels: $b^{2+}_{isk}$ denotes the battery level of \nth{2}-echelon vehicle $k\in K^2$ launched from satellite $s\in S$ when entering the node $i$ in the second echelon; similarly, $b^{2-}_{isk}$ indicates the battery level of \nth{2}-echelon vehicle $k\in K^2$ launched from satellite $s\in S$ when leaving the node $i$ in the second echelon. The model's objective function is to minimize the total delivery costs, including traveled distance by vehicles, usage cost of vehicles, and handling cost in the satellites while satisfying several operational constraints. To ease readership, Table \ref{notations} summarizes all notations, parameters, and decision variables of the basic 2E-EVRP model. Given these notations, the mixed-integer linear programming model, inspired by the arc-based formulation developed by \cite{jie2019two}, is presented as follows.

\begin{table}[ht!]
\centering
\caption{Notations, parameters, and decision variables of the basic 2E-EVRP model}\label{notations}
\begin{adjustbox}{max width=\textwidth}
\begin{tabular}{ll}
    \hline
      \textbf{Set}   & \textbf{Definition}  \\
      \hline
       $V$  & All nodes; $V=D\cup C\cup S\cup F$\\
       $D$ & Depot; $D=\{0\}$ \\
       $C$ & Customer nodes\\
       $S$ & Satellite nodes \\
       $S'$ & Dummy satellites (for each satellite node $s\in S$, there is a corresponding dummy satellite $s'\in S'$)\\
       $F$ & Recharging stations \\
       $K^1$ & \nth{1}-echelon vehicles located at the depot\\
       $K^2$ & \nth{2}-echelon vehicles located at the satellites; $K^2=\sum_{s\in S}K^2_s$\\
       $K^2_s$ & \nth{2}-echelon vehicles located at satellite $s\in S$\\
       $A^1$ & Set of arcs in the first echelon; $A^1=\{(i,j)|i,j\in D\cup S
       ,i\neq j\}$\\
       $A^2$ & Set of arcs in the second echelon; $A^2=\{(i,j)|i,j\in C\cup S\cup F\cup S',i\neq j\}$\\
       \hline
\textbf{Parameters} & \\
\hline
%$|K^2|$ & Number of all \nth{2}-echelon vehciles\\
$q_i$ & Demand of customer $i\in C$\\
$Q^1$ & \nth{1}-echelon vehicle load capacity\\
$\mathcal{F}^1$ & Fixed usage cost of \nth{1}-echelon vehicle\\
$Q^2$ & \nth{2}-echelon vehicle load capacity\\
$B^2$ & \nth{2}-echelon vehicle battery capacity\\
$h^2$ & \nth{2}-echelon vehicle battery consumption rate\\
$\mathcal{F}^2$ & Fixed usage cost of \nth{2}-echelon vehicle\\
%$Q_s$ & Satellite $s\in S$ accommodation capacity\\
$d_{ij}$ & Distance between node $i,j\in V$\\
$h_s$ & The unit cost of handling each parcel in the satellite $s\in S$\\
$Q'_s$ & Capacity of satellite $s\in S$, indicating the number of demands it can satisfy\\
\hline
\textbf{Decision variables} & \\
\hline
$x_{ijk}$ & One if the arc $(i,j)\in A^1$ is traveled by the \nth{1}-echelon vehicle $k\in K^1$; 0 otherwise.\\
$y_{ijks}$ & One if the arc $(i,j)\in A^2$ is traveled by the \nth{2}-echelon vehicle $k\in K^2$ dispatched from satellite $s\in S$\\
$\lambda_s$ & Total demands fulfilled by the satellite $s\in S$\\
$w_{sk}$ & Total demands shipped to the satellite $s\in S$ by \nth{1}-echelon vehicle $k\in K^1$\\
$f_{ijk}$ & Remaining load on \nth{1}-echelon vehicle $k\in K^1$ departing from node $i$ to node $j$ ($(i,j)\in A^1$)\\
$f'_{ijks}$ & Remaining load on \nth{2}-echelon vehicle $k\in K^2$ dispatched from satellite $s\in S$ while traveling from node $i$ to node $j$ ($(i,j)\in A^2$)\\
$b^{2+}_{isk}$ & Remaining battery of \nth{2}-echelon vehicle $k\in K^2$ dispatched from satellite $s\in S$ when arriving at node $i\in C\cup S\cup F$\\
$b^{2-}_{isk}$ & Remaining battery of \nth{2}-echelon vehicle $k\in K^2$ dispatched from satellite $s\in S$ upon departure from node $i\in C\cup S\cup F$\\
\hline
\end{tabular}
\end{adjustbox}
\end{table}

\begin{equation}\label{eq1}
    Min. \sum_{(i,j)\in A^1}\sum_{k\in K^1}d_{ij}x_{ijk}
    +\sum_{(i,j)\in A^2}\sum_{k\in K^2}\sum_{s\in S}d_{ij}y_{ijks}+\sum_{j\in S}\sum_{k\in K^1}\mathcal{F}^1x_{0jk}+\sum_{s\in S}\sum_{j\in C\cup F}\sum_{k\in K^2}\mathcal{F}^2y_{sjks} + \sum_{s\in S}h_s\lambda_s
\end{equation}

s.t.,

\begin{equation}\label{eq2}
 \sum_{j\in D\cup S: j\neq i}x_{ijk}=\sum_{j\in D\cup S: j\neq i}x_{jik}, \quad \quad \forall i\in D\cup S, \forall k\in K^1   
\end{equation}

\begin{equation}\label{eq3}
\sum_{j\in D\cup S: j\neq s}x_{sjk}\le 1,\quad \quad \forall s\in S, \forall k\in K^1 
\end{equation}

\begin{equation}\label{eq4}
w_{sk}=\sum_{j\in D\cup S: j\neq i}f_{jsk}-\sum_{j\in D\cup S: j\neq i}f_{sjk}, \quad \quad \forall s\in S, \forall k\in K^1   
\end{equation}

\begin{equation}\label{eq5}
f_{ijk}\le Q^1x_{ijk}\quad \quad \forall i,j\in D\cup S, i\neq j, \forall k\in K^1    
\end{equation}

\begin{equation}\label{eq6}
\sum_{s\in S}w_{sk}\le Q^1, \quad \quad \forall k\in K^1
\end{equation}

\begin{equation}\label{eq7}
\sum_{k\in K^1}w_{sk}=\lambda_s, \quad \quad \forall s\in S     
\end{equation}

\begin{equation}\label{eq8}
\sum_{k\in K^2}\sum_{s\in S}\sum_{j\in C\cup S\cup F:j\neq i}y_{ijks}=1, \quad \quad \forall i\in C    
\end{equation}

\begin{equation}\label{eq9}
\sum_{j\in C\cup S\cup F:j\neq i}y_{ijks}=\sum_{j\in C\cup S\cup F:j\neq i}y_{jiks}, \quad \quad \forall i\in C, \forall s\in S, \forall k\in K^2    
\end{equation}

\begin{equation}\label{eq10}
\sum_{\hat{s}\in S:\hat{s}\neq s}(\sum_{j\in C\cup S\cup F:j\neq s}y_{sjk\hat{s}}+\sum_{i\in C\cup S\cup F:i\neq s'}y_{is'k\hat{s}})=0,\quad \quad \forall s\in S, \forall s' \in S', \forall k\in K^2    
\end{equation}

\begin{equation}\label{eq11}
 \sum_{s\in S}\sum_{k\in K^2}y_{ss'ks}=0   
\end{equation}

\begin{equation}\label{eq12}\sum_{j\in C\cup S\cup F:j\neq s}\sum_{k\in K^2}y_{sjks}\le |K^2_s|,\quad \quad \forall s\in S\end{equation}

%\begin{equation}\label{eq14}\sum_{s\in S}\sum_{j\in C\cup S:j\neq s}\sum_{r\in R}y_{sjrs}\le |R|\end{equation}

\begin{equation}\label{eq13}
\sum_{j\in C\cup S\cup F:j\neq i}\sum_{s\in S}\sum_{k\in K^2}f'_{ijks}=\sum_{j\in C\cup S\cup F:j\neq i}\sum_{s\in S}\sum_{k\in K^2}f'_{jiks}-q_i,\quad \quad \forall i\in C
\end{equation}

\begin{equation}\label{eq14}
f'_{ijks}\le Q^2y_{ijks},\quad \quad \forall i,j \in C\cup S\cup F, i\neq j, \forall s\in S, \forall k\in K^2
\end{equation}

\begin{equation}\label{eq15}
\sum_{i,j\in C\cup S\cup F:i\neq j}\sum_{k\in K^2}q_iy_{ijks}=\sum_{k\in K}w_{sk}, \quad \quad \forall s\in S    
\end{equation}

\begin{equation}\label{eq15.1}
\sum_{i,j\in C\cup S\cup F:i\neq j}\sum_{k\in K^2}q_iy_{ijks}\le Q'_s, \quad \quad \forall s\in S    
\end{equation}

\begin{equation}\label{eq16}
    b^{2-}_{isk}=B^2, \quad \quad \forall s\in S, \forall i\in S\cup F, \forall k\in K^2
\end{equation}

\begin{equation}\label{eq17}
    b^{2-}_{isk}=b^{2+}_{isk}, \quad \quad \forall s\in S, \forall i\in C, \forall k\in K^2
\end{equation}

\begin{equation}\label{eq18}
    b^{2+}_{jsk}\le b^{2-}_{isk}-h^2d_{ij}y_{ijks}+B^2(1-y_{ijks}), \quad \quad \forall (i,j)\in A^2, s\in S, k\in K^2
\end{equation}

\begin{equation}\label{eq19}
x_{ijk}\in \{0,1\},f_{ijk}\ge 0,\quad \quad \forall (i,j)\in A^1, \forall k\in K^1
\end{equation}

\begin{equation}\label{eq20}
y_{ijks}\in \{0,1\},f'_{ijks}\ge 0,\quad \quad \forall (i,j)\in A^2, \forall k\in K^2
\end{equation}

\begin{equation}\label{eq21}
\lambda_s,w_{sk}\ge 0,\quad \quad \forall s\in S, \forall k\in K^1
\end{equation}

\begin{equation}\label{eq22}
    b^{2-}_{isk},b^{2+}_{isk}\ge 0, \quad \quad \forall s\in S, \forall i\in S\cup C\cup F\cup S', \forall k\in K^2
\end{equation}

This model's objective function, as shown in equation \eqref{eq1}, contains five terms. The first term minimizes the distances traveled by the \nth{1}-echelon vehicles. Similarly, the second term minimizes the total traveled distances by the \nth{2}-echelon vehicles (e.g., EVs). The third and fourth terms minimize the total number of used \nth{1}- and \nth{2}-echelon vehicles (or the fixed usage cost of using these vehicles). The fifth term minimizes the total handling cost of loading/unloading the parcels in the satellites. Constraint \eqref{eq2} is known as vehicle flow balance constraint, ensuring that the number of vehicles arriving at and departing from each node in the first echelon (i.e., main depot and satellites) must be equal. Constraint \eqref{eq3} enforces that each \nth{1}-echelon vehicle could not be used more than once. Constraint \eqref{eq4} holds the balance between the cargo level of the \nth{1}-echelon vehicle before and after delivering to a satellite. It also finds the number of parcels delivered to every satellite in the first echelon. Constraint \eqref{eq5} guarantees that the load capacity of \nth{1}-echelon vehicles must not be violated. Constraint \eqref{eq6} enforces that a \nth{1}-echelon vehicle could not deliver more than its load capacity to the satellite. Note that constraints \eqref{eq4}-\eqref{eq6} track the cargo level of the \nth{1}-echelon vehicles. Constraint \eqref{eq7} calculates the total number of parcels delivered to each satellite by the \nth{1}-echelon vehicles. Constraints \eqref{eq8}-\eqref{eq12} are vehicle flow constraints in the second echelon. In particular, constraint \eqref{eq8} ensures that each custom node must be satisfied or served by precisely one \nth{2}-echelon EV. Constraint \eqref{eq9} determines that the number of arcs (vehicles) arriving at a node in the second echelon must equal the number of arcs (vehicles) departing from that node. Constraint \eqref{eq10} eliminates the arcs connecting a pair of satellites in the second echelon route. It also ensures that each \nth{2}-echelon EV must start and finish at the same satellite. Constraint \eqref{eq11} eliminates the arc from a satellite to its dummy. Constraint \eqref{eq12} indicates that the number of EVs dispatched from a satellite must not exceed the maximum available number of EVs located at that satellite. Constraints \eqref{eq13}-\eqref{eq15} tack the cargo level of \nth{2}-echelon EVs while visiting the nodes. Particularly, constraint \eqref{eq13} indicates that the difference between the cargo level of a \nth{2}-echelon EV before and after serving a customer node must equal the demand quantity at that customer node. Constraint \eqref{eq14} ensures that the load capacity of \nth{2}-echelon EVs must not be violated. Constraint \eqref{eq15} ensures that the total demands a satellite satisfies must equal the number of parcels delivered to it via the \nth{1}-echelon vehicles. Constraint \ref{eq15.1} ensures that each satellite can not exceed its maximum capacity to fulfill the assigned demands. Constraints \eqref{eq16}-\eqref{eq18} track the battery or energy level of \nth{2} echelon EVs while traveling the arcs. Specifically, constraint \eqref{eq16} ensures that an EV starts its trip at the whole battery level. Constraint \eqref{eq17} indicates that the battery level of an EV is not changed before and after visiting a customer node. Constraint \eqref{eq18} prevents the battery capacity of a \nth{2}-echelon EV from being violated during its trip. Finally, constraints \eqref{eq19}-\eqref{eq22} denote the model's decision variables' type and domain.     

\subsection{Variants}\label{variants}
The basic 2E-EVRP model provides a foundation for optimizing delivery operations across two stages, typically from depots to satellite facilities and from satellites to customers. However, real-world systems are rarely straightforward; diverse factors such as time windows, heterogeneous fleets, environmental constraints, and EV-specific requirements like various recharging technologies influence real-world systems. The variants of 2E-EVRP in the literature are investigated in the following.

\subsubsection{2E-EVRP with first echelon EVs}\label{2e-evrp-first-evs}
\begin{figure}
\centering
\fbox{\includegraphics[width=0.7\linewidth]{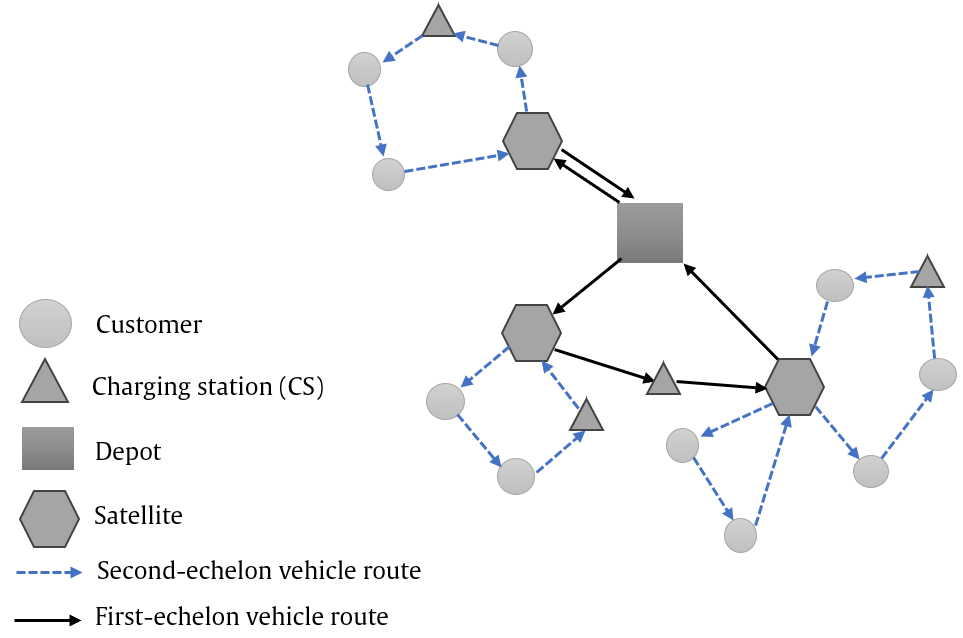}}
    \caption{Schematic example for 2E-EVRP with \nth{1}-echelon EVs}
    \label{fig:2e-evrp-evs}
\end{figure}

As mentioned earlier, EVs emerge as a transformative solution for enhancing sustainability and efficiency in last-mile delivery. Utilizing EVs, even as the first echelon delivery vehicles, could reduce greenhouse gas emissions, lowering noise pollution and minimizing dependency on fossil fuels. Including the EVs in the first echelon is represented by Fig. \ref{fig:2e-evrp-evs}, in which the EVs at both echelons need to be recharged at designated CSs. Their operational advantages, such as lower maintenance costs and navigating urban restrictions like low-emission zones, make them highly suitable for urban delivery networks. Large EVs with a high cargo capacity like Purolator e-truck \citep{e-truck-puro} and Ford e-transit \citep{ford.ca} can optimize delivery routes, reduce energy consumption, and improve service quality, positioning EVs as a cornerstone for sustainable last-mile logistics. To consider the EVs with limited battery capacity and required recharging infrastructure in the basic model of the 2E-EVRP, the arc set in the first echelon is modified by adding the recharging stations: $A^1=\{(i,j)|i,j\in D\cup S\cup F,i\neq j\}$. Also, the following constraints \eqref{eq23}-\eqref{eq26} could be added to include the battery capacity of EVs in the first echelon. $B^1$ and $h^1$ are the battery capacity and consumption rates of the \nth{1}-echelon EVs. Also, two decision variables $b^{1+}_{ik}$ and $b^{1-}_{ik}$ are defined as the remaining battery of the \nth{1}-echelon EV $k\in K^1$ while reaching and leaving the node $i\in D\cup S\cup F$, respectively. 

\begin{equation}\label{eq23}
    b^{1-}_{ik}=B^1, \quad \quad \forall  i\in D\cup F, \forall k\in K^1
\end{equation}

\begin{equation}\label{eq24}
    b^{1-}_{ik}=b^{1+}_{ik}, \quad \quad \forall i\in S, \forall k\in K^1
\end{equation}

\begin{equation}\label{eq25}
    b^{1+}_{jk}\le b^{1-}_{ik}-h^1d_{ij}x_{ijk}+B^1(1-x_{ijk}), \quad \quad \forall (i,j)\in A^1, k\in K^1
\end{equation}

\begin{equation}\label{eq26}
    b^{1-}_{ik},b^{1+}_{ik}\ge 0, \quad \quad \forall i\in D\cup S\cup F, \forall k\in K^1
\end{equation}

The constraint \eqref{eq23} ensures that the \nth{1}-echelon EV leaves the depot or CS with the whole battery level. Constraint \eqref{eq24} determines that the \nth{1}-echelon EV battery level is not changed after serving a customer node. Constraint \eqref{eq25} ensures that the battery capacity of the \nth{1}-echelon EV is not violated while traveling its route. Constraint \eqref{eq26} determines the domain of $b^{1-}_{ik}$ and $b^{1+}_{ik}$. Only two papers in the literature consider the EVs in the first and second echelons. In the first work, \citep{jie2019two} addressed the 2E-EVRP with battery swapping stations (2E-EVRP-BSS), where goods are delivered from a depot to customers via satellites. The system uses two types of EVs across two echelons. The battery swapping stations (BSSs) are located at fixed sites and have no service limitations, meaning they can serve both types of EVs. The first echelon consists of large EVs with long battery ranges, while the second uses smaller ones. The vehicle routing and battery swapping for both echelons are optimized simultaneously in the proposed model. An important feature is that loops are allowed in the BSS network, meaning an EV can visit a BSS multiple times. One small EV must visit each customer, but not all satellites must be used. The assumption is that the large and small EVs at the depot and satellites, respectively, start with fully charged batteries and can be recharged at BSSs during their journey, with a swapping cost incurred at each stop. Since the EVs have limited battery capacities, multiple satellite visits may be required to meet customer demands without exceeding the vehicles' capacities. 

In the second work, \cite{caggiani2021green} suggested a 2E-EVRP that incorporates time windows and partial recharging (2E-EVRPTW-PR), modeled as a mixed-integer program. In this model, the first echelon corresponds to an urban zone, while the second echelon comprises a restricted traffic area. The two echelons are connected by transmission points that link the e-vans (first echelon), delivering goods to the e-bikes (second echelon). These transmission points are dummy stations, allowing multiple visits to the CSs. At each transmission point, the quantity of available goods equals the total amount delivered to customers, ensuring there is no daily inventory cost at these points. The objective of the model is to minimize the travel costs, the initial investment in vehicles, driver salaries, and the transmission point costs across both echelons.

Both \citep{jie2019two} and \cite{caggiani2021green} aim to optimize the EVRP in a two-echelon framework but differ in their approach to battery management and charging strategies. \citep{jie2019two}  emphasize BSS as a key component for maintaining vehicle range, where multiple visits to BSSs can be necessary due to the limited battery capacities of smaller EVs. In contrast, \cite{caggiani2021green} introduces partial recharging and time windows, which allows for more flexible recharging strategies and accommodates different types of EVs across both echelons. Both models aim to minimize the total operational cost, but \citep{jie2019two} focus more on the handling and swapping costs at BSSs. At the same time, \cite{caggiani2021green} incorporate additional cost factors such as initial vehicle investments and driver salaries. The main difference lies in detail regarding battery management strategies, with \citep{jie2019two} proposing a swapping-based solution and \cite{caggiani2021green} using a partial recharging approach combined with time constraints.

\subsubsection{2E-EVRP with limited number of \nth{1}- and \nth{2}-echelon vehicles}
Two main approaches address the fleet size of the EVs used in the 2E-EVRP: minimizing the number of EVs used while imposing an upper limit on their quantity. The former approach is considered in the general mathematical model \eqref{eq1}-\eqref{eq22} for the basic 2E-EVRP, where the second and fourth terms of the objective function \eqref{eq1} are minimizing the EVs fixed usage costs (or the number of used EVs). The following papers addressed the 2E-EVRP with limited fleet size for \nth{1} and \nth{2} echelons: \cite{breunig2017two,agardi2019two,breunig2019electric,jie2019two,zijlstra2021integrating,akbay2022variable}. The following constraints \eqref{eq27}-\eqref{eq28} could be added to the basis model of the 2E-EVRP to set an upper limit on the number of vehicles used in both the first and second echelons. The number of vehicles used in the first and second echelons is capped at a maximum, denoted by $\mathcal{M}^1$ and $\mathcal{M}^2$, respectively. Using these constraints, the third and fourth terms of the objective function \eqref{eq1} could be removed from the model.  

\begin{equation}\label{eq27}
\sum_{j\in S}\sum_{k\in K^1}x_{0jk}\le \mathcal{M}^1
\end{equation}

\begin{equation}\label{eq28}
\sum_{s\in S}\sum_{j\in S\cup C\cup F:j\neq s}\sum_{k\in K^2}y_{sjks}\le \mathcal{M}^2
\end{equation}

\subsubsection{2E-EVRP with time windows}
Customer or satellite time windows are critical to a two-tier parcel delivery system, ensuring timely and reliable service while enhancing customer satisfaction. By adhering to specific time slots requested by customers, delivery services can meet expectations for convenience and flexibility, which are especially important in the era of e-commerce and rapid same-day delivery demands. Time windows also help reduce delivery failures, such as missed deliveries, by aligning drop-offs with customer availability. Additionally, respecting customer time windows fosters trust and loyalty, giving companies a competitive edge in a market driven by service quality and customer-centric logistics. In the literature of 2E-EVRP, four papers were found to consider time windows-related constraints in their problems, described as follows.

\cite{wang2019two} presented a 2E-EVRP with time windows and BSS (2E-EVRPTW-BSS). The system involves one central depot, satellite locations, customers, and BSSs. Deliveries in the first echelon can be split, but customer demands are known and indivisible, each with a specified (hard) time window for delivery. In the first echelon, vehicles transport demand from the depot to the satellites, from where they proceed to serve the customers. A delivery deadline is associated with each satellite, which must be satisfied in the first echelon (i.e., \nth{1}-echelon vehicles must deliver the parcels to a satellite before the deadline of that satellite). In the second echelon, EVs serve the customers at their homes and may visit BSSs to change the batteries before they are completely drained. Each satellite has a CS, which various EVs can visit multiple times. They assume that service times at the BSSs and customer nodes are negligible. The objective is to minimize the overall expenses, including costs associated with transportation,  BSSs (cost of swapping the batteries), satellite handling charges, vehicle usage costs, and penalties for delayed customer deliveries. 

Although \cite{wang2019two} claimed that they minimized the penalties for delayed deliveries, they minimized the cost of excess inventory at the satellites. To incorporate the customers' and satellites' time windows in the model, they defined two non-negative continuous decision variables $\tau_{ik}$, representing the \nth{1}-echelon vehicle's $k\in K^1$ arrival time at node $i\in D\cup S\cup S'$; and $\tau'_{ik}$, indicating the arrival time \nth{2}-echelon vehicle $k\in K^2$ to the node $i\in C\cup F\cup S\cup S'$. Also, for the depot $D=\{0\}$, a copy (dummy node) is created as $D'$. For each customer $i \in C$, the time windows' lower and upper bounds are given by the earliest and latest service times, denoted by $e_i$ and $l_i$, respectively. In addition, the delivery deadline for a satellite $s\in S$ is denoted by parameter $l_s$. To track the time of vehicle arrivals, the travel time of each \nth{1}- or \nth{2}-echelon vehicle between nodes $i$ and $j$ is represented by $t_{ij}$, which is obtained as $d_{ij}/v$, where $v$ is the speed of the vehicle. To generalize the definition of these parameters to include in our basic model \eqref{eq1}-\eqref{eq22}, the traveling time between nodes could be defined according to the vehicles in the first and second echelons: $t^1_{ij}=d_{ij}/v^1, \forall (i,j)\in A^1$ and $t^2_{ij}=d_{ij}/v^2, \forall (i,j)\in A^2$ where $v^1$ and $v^2$ are the \nth{1}- and \nth{2}-echelon vehicles' driving speed. In addition, the service time at customer and satellite nodes is denoted by $s_i$. The service time at customer nodes means the time for a \nth{2}-echelon vehicle to deliver a parcel to that customer; however, for a satellite, this could be seen as a collected duration of parcel unloading time in that satellite. Also, the recharging rate of an EV battery in CS per unit time is denoted by $g^2$. The constraints related to customers and satellite time windows are presented below. 

\begin{equation}\label{eq29}
\tau_{ik}+(t^1_{ij}+s_i)x_{ijk}\le \tau_{jk}+l_0(1-x_{ijk}), \quad \quad \forall i \in D\cup S\cup S', j \in S\cup S' \cup D', k\in K^1
\end{equation}

\begin{equation}\label{eq30}
    0\le \tau_{sk}\le l_s, \quad \quad \forall s \in S, k\in K^1
\end{equation}

\begin{equation}\label{eq31}
    \tau^1_{sk}\le \tau^2_{sk'}, \quad \quad \forall s \in S, k\in K^1,k' \in K^2
\end{equation}

\begin{equation}\label{eq32}
    \tau'_{ik}+(t^2_{ij}+s_i)y_{ijks}\le \tau'_{jk}+l_s(1-y_{ijks}), \quad \quad \forall s\in S, i \in C\cup S, j \in C\cup F\cup S', k\in K^2
\end{equation}

\begin{equation}\label{eq33}
\tau'_{ik}+t^2_{ij}y_{ijks}+g^2(B^2-b^{2+}_{isk})\le \tau'_{jk}+(l_s+g^2B^2)(1-y_{ijks}), \quad \quad \forall s\in S, i \in F, j \in C\cup F\cup S', k\in K^2
\end{equation}

\begin{equation}\label{eq34}
    e_i\le \tau'_{ik}\le l_i, \quad \quad \forall i \in C, k\in K^2
\end{equation}

The arrival and departure time of \nth{1}-echelon vehicles could be tracked by constraint \eqref{eq29}. In this constraint, $l_0$ means the latest time a \nth{1}-echelon vehicle could return to the depot after delivery. Constraint \eqref{eq30} guarantees that a satellite must be visited by a \nth{1}-echelon vehicle before the satellite delivery deadline (when the satellite may be closed). It also shows the non-negativity of variable $\tau_s$ for all $s\in S$. Constraints \eqref{eq31} ensure that the \nth{2}-echelon can not start its trip from a satellite unless that satellite has received its deliveries from the \nth{1}-echelon vehicles. It synchronizes the moving time of the \nth{1}- and \nth{2}-echelon vehicles at the satellites. Constraints \eqref{eq32}-\eqref{eq33} track the arrival and departure time of the \nth{2}-echelon vehicles when finishing the delivery service at a customer or recharging at a CS node. Constraint \eqref{eq34} guarantees that a customer node must be served within his/her time windows via a \nth{2}-echelon vehicle.

%Moreover, \cite{bakach2021two} proposed a two-tier urban vehicle routing problem with time windows, where a truck departs from a depot to deliver goods to the robot (as the \nth{2}-echelon delivery vehicle) hubs. The assumption is that the robot hubs (i.e., satellites) are sufficiently intelligent to assign packages to individual robots, delivering the packages to customers within specified time windows. The objective is to minimize the number of robot hubs needed and the shipping costs of robots and first-echelon trucks. Additionally, the number of robots at each robot hub is limited. Each robot has a limited route duration (due to a maximum battery range) and length (due to limited hours in a shift) and can operate daily, assuming that each robot delivers only one package per trip. The robots are fully recharged upon arrival at the robot hub. 

Furthermore, \cite{caggiani2021green} proposed a 2E-EVRPTW-PR that addresses a 2E-EVRP accounting for time windows and partial recharging. They also assumed two sets of customers related to the first and second echelons. A central aspect of this model is the incorporation of strict time window constraints, where each customer in the first and second echelons must receive their delivery within a specific time frame. This means that not only must the routes be optimized for efficiency, but the timing of each delivery must be carefully managed to ensure that all time windows are met while balancing travel and recharging needs. In addition,  

%In addition, \cite{liu2021hybrid}, a two-tier distribution network is proposed, involving depots, satellites, customers, and a mix of vehicle fleets. In the first tier, vans transport parcels from depots to satellites, where they are unloaded before returning to the depot. In the second tier (driverless delivery robots), DDRs deliver parcels from satellites to customers, returning to the same satellite. Customers are assigned to satellites based on proximity, demand, and time windows, prioritizing those with stricter time windows and higher demands. Due to capacity constraints, some orders may be denied and rescheduled if the penalty is lower than the cost of opening new satellites. Each vehicle is assigned a daily route, starting and ending at the same depot or satellite. One vehicle fulfills customer and satellite demands, and packages cannot be split across multiple vehicles. Operational and emissions costs are calculated per mile, with handling costs per package. The model assumes homogenous packages, fixed service times, and known vehicle capacities, including DDR battery life. Direct shipments from depots to customers are not allowed, and soft time windows offer some flexibility in scheduling.

In other work, \cite{wang2021two} proposed 
a 2E-EVRPTW-BSS model features a two-tiered logistics structure, where the first tier uses heavy-duty internal combustion vehicles (ICVs) to transfer goods from a central depot to satellite locations, and the second tier employs EVs to deliver goods from satellites to customers. Constrained by their limited battery capacity, EVs must visit BSSs to recharge during their journeys. To maintain punctual deliveries, time windows for customers are incorporated, ensuring that ICVs reach satellites within specified time frames. The primary objective is minimizing overall costs, including transportation, vehicle maintenance, and battery swapping expenses. The model operates under several key constraints: direct deliveries from depots to customers are prohibited. ICVs may make multiple stops at satellites in the first tier, and only the first tier allows split deliveries. Additionally, strict adherence to time window requirements is necessary, and direct travel between satellites is not allowed in the second tier. Assumptions include a sufficient supply at the depot to meet demand, and loading and service times at customer sites are integrated into the travel time.

\cite{akbay2022variable} described a system comprised of four key components: central warehouses, satellites, CSs, and customers. These elements are interconnected through two distinct echelons. The first echelon involves transporting products from the central warehouses to the satellites. In contrast, the second echelon focuses on delivering goods from the satellites to customers, possibly including CSs along the route. Each of these echelons incurs specific travel costs or distances, and customers are characterized by their demand levels and specific time windows, which define when deliveries must occur. The model duplicates each station to accommodate multiple visits to CSs, allowing EVs to recharge as needed. The number of replicas of every station is a critical parameter—too few copies could prevent optimal solutions, while too many could increase computational complexity. Based on initial experimentation, the number of copies has been set to three. The model is formulated as a three-index integer programming problem, with each vehicle type treated as a separate index to account for the heterogeneous nature of the fleets. This allows the model to effectively handle the varied characteristics of different vehicle types. The goal of the model is to minimize the aggregate traveled distance by all vehicles across both echelons, optimizing the efficiency of the entire transportation network.

The studies discussed in this section on the 2E-EVRP-TW share a common structure but differ in their approaches and focus. \cite{wang2021two} use EVs and recharge stations in the second echelon, emphasizing cost minimization, but \cite{wang2021two} adds stricter time window constraints. \cite{wang2019two} introduces BSSs and split deliveries, targeting cost reduction across transportation, handling, and penalties for delayed deliveries. \cite{caggiani2021green} focuses on green logistics with e-vans and e-bikes, integrating strict time windows. Lastly, \cite{akbay2022variable} optimizes the number of CSs by replicating them, offering a flexible approach to managing charging infrastructure, a feature less emphasized in other models. Each study introduces unique components, such as autonomous robots or flexible charging management, reflecting varying trade-offs in logistics optimization, sustainability, and cost efficiency. Despite their differences, all the studies aim to optimize vehicle routing and reduce overall costs in a two-echelon system involving electric or autonomous vehicles, focusing on time windows, customer demand, and charging or recharging infrastructure.

\subsubsection{2E-EVRP with partial recharging}
Partial recharging is a critical concept in EVRPs as it allows EVs to recharge only to the extent necessary to complete their current route or reach the next CS \citep{keskin2016partial}. This flexibility significantly enhances EV fleets' operational efficiency by reducing waiting times for full recharges and enabling vehicles to prioritize urgent deliveries \citep{zhou2021electric}. Additionally, partial recharging supports dynamic and real-time decision-making, allowing EVs to adapt to varying energy demands and unforeseen route changes, thereby improving service reliability and reducing overall energy consumption and costs in last-mile delivery operations. Only one paper has addressed the partial recharging in the 2E-EVRP in the literature. 

In \cite{caggiani2021green}, a 2E-EVRPTW-PR model is proposed to address 2E-EVRP-TW and partial recharging, introducing a more flexible approach to managing vehicle energy consumption. The model operates in two layers: the first echelon uses e-vans to transport goods within an urban area. In contrast, the second echelon employs e-bikes to serve deliveries in a restricted traffic zone. A unique feature of this model is its incorporation of multiple charging visits at designated transmission points that connect the two layers, allowing vehicles to recharge partially as needed during their routes. This flexibility in charging is crucial, especially for vehicles with limited battery capacities, as it ensures that they can continue their routes without running out of power. An essential aspect of the model is the enforcement of strict time window constraints for each customer, meaning that every delivery must be completed within a specified timeframe. This requirement necessitates careful management of both the travel routes and the timing of deliveries to minimize delays and penalties while considering the need for recharging. The model aims to reduce overall operational costs by balancing travel, vehicle maintenance, and recharging costs while ensuring that all customer time windows are met and energy needs are managed efficiently. By incorporating partial recharging, this approach provides a more practical and cost-effective solution for managing EVs in urban logistics, particularly in areas with restrictive traffic conditions. Notably, to include the partial recharging at the CSs, a new decision variable $Y^2_{ik}$ is defined, denoting the state of battery charge of \nth{2}-echelon vehicle $k\in K^2$ after leaving the node $i\in C\cup F\cup S$.  

\subsubsection{2E-EVRP with battery swapping stations}
BSSs differ primarily from traditional EV recharging stations in their approach to recharging EVs. While CSs require vehicles to plug in and wait for their batteries to recharge, which can take anywhere from 30 minutes to several hours, depending on the charging speed, BSSs offer a quick and seamless solution by replacing a depleted battery with a fully charged one in just a few minutes \citep{wang2019two}. This eliminates the wait time associated with charging and is particularly advantageous for time-sensitive operations like last-mile delivery or ride-sharing services. Additionally, swap stations often utilize standardized, shared batteries, enabling efficient battery management and reducing the upfront cost of EV ownership. In contrast, CSs are more suitable for private EV users and locations where vehicles can remain parked for extended periods, such as homes, offices, or public parking spaces. In the 2E-EVRPs in the literature, three papers were found to address the 2E-EVRP with BSSs, explained as follows. Note that the work by \cite{jie2019two} has studied the 2E-EVRP with BSS and \nth{1}-echelon EVs, described in Section \ref{2e-evrp-first-evs}. 

In \cite{wang2019two}, a 2E-EVRP-TW with battery swapping stations (2E-EVRPTW-BSS) is presented. This model addresses the limited EV battery life challenge by incorporating BSSs in the second echelon. As EVs deliver goods from satellites to customers, they visit BSSs to replace their depleted batteries, allowing them to continue their routes without long charging delays. The system ensures strict adherence to customer time windows while optimizing the number of battery swaps to reduce costs, including transportation, battery swapping, and maintenance. They defined parameter $c_b$ as the battery swapping cost of \nth{2}-echelon vehicles at each visited BSS while ignoring the battery swapping time.

Similarly, \cite{wang2021two} proposes a two-tier logistics structure where heavy-duty ICVs transport goods from depots to satellites, and EVs deliver goods from satellites to customers. Due to the limited range of EVs, they must visit BSSs for battery replacements during their deliveries. The model optimizes battery swap strategies to ensure punctual deliveries within specified time windows, aiming to minimize overall costs, including transportation, vehicle maintenance, and swapping expenses. The efficient management of BSSs is crucial for maximizing the system's operational efficiency and cost-effectiveness. The model operates under several key constraints: direct deliveries from depots to customers are prohibited. ICVs may make multiple stops at satellites in the first tier, and only the first tier allows split deliveries. Additionally, strict adherence to time window requirements is necessary, and direct travel between satellites in the second tier is not permitted. Assumptions include a sufficient supply at the depot to meet demand, and loading and service times at customer sites are integrated into the travel time.

Both \cite{wang2019two} and \cite{wang2021two} address two-echelon electric vehicle routing problems with time windows and battery swapping stations (BSSs), but they differ in their vehicle fleet composition and logistics structure. The former focuses on a system where EVs deliver goods from satellites to customers, requiring battery swaps at BSSs during deliveries. In this model, customer demands are indivisible, and the first echelon allows for split deliveries. On the other hand, the latter model uses ICVs in the first tier to transfer goods to satellites and EVs in the second tier to complete customer deliveries, emphasizing vehicle routing and battery swapping optimization. The key similarity is the focus on battery swapping as a critical element for ensuring the EVs can continue their routes without delays, optimizing operational efficiency and cost.

\subsubsection{2E-EVRP with simultaneous pickup and delivery}
``Simultaneous pickup and delivery (SPD)" in urban parcel logistics involves optimizing the transportation of goods to customers while simultaneously collecting items for return or further processing \citep{kocc2020review}. This approach addresses the growing demand for efficient and sustainable logistics in dense urban environments, where minimizing vehicle trips and operational costs is critical. By combining deliveries and pickups within the same route, companies can reduce empty vehicle mileage, improve resource utilization, and lower carbon emissions \citep{zhou2024pickup}. Advanced routing algorithms, such as those solving the SPD in 2E-EVRP, play a pivotal role in planning efficient routes, considering constraints like time windows, vehicle capacity, and customer priorities. This method supports green logistics initiatives and enhances customer satisfaction by offering streamlined, cost-effective services. Nevertheless, literature on 2E-EVRP with SPD is rare. Only one paper has studied this delivery approach. This shows a promising study gap in this research area for future work. 

\cite{akbay2023application} focused on the 2E-EVRP that involves SPD, called 2E-EVRP-SPD. In a two-tiered distribution framework, larger vehicles carry products from central warehouses to satellite locations, while smaller, environmentally-friendly vehicles handle the distribution of these goods to end customers. The problem also incorporates SPD requirements, which often arise from practices in reverse logistics. It considers SPD requirements for transporting goods to end customers. Taking SPD constraints into account, each customer might face two separate requirements: (i) the products that need to be transported to the demand location (i.e., delivery demand), and (ii) the products that have to be picked up from the demand location (i.e., pickup demand). Therefore, when a vehicle stops by a specific customer, both of these demands must be satisfied at the same time. To include pickup and demand in the 2E-EVRP, another type of demand is added to $q_i$ (i.e., the demand delivery) for each customer node $i\in C$ as $p_i$, representing the demand that must be picked up from their homes. In addition to $f_{ijk}$ and $f'_{ijks}$, two other decision variables $u_{ijk}$ and $u'_{ijks}$ representing the collected demands that picked up by the \nth{1}- and \nth{2}-echelon vehicles ($k\in K^1$ and $k\in K^2$, respectively) when traveling on arc $(i,j)\in A^1$ and $(i,j)\in A^2$, respectively. Also, two new decision variables on each satellite are defined: (i) $P_{sk}$ is the total demands the \nth{1}-echelon vehicle $k\in K^1$ must pick up from the satellite $s\in S$, and (ii) $\mu_s$ is the total demands that must be picked up by (collected to) satellite $s\in S$. Given these variables, to ensure the delivery and pickup demands are met, the following constraints can be incorporated into the 2E-EVRP model.

\begin{equation}\label{eq35}
P_{sk}=\sum_{j\in D\cup S: j\neq i}u_{sjk}-\sum_{j\in D\cup S: j\neq i}u_{jsk}, \quad \quad \forall s\in S, \forall k\in K^1   
\end{equation}

\begin{equation}\label{eq36}
u_{ijk}+f_{ijk}\le Q^1x_{ijk}\quad \quad \forall i,j\in D\cup S, i\neq j, \forall k\in K^1    
\end{equation}

\begin{equation}\label{eq37}
u_{0jk}=0, \quad \quad \forall j\in S, k\in K^1
\end{equation}

\begin{equation}\label{eq38}
\sum_{k\in K^1}P_{sk}=\mu_s, \quad \quad \forall s\in S     
\end{equation}

\begin{equation}\label{eq39}
\sum_{j\in C\cup S\cup F:j\neq i}\sum_{s\in S}\sum_{k\in K^2}u'_{jiks}=\sum_{j\in C\cup S\cup F:j\neq i}\sum_{s\in S}\sum_{k\in K^2}u'_{ijks}-p_i,\quad \quad \forall i\in C
\end{equation}

\begin{equation}\label{eq40}
u'_{ijks}+f'_{ijks}\le Q^2y_{ijks},\quad \quad \forall i,j \in C\cup S\cup F, i\neq j, \forall s\in S, \forall k\in K^2
\end{equation}

\begin{equation}\label{eq41}
u'_{sjks}=0, \quad \quad \forall s\in S, j\in C\cup F, k\in K^2
\end{equation}

\begin{equation}\label{eq42}
\sum_{i,j\in C\cup S\cup F:i\neq j}\sum_{k\in K^2}p_iy_{ijks}=\mu_s, \quad \quad \forall s\in S    
\end{equation}

Constraint \eqref{eq35} is the (pickup) parcel flow balance at the satellite in the first echelon. Constraint \eqref{eq36} and \eqref{eq40} could be replaced by the constraint \eqref{eq5} and \eqref{eq14}, respectively, ensuring that the total pickup and delivery demands by a \nth{1}- and \nth{2}-echelon vehicle must not exceed their load capacity. Constraint \eqref{eq37} and \eqref{eq41} enforce that the collected pickup demands by \nth{1}- and \nth{2}-echelon vehicles starting from the depot and satellites must be zero, respectively. Constraint \eqref{eq38} calculates the total demands picked up by each satellite, while constraint \eqref{eq42} connects the \nth{1}-echelon to \nth{2}-echelon by finding the total demands that must be picked up by a satellite in the \nth{2}-echelon. Constraints \eqref{eq38} and \eqref{eq42} guarantee that the number of pickup demands gathered from a satellite must equal the number of pickup demands that must be fulfilled by that satellite. Constraint \eqref{eq39} is the parcel flow conservation at customer nodes where the demand type is pickup.

\subsection{2E-EVRP with customer delivery in the first echelon}
\begin{figure}[ht!]
    \centering
    \fbox{\includegraphics[width=0.7\linewidth]{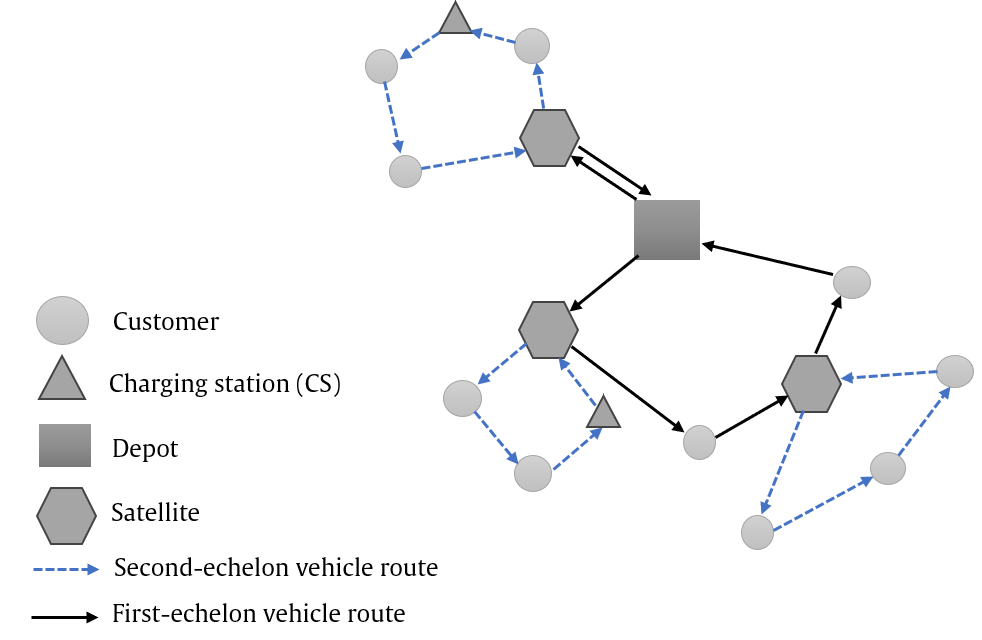}}
    \caption{Schematic example for 2E-EVRP with \nth{1} echelon customers}
    \label{fig:2E-EVRP-1st-customers}
\end{figure}

As given in Fig. \ref{fig:2E-EVRP-1st-customers}, the \nth{1}-echelon may also serve the customer directly in addition to home-attend delivery service via \nth{2}-echelon EVs. This approach is helpful when the demand is high, the satellites have limited EVs, and customer satisfaction is critical. It also uses first- and second-echelon vehicles to serve customers, making the deliveries fast, especially for customers on the way to the \nth{1}-echelon vehicles. Among the papers, only \cite{caggiani2021green} has considered this assumption in their 2E-EVRP. They assume that the customers that must be visited in the first and second echelons are known in advance. So, the \nth{1}-echelon customer can not be served by a \nth{2}-echelon vehicle and vice versa. They also consider time window constraints for the customers in both echelons. To include the customers in the first echelon routes, the arcs related to the first echelon must be changed to $A^1=\{(i,j)|i,j\in D\cup C\cup S,j\neq i\}$. Also, the expression of $\sum_{j\in C\cup S}\sum_{k\in K}\mathcal{F}^1x_{0jk}$ can be replaced the third term of objective function \eqref{eq1}. The following constraints are defined for the new assumption.

\begin{equation}\label{eq43}
 \sum_{j\in D\cup C\cup S: j\neq i}x_{ijk}=\sum_{j\in D\cup C\cup S: j\neq i}x_{jik}, \quad \quad \forall i\in D\cup C\cup S, \forall k\in K^1   
\end{equation}

\begin{equation}\label{eq44}
\sum_{j\in D\cup C\cup S: j\neq s}x_{sjk}\le 1,\quad \quad \forall s\in S, \forall k\in K^1 
\end{equation}

\begin{equation}\label{eq45}
w_{sk}=\sum_{j\in D\cup C\cup S: j\neq i}f_{jsk}-\sum_{j\in D\cup C\cup S: j\neq i}f_{sjk}, \quad \quad \forall s\in S, \forall k\in K^1   
\end{equation}

\begin{equation}\label{eq46}
f_{ijk}\le Q^1x_{ijk}\quad \quad \forall i,j\in D\cup C\cup S, i\neq j, \forall k\in K^1    
\end{equation}

\begin{equation}\label{eq47}
\sum_{k\in K^1}\sum_{j\in D\cup C\cup S:j\neq i}x_{ijk}+\sum_{k\in K^2}\sum_{s\in S}\sum_{j\in C\cup S\cup F:j\neq i}y_{ijks}=1, \quad \quad \forall i\in C    
\end{equation}

\begin{equation}\label{eq48}
q_{i}=\sum_{j\in D\cup C\cup S: j\neq i}f_{jik}-\sum_{j\in D\cup C\cup S: j\neq i}f_{ijk}, \quad \quad \forall i\in C, \forall k\in K^1   
\end{equation}

The constraints \eqref{eq43}-\eqref{eq47} have the functions of constraints \eqref{eq2}, \eqref{eq3}, \eqref{eq4}, \eqref{eq5}, and \eqref{eq8}, so they may replace them. Constraint \eqref{eq48} is the parcel flow conservation at customer nodes in the first echelon, ensuring that the difference between the cargo level of the \nth{1}-echelon vehicles when arriving and leaving a customer's home equals the demand quantity of that customer.    

%\subsection{2E-EVRP with delivery robots}Robot delivery vehicles are revolutionizing last-mile delivery by offering an efficient, cost-effective, and environmentally friendly alternative to traditional methods \citep{jennings2019study}. These autonomous or semi-autonomous robots are designed to navigate urban and suburban environments, delivering parcels directly to customers' doorsteps. They reduce dependency on human couriers and fossil-fuel-powered vehicles, lowering operational costs and carbon emissions \citep{srinivas2022autonomous}. Robot delivery vehicles excel in addressing challenges like traffic congestion, narrow streets, and limited parking spaces, making them ideal for densely populated areas \citep{moradi2024covering}. Additionally, their compact size and advanced navigation systems allow for faster, contactless deliveries, enhancing customer convenience and safety. As technology advances, integrating robot delivery vehicles into logistics networks promises to reshape last-mile delivery, fostering innovation and sustainability in urban logistics \citep{alverhed2024autonomous}.\cite{poeting2019comprehensive}, \cite{poeting2019simulation},\cite{bakach2021two}\cite{liu2021hybrid}\cite{alfandari2022tailored}\citep{liu2022physical}\cite{yu2024collaborative}

%\subsubsection{2E-EVRP with autonomous delivery vehicle}\cite{liu2020two} (similar to 2E-VRP)

\section{Objective Functions}\label{obj-functions-section}

In the context of the 2E-EVRP, several objective functions (as presented in Table \ref{objective-functions-table}) are utilized to evaluate and minimize various costs. These objective functions are explained in more detail in the following. 

\subsection{Transportation (Shipping) costs}
When moving goods between different locations, transportation (shipping) costs are incurred for first—and second-echelon vehicles. These include fuel, vehicle wear and tear, and third-party shipping fees. EV routing or shipping costs account for electricity consumption and the associated charging infrastructure. Minimizing transportation costs is essential for optimizing the two-echelon system, where goods are first transported to intermediate depots and then delivered to final destinations. This objective function, denoted by $\mathcal{O}_1$, is the most used one, widely discussed in the literature, including works such as \cite{breunig2017two,agardi2019two,breunig2019electric,jie2019two,wang2019two,affi2020general,wang2021two,zijlstra2021integrating,wu2023branch}. These studies explore methods to minimize transportation costs in the first and second echelons routing by optimizing routes, fleet size, and satellite utilization. To include this cost in the basic 2E-EVRP model, the first and second terms of the objective function \eqref{eq1} could be multiplied by parameters $c^1_{ij}$ and $c^2_{ij}$, respectively, to calculate the transportation cost related to the first- and second-echelon vehicles. These parameters are the vehicle traveling cost on arc $(i,j)$ per unit distance in the first ($A^1$) and second ($A^2$) echelons. Therefore, the transportation or shipping cost of \nth{1}- and \nth{2}-echelon vehicles could be represented by the following objective function \eqref{eq50}.

\begin{equation}\label{eq50}
    Min. \sum_{(i,j)\in A^1}\sum_{k\in K^1}c^1_{ij}d_{ij}x_{ijk}
    +\sum_{(i,j)\in A^2}\sum_{k\in K^2}\sum_{s\in S}c^2_{ij}d_{ij}y_{ijks}
\end{equation}

\subsection{Traveled distances}
Traveled distances are the total distance vehicles travel within the routing system. In the 2E-EVRP, minimizing traveled distances is crucial for reducing fuel consumption (or electricity usage) and time and ensuring system efficiency. By optimizing the routing of EVs, the distance between intermediate depots and customer locations can be reduced, leading to cost savings and operational improvements. It may also be treated as a particular instance of transportation cost when the $c^1_{ij}=c^2_{ij}=1$. The first and second terms of the objective function \eqref{eq1} represent the traveled distance by the \nth{1}- and \nth{2}-echelon vehicles. Two works by \cite{akbay2022variable,akbay2023application} explore this objective function.

\subsection{Fixed vehicle usage costs}
Fixed vehicle usage costs refer to expenses incurred when using a vehicle in the system, irrespective of the distance traveled. These include maintenance, vehicle depreciation, and fixed operational costs like insurance and registration. For EVs, additional fixed costs may include battery maintenance and charging infrastructure. Studies such as \cite{wang2019two,affi2020general,caggiani2021green,wang2021two,akbay2023application,wu2023branch} address the minimization of fixed usage costs of \nth{1}- and \nth{2}-echelon vehicles. Unlike the transportation cost of vehicles, which is operational or tactical, the fixed usage costs of vehicles are strategic and have a long-term impact on the total costs. The third and fourth terms of the objective function \eqref{eq1} represent this cost.  

\subsection{Fixed satellite utilization costs}
Fixed satellite utilization costs represent the expenses of maintaining and operating satellites in the 2E-EVRPs. These costs include each satellite's rent, maintenance, and staffing expenses, regardless of the volume of goods processed. Optimizing the number of satellites and their usage can minimize these costs. This cost is only explored in \cite{caggiani2021green} defining a parameter $c_m$ as a transshipment point cost, which was added to the objective function.

\subsection{Handling costs}
Handling costs refer to the expenses associated with physically handling goods at satellites, including loading and unloading operations, storage, and related labor costs. Minimizing these costs can increase overall system efficiency, especially in a two-echelon system with multiple satellite facilities. This is one of the main cost factors in 2E-EVRPs as a satellite-related cost. This cost is considered in works by \cite{jie2019two,wang2019two}, where handling costs are incorporated into optimization models for minimizing overall operational expenses. The last term of the objective function \eqref{eq1} represents this cost. 

\subsection{Energy/Battery consumption costs}
Energy or battery consumption costs reflect the expenses related to the battery swapping, recharging, and energy consumed by EVs while traveling and operating within the 2E-EVRPs. These costs can vary based on the vehicle's energy efficiency, the distance traveled, and the charging infrastructure used. Minimizing energy consumption is vital for reducing operational costs and ensuring the system's sustainability. This cost function is defined as a battery swapping cost in the BSSs in studies such as \cite{jie2019two,wang2019two,wang2021two}, where they defined two parameters $c^1_b$ and $c^2_b$ as the battery swapping cost in the station $b\in F$ at the first and second echelons, respectively. Multiplying these parameters to the times that a BSS is used in the first and second echelons obtains the total battery swapping cost, presented by the objective function \eqref{eq51}.

\begin{equation}\label{eq51}
    Min. \sum_{b\in F}\sum_{(b,j)\in A^1}\sum_{k\in K^1}c^1_{b}x_{bjk}
    +\sum_{b\in F}\sum_{(b,j)\in A^2}\sum_{k\in K^2}\sum_{s\in S}c^2_{b}y_{bjks}
\end{equation}

In addition, the energy consumption cost is considered by \cite{caggiani2021green} as the total electric energy cost. To do this, they defined the parameters $c^1_e$ and $c^2_e$ as the electric energy cost of \nth{1}- and \nth{2}-echelon vehicles per unit distance, respectively. Then, multiplying these parameters by the total traveled distances, they calculated the total consumed (electric) energy cost. The objective function \eqref{eq52} finds this cost in both echelons.

\begin{equation}\label{eq52}
    Min. \sum_{(i,j)\in A^1}\sum_{k\in K^1}c^1_{e}d_{ij}x_{ijk}
    +\sum_{(i,j)\in A^2}\sum_{k\in K^2}\sum_{s\in S}c^2_{e}d_{ij}y_{ijks}
\end{equation}

\subsection{Drivers' wages}
Drivers' wages refer to the salaries paid to the system's drivers, who operate the vehicles. These wages are often a significant part of the operational costs, and minimizing the number of drivers required or optimizing their schedules can reduce this cost. Careful scheduling of drivers and vehicles can help lower this expense in two-echelon systems. This cost was only assumed in work by \cite{caggiani2021green}, where the authors incorporate drivers' wages into the optimization function to improve cost efficiency. They defined the parameters $c^1_d$ and $c^2_d$ as the fixed wage paid to the drivers per trip in the first and second echelons, respectively. The objective function \eqref{eq53} calculates the drivers' wages in the 2E-EVRPs.

\begin{equation}\label{eq53}
    Min. \sum_{j\in S}\sum_{k\in K^1}c^1_dx_{0jk}+\sum_{s\in S}\sum_{j\in C\cup F}\sum_{k\in K^2}c^2_dy_{sjks} + \sum_{s\in S}h_s\lambda_s
\end{equation}

\section{Solution methodologies}\label{solution-methodology}
This section reviews the selected works on 2E-EVRP, including their solution methodology and techniques to model and tackle small and large instances. It also discusses the scalability of the proposed methods by examining the size of the solved cases and comparing the algorithms in the related benchmarks. The solution methodologies are categorized into four classes: (i) mathematical programming-based methods, which are mainly solved by the proposed model by a commercial exact solver like Gurobi or CPLEX (Section \ref{math-model-based}); (ii) exact methods, which are developed specifically for the problem e.g., branch-and-price (BP) (Section \ref{exact-methods}); (iii) heuristic/metaheuristics-based methods, in which a problem-tailored metaheuristic solution method is designed (Section \ref{Metaheuristic-methods}); and (iv) metaheuristic-based methods which combine exact methods with metaheuristics (Section \ref{matheuristic-method}).  

\subsection{Mathematical programming-based}\label{math-model-based}
In the literature on 2E-EVRP, the following papers have used a commercial exact solver to tackle their problems: Gurobi \citep{wang2019two,wang2021two}, CPLEX \citep{jie2019two,caggiani2021green,akbay2022variable,akbay2023application,wu2023branch}. Mainly, the commercial exact solvers were used to verify the performance of the proposed solution methods (e.g., metaheuristics) on the small-sized problem instances. Table \ref{excat-solver-comp} summarizes the exact solvers used in each paper alongside the benchmarks, the size of the solved instances with those solvers, and the related operating machines. In this table, $\mathcal{N}$, $\mathcal{N}^*$, $\mathcal{T}$ are the number of instances solved by the exact solver, the number of instances solved to optimality, and the maximum running time of the solver in seconds. This table shows that after publishing the work by \cite{breunig2017two} and introducing the 2E-EVRP dataset by them, most papers use this dataset to solve and evaluate the performance of the proposed methods. Also, the 2E-EVRP instances are generated by modifying the 2E-VRP and EVRP-TW cases in the literature. The largest size for a 2E-EVRP instance tackled by an exact solver is $|V|=27$ by \cite{akbay2023application}; however, they do not present an optimal solution. Also, \cite{wang2019two,caggiani2021green,wu2023branch} found the exact solution for all instances solved. The maximum running time of 7200 seconds (two hours) is the most common setting in the literature for exact solvers.     

\begin{table}[ht!]
\centering
\caption{Summary of work which applied a commercial exact solver}\label{excat-solver-comp}
\begin{adjustbox}{max width=\textwidth}
\begin{tabular}{llllllll}
\hline
\textbf{Ref.}& \textbf{Exact solver} & \textbf{Operating machine}& \textbf{Dataset}& \textbf{Size of instances} & $\mathcal{N}$ & $\mathcal{N}^*$ & $\mathcal{T}$ \\
\hline
\citep{wang2019two}    & Gurobi       & Intel Xeon E3-123v3 with 8 GB RAM& Randomly-generated& $|V|=[7,23]$    & 14& 14& N/A      \\
\citep{jie2019two}        & CPLEX        & N/A& 2E-VRP \citep{perboli2011two}         & $|V|=[14,19]$   & 19 & 14& 7200     \\
\citep{wang2021two}       & Gurobi       & Intel Xeon E3- 123v3 with 8 GB RAM& 2E-EVRP \citep{breunig2019electric}   & $|V|=[7,20]$    & 15& 10& 4500     \\
\citep{caggiani2021green} & CPLEX        & Intel(R) Core (TM) i7-8550U CPU (1.80GHz) and 16GB of RAM & EVRP-TW \citep{goeke2019vrptw}& $|C|=[10,15]$   & 12& 12 & 7200     \\
\citep{akbay2022variable}& CPLEX        & Intel Xeon CPU 5670 CPUs with 12 cores of 2.933 GHz and 32 GB RAM & EVRP-TW \citep{schneider2014electric} & $|C|=[5,15]$   & 36 & 21 & 7200     \\
\citep{akbay2023application}& CPLEX        & Intel©R Xeon©R 5670 CPUs with 12 cores of 2.933 GHz and 32 GB RAM & 2E-EVRP \citep{breunig2019electric}   & $|V|=27$            & 12& 0& 43200    \\
\citep{wu2023branch} & CPLEX        & Inter(R) Core(TM) i7-9700K CPU (3.6GHz) and 32GB of RAM & 2E-EVRP \citep{breunig2019electric}   & $|V|=14$& 12& 12& 3600\\
\hline
\end{tabular}
\end{adjustbox}
\end{table}

\subsection{Exact Methods}\label{exact-methods}
Among the papers, only two works have proposed exact methods: (i) \cite{breunig2019electric} developed a mathematical programming algorithm, adapted from \cite{baldacci2013exact}, and (ii) a BP algorithm developed by \citep{wu2023branch}. 

In \cite{breunig2019electric}, the solution methodology for solving the 2E-EVRP to optimality is built upon an exact method adapted from \cite{baldacci2013exact}. It is specifically tailored to handle the problem's multigraph representation, incorporating the constraints and complexities of EV routing. The method is divided into two main steps: (i) ``Generation of First-Level Routes and Initial Lower Bound": First, all potential first-level routes are enumerated, and a lower bound (LB0) is calculated using an integer relaxation of the problem, which is reformulated as a multiple-choice knapsack problem. Second, the relaxation extends the ng-route relaxation to account for the unique characteristics of the multigraph representation, incorporating the limitations on battery capacity and recharging stations. (ii) ``Evaluation of First-Level Route Subsets": At first, for every subset of first-level routes, the following sequence is executed: ``Subset-Specific Lower Bound (LB1)," in which the selected set of first-level routes is fixed, and a new lower bound is calculated. If this bound exceeds the best-known solution, the subset is discarded. Then, in the ``Solution to Second-Level Routing Problem" step, if the subset is feasible, the remaining problem (a capacitated multi-depot vehicle routing problem, or MDCVRP) is solved optimally using exact methods such as dynamic programming. The optimal 2E-EVRP solution is obtained by minimizing the total cost across all subsets.

Throughout this process, specialized algorithms efficiently generate elementary routes and evaluate subsets. The method also integrates bounding functions to reduce the computational burden by pruning unpromising subsets early in the search. It relies on practical lower and upper bounds to guide the enumeration process and minimize computational costs. For larger instances, the approach benefits from heuristic initialization called ``large neighborhood search (LNS)" (as described in Section \ref{Metaheuristic-methods}) to provide strong upper bounds that enhance overall efficiency.

Moreover, the BP (branch-and-price) algorithm presented in \cite{wu2023branch} addresses the 2E-EVRP-BSS by combining a master problem and a pricing subproblem to achieve optimal solutions. The methodology leverages the arc flow model, decomposing the problem into manageable components through Dantzig–Wolfe decomposition. It introduces a restricted master problem (RMP) and applies column generation (CG) to refine solutions iteratively. The algorithm operates in three main stages: (i) Initialization: Feasible first-echelon routes are enumerated using heuristics like the insertion heuristic to generate initial subsets of second-echelon routes. These initial routes serve as columns in the RMP, which is solved to provide an initial lower bound. (ii) CG: The master problem is solved iteratively, with dual variables extracted to construct and solve a pricing subproblem for generating columns with reduced costs. Then, the pricing subproblem for the second echelon is formulated as an Elementary Shortest Path Problem with Resource Constraints (ESPPRC), solved using dynamic programming and a bidirectional labeling algorithm to accelerate the computation. (iii) BP framework: If the solution to the RMP is fractional, the algorithm applies a branching strategy to divide the problem into subproblems, creating a search tree. The branching rules involve constraints on variables, vehicle assignments, or arcs in the solution, ensuring integer solutions in subsequent iterations. Then, depth-first search minimizes memory usage while exploring the search tree.

The BP algorithm by \cite{wu2023branch} integrates acceleration strategies such as state-space relaxation and ng-route relaxation to reduce computational complexity and enhance scalability. Numerical experiments demonstrate that the BP algorithm is highly effective, outperforming traditional solvers like CPLEX on small- and medium-sized test instances, achieving optimal solutions with smaller computational gaps and reasonable runtime. The BP solved the newly generated 2E-EVRP-BSS instances, adapted from 2E-EVRP \citep{breunig2019electric} with up to 20 customers.

\subsection{Metaheuristic Approaches}\label{Metaheuristic-methods}
Metaheuristics are advanced optimization techniques designed to find near-optimal solutions for complex problems within a reasonable timeframe, particularly when exact methods become computationally infeasible. These methods have been successfully developed for combinatorial optimization problems like VRPs \citep{prodhon2016metaheuristics} EVRPs \citep{reddy2022meta}, and 2E-VRPs \citep{crainic2008clustering,sluijk2023two}. Remarkably, the following metaheuristic-based algorithms have been proposed for 2E-EVRPs: variable neighborhood search (VNS) \citep{wang2021two,akbay2022variable}, general VNS (GVNS) \citep{affi2020general}, genetic algorithm (GA) \citep{agardi2019two}, LNS \citep{breunig2017two,breunig2019electric,zijlstra2021integrating}, and hybrid metaheuristic Construct, Merge, Solve \&
Adapt (Adapt-CMSA) \citep{akbay2023application}. The 2E-EVRPs often involve large solution spaces, intricate constraints, and multiple objectives, such as minimizing costs, energy consumption, or environmental impact. Metaheuristics provide flexibility and scalability to address such challenges, enabling researchers to model realistic scenarios, integrate real-life constraints, and explore trade-offs between efficiency and sustainability in modern logistics and transportation systems. These metaheuristics for the 2E-EVRPs are explained in more detail in the following. In addition, the proposed metaheuristics in the literature are summarized in Section \ref{summary-meta}, discussing their similarities and distinctions with the datasets (instances) they solved.

\cite{breunig2017two} is the first work that proposed a metaheuristic-based solution method for the basic 2E-EVRP. In this work, the proposed methodology relies on a metaheuristic based on LNS \citep{shaw1998using}. This approach iteratively refines solutions through partial destruction (destroy) and reconstruction (repair), combined with a labeling algorithm for CS optimization and a local search for route enhancement. The process begins with generating an initial feasible solution using a repair operator that meets capacity constraints. The CS placements are initially disregarded and are later refined using a dynamic programming-based labeling algorithm. The destruction phase employs five operators that strategically dismantle parts of the solution. For example, related customers are removed, random routes are eliminated, or satellites are temporarily closed to allow their reassignment in subsequent repair phases. Rebuilding the solution employs a simplified ``Cheapest Insertion Algorithm" that reinserts customers based on cost considerations, with adjustments made for high-demand customers to ensure feasibility.

Once second-level routes are finalized, the methodology constructs first-level routes to supply satellites with the necessary goods. In this phase, split deliveries are allowed, and preprocessing steps segment large demands into manageable truckloads. Local search techniques, such as 2-opt and relocate, further refine the solution by optimizing route quality while dynamically adjusting CS placements when potential improvements are identified. The labeling algorithm identifies optimal CS positions by evaluating non-dominated paths, ensuring energy constraints are met while minimizing costs. This iterative process alternates between destruction, repair, and optimization phases until a termination criterion, such as a time limit, is reached. It is designed to handle large instances efficiently by concentrating computational effort on the most promising parts of the solution space. Including multiple destroy and repair operators ensures diverse solution exploration, avoiding premature convergence and enhancing robustness. 

Moreover, the solution methodology in \cite{agardi2019two} employs a combination of ``construction and improvement heuristics'' to solve the 2E-EVRP. The process begins with ``construction algorithms'' inspired by the traveling salesman problem (TSP) \citep{flood1956traveling}. These algorithms iteratively build initial solutions by selecting customers and creating routes. The ``Nearest Neighbor algorithm" starts from a random customer and selects the nearest unvisited customer at each step, minimizing local distances. The ``Arbitrary Insertion algorithm" randomly chooses a customer and inserts them between two selected nodes in a way that reduces the route cost. These algorithms produce a complete, albeit initial, set of routes for further refinement. Improvement algorithms are then applied to enhance the initial solutions. They used a GA \citep{holland1992genetic} simulating natural evolutionary processes. It operates on a population of solutions, evaluating them based on a fitness function. New solutions (offspring) are generated and iteratively improved through operators like crossover and mutation. The study applies specific crossover techniques such as ``order crossover (OX)," ``cycle crossover (CX)," and ``partially matched crossover (PMX)," as well as the ``2-opt operator" for mutations. Another improvement method is the ``Hill Climbing algorithm," which iteratively improves a single solution by evaluating its neighbors and replacing it with a better neighbor if one exists. The 2-opt operator is also used to optimize edge exchanges in routes.

The methodology integrates these heuristics to solve the problem in phases: generating initial solutions, improving them through iterative refinements, and optimizing route assignments for recharge stations and satellites. Computational tests show that combining construction heuristics with improvement algorithms, particularly GA-enhanced solutions, yields superior performance in minimizing travel costs and achieving feasible routes under given constraints. 

Furthermore, \citep{breunig2019electric} proposed an LNS-based algorithm, taken from \cite{shaw1998using}, to provide good upper bounds for their exact method. Their approach iteratively improves solutions by alternating between destruction, reconstruction, and optimization phases. The general structure of the process follows a ``ruin and recreate" principle, which systematically destroys parts of the solution and repairs it using sophisticated operators and local search techniques. The process begins with a destruction phase, where three leading operators are used to partially dismantle the solution. These include removing related customers, random routes, and temporarily closing satellites. Additional operations may involve reopening all satellites or removing inefficient single-customer routes. Following the destruction, the reconstruction phase reinserts customers into routes, reconstructs first-level and second-level routes, and optimally places CSs. This reconstruction uses a simplified, cheapest insertion algorithm for node placement and dynamic programming to address battery constraints. After reconstruction (or repair), a local search phase is applied to refine the solution using neighborhood moves such as 2-opt, relocate, and swap. These moves are systematically evaluated for cost improvements, and adjustments are made to CS placements to ensure feasibility. The local search emphasizes computational efficiency by filtering moves before detailed evaluations. The proposed LNS method iterates until a termination criterion, such as a time limit or lack of improvement over several iterations, is met. The algorithm is characterized by its simplicity, robustness, and efficiency, leveraging a small set of operators without relying on adaptive mechanisms. This approach has demonstrated strong performance in producing high-quality solutions for large 2E-EVRP instances.

In addition, \cite{affi2020general} applies a GVNS-based algorithm, adapted from \cite{mladenovic2012general}, to address the 2E-EVRP. This approach integrates various neighborhood structures for local search and perturbation to iteratively improve solution quality. The methodology is structured around exploring first-echelon, second-echelon, and two-echelon neighborhoods. The GVNS starts with a randomly generated initial solution and leverages ten neighborhood structures. These structures are categorized into three groups: first-echelon (e.g., shifting and swapping satellites), second-echelon (e.g., shifting and swapping customers, modifying CSs), and two-echelon (e.g., opening, closing, or replacing satellites). Each neighborhood explores potential improvements by shifting, swapping, adding, or removing elements within the solution. The local search component uses a ``variable neighborhood descent (VND)" algorithm \citep{duarte2018variable}, which systematically explores the defined neighborhoods until no further improvement can be achieved. This involves sequentially applying descent methods to the first, second, and two-echelon neighborhoods. In the shaking phase, perturbations are introduced by selecting and repositioning customer nodes to escape local optima. The GVNS emphasizes flexibility and robustness by applying diverse neighborhood structures and maintaining a systematic exploration strategy. Computational experiments show that it outperforms LNS \citep{breunig2019electric} in both solution quality and computational efficiency.

Moreover, the solution methodology presented by \cite{wang2021two} addresses the 2E-EVRPTW-BSS using a metaheuristic based on the VNS algorithm \citep{mladenovic1997variable}. This approach optimizes routes for both echelons by iteratively exploring and improving solutions. The algorithm starts with generating an initial solution for the second echelon. Routes are constructed by assigning customers to satellites and iteratively adding them to routes based on proximity and angular sorting. BSSs are inserted into routes as needed to ensure compliance with energy constraints, and additional routes are created when capacity violations occur. Time window requirements are also integrated during this phase to maintain feasibility. The algorithm applies a series of inter-route and intra-route neighborhood structures to generate diverse solutions in the shaking phase. These structures include operators like 2-opt*, relocate, exchange, swap, shift, and satellite-change, each designed to modify specific elements of the solution for enhanced exploration. The shaking phase aims to escape local optima by perturbing the current solution. The local search phase refines the solutions generated in the shaking phase by systematically cyclically applying neighborhood operators. These operations optimize customer allocations and BSS placements while respecting time windows, capacity, and energy constraints. The process continues until no further improvements can be made within each neighborhood. The VNS algorithm's iterative nature ensures a thorough exploration of the solution space through a decoding process; the approach integrates first and second-echelon routes, linking the required demands and time constraints of satellites with the optimized routes of the second echelon. This interconnection enhances the solution's overall efficiency and feasibility. The algorithm terminates once all shaking neighborhoods have been explored or a predefined termination criterion, such as maximum runtime, is reached. Computational results demonstrate the proposed VNS's effectiveness in generating high-quality solutions within a reasonable timeframe, outperforming another metaheuristic proposed by \cite{wang2017matheuristic} in scalability and computational efficiency.

Furthermore, \cite{zijlstra2021integrating} applies an extended LNS heuristic, denoted as LNS-E2E*, to solve the 2E-EVRP with considerations for non-linear charging functions and partial recharging. The methodology consists of two main phases: replicating the original LNS-E2E algorithm proposed by \cite{breunig2019electric} and extending it to incorporate charging costs. In the first phase, the replication involves the implementation of LNS-E2E*, a modified version of the original algorithm that does not include local search enhancements due to time constraints. The LNS-E2E* algorithm operates iteratively, using destruction and reconstruction phases to explore the solution space. During the destruction phase, customers and routes are strategically removed using operators such as related node removal, random route removal, and satellite closure. In the reconstruction phase, removed customers are reinserted using the cheapest insertion heuristic to minimize additional costs. This process ensures feasible solutions by adhering to capacity and battery constraints. The second phase extends the methodology by integrating non-linear charging costs and partial recharging. A non-linear charging function is adopted to reflect realistic battery charging behaviors, where charging speed decreases as the battery level increases. Partial charging strategies are incorporated to optimize the charging time and cost by allowing vehicles to charge only to the level needed for the next route segment rather than fully charging. This modification is implemented during the reconstruction phase of LNS-E2E*, enabling the algorithm to evaluate and choose routes based on total costs, including distance, usage, and charging costs. Since the VNS developed by \cite{zijlstra2021integrating} is less advanced than the VNS by \cite{breunig2019electric} (i.e., there are no local search algorithms within the VNS compared to VNS by \cite{breunig2019electric}), the solution quality of the instances reported by \cite{zijlstra2021integrating} are worst than the related results found by \cite{breunig2019electric}.

Moreover, The solution approach by \cite{akbay2022variable} introduces a hybrid methodology that combines an extended ``Clarke and Wright Savings Heuristic" with a VNS to solve the 2E-EVRP-TW. The methodology addresses the complexities of multi-tier logistics, time window constraints, and EV limitations. The approach begins with the ``Clarke and Wright Savings Heuristic," which generates an initial solution by assigning customers to their nearest satellites and constructing routes for both echelons. For second-echelon routes, it creates direct tours for customers, inserting CSs when necessary to address battery feasibility. A savings-based merging mechanism iteratively consolidates routes, guided by distance and demand-based savings metrics. Any unvisited customers are handled through a greedy insertion operator to ensure completeness. This initial solution is the starting point for the VNS algorithm, which systematically explores the solution space using multiple neighborhood structures. The algorithm alternates between a shaking phase, where routes are perturbed using operators like random customer or route removal, and a local search phase, which applies inter-route and intra-route operators such as relocation, 2-opt, and CS reinsertion. To escape local optima and explore diverse solutions, the shaking phase employs large-scale removal and reinsertion operators inspired by LNS.

The methodology uses an extended objective function to address infeasible solutions by incorporating penalty terms for capacity, battery, and time window constraint violations. During the VNS process, solutions are evaluated and accepted probabilistically, with a dynamic temperature-based mechanism inspired by simulated annealing to control exploration and exploitation. The hybrid approach effectively balances intensification and diversification in the search process, outperforming the exact solve CPLEX and Clarke and Wright savings heuristic on benchmark instances of the 2E-EVRP-TW. Combining savings-based initialization with VNS-driven refinement allows the methodology to handle the problem’s multi-tier structure and operational constraints efficiently.

Recently, the Adapt-CMSA algorithm detailed in \cite{akbay2023application} integrates problem decomposition with adaptive processes to solve the 2E-EVRP-SPD. This approach systematically constructs and improves solutions by focusing on subproblem optimization and adaptive learning for dynamic solution refinement. The methodology begins by decomposing the problem into smaller, manageable subproblems, each representing a portion of the overall routing system. These subproblems are solved independently using a heuristic solver, generating initial feasible solutions for both echelons. In the construction phase, initial solutions are created for the problem's routing and allocation components, ensuring capacity, battery, and SPD constraints are respected. The merging phase combines subproblem solutions to form a comprehensive solution for the entire routing system. This process leverages heuristics to ensure route consistency while minimizing costs and satisfying operational constraints. In this step, the algorithm identifies synergies between solutions to improve the overall objective. The adaptation phase involves iterative improvement through learning and modification. The algorithm applies adaptive techniques to refine solutions based on performance feedback. It identifies critical areas for improvement, such as vehicle assignments, time windows, or CS placements, and adjusts them dynamically. This phase also incorporates a feedback loop where successful modifications are reinforced to guide future iterations.

The solving phase integrates exact and heuristic methods to optimize the final combined solution. A mixed-integer programming model or a local search heuristic resolves inconsistencies or inefficiencies, ensuring the solution is globally feasible and cost-efficient. Adapt-CMSA demonstrates significant flexibility by adapting to evolving problem constraints and improving solution quality iteratively. Its ability to decompose and reassemble solutions while adapting dynamically allows it to effectively address the complexities of the 2E-EVRP-SPD, resulting in high-quality solutions while outperforming the probabilistic Clark and Wright savings algorithm (pC\&W) and CPLEX in medium to large-size instances.

\subsubsection{Summary on metaheuristics}\label{summary-meta}
This section compares the proposed metaheuristics proposed for 2E-EVRPs in the literature. These comparisons are given in Tables \ref{meta-comp-table} and \ref{meta-problems-table}. Table \ref{meta-comp-table} presents the similarities and differences between the proposed metaheuristics regarding the method, initial solution construction, perturbation or shaking phase, improvement phase, stopping criteria, and constraints handling approaches. A significant similarity among the methods lies in their approaches to generating initial solutions. Most studies rely on structured or heuristic-based techniques, such as the cheapest insertion algorithm, as seen in \cite{breunig2017two,breunig2019electric,zijlstra2021integrating}. Others employ randomized strategies, including the approaches in \cite{agardi2019two,affi2020general,akbay2023application}. Also, the shaking phase in the metaheuristics shows a trend of using techniques designed to introduce diversity into the search process. Commonly employed methods include node removal, route removal, and satellite-based operations, evident in studies like \cite{breunig2017two,zijlstra2021integrating,akbay2022variable}. Improvement phases often utilize classic operators such as 2-opt, relocate, swap, and exchange. Advanced strategies such as integrating set-covering MILP models in \cite{akbay2023application} demonstrate a focus on enhancing solution quality through more targeted refinement. In terms of handling constraints, the use of repair operators is widespread. These are utilized to ensure the feasibility of solutions after perturbation, as seen in \cite{breunig2017two,zijlstra2021integrating,akbay2023application}. Alternatively, some methods employ penalty mechanisms for infeasible solutions, as shown in \cite{wang2021two,akbay2022variable}. Despite these similarities, the metaheuristics differ in several ways. For instance, the choice of metaheuristic framework varies significantly, with methods LNS \citep{breunig2017two,zijlstra2021integrating}, GA \citep{agardi2019two}, VNS \citep{wang2021two,akbay2022variable}, and Adapt-CMSA \citep{akbay2023application}. These differences reflect varying levels of complexity, flexibility, and suitability for different problem scales and constraints. The repair phases also show diversity, with some methods emphasizing dynamic programming or reinsertion algorithms \citep{breunig2017two,breunig2019electric}, while others employ simpler greedy or probabilistic strategies \citep{akbay2022variable,akbay2023application}. Local search also varies, with some methods relying on basic operators and others integrating advanced techniques, such as systematic CS relocations in \cite{breunig2019electric} or MILP-based refinements in \cite{akbay2023application}. 

\begin{sidewaystable}
\centering
\caption{Comparison of proposed metaheuristics for 2E-EVRPs in the literature}\label{meta-comp-table}
\begin{adjustbox}{max width=\textwidth}
\begin{tabular}{|l|l|p{3cm}|p{4cm}|p{4cm}|p{4cm}|p{3cm}|p{4cm}|}
\hline
\textbf{Ref.}  & \textbf{Metaheuristic} &\textbf{Initial solution generation}& \textbf{Shaking, destroy, or perturbation phase} & \textbf{Repair phase} & \textbf{Improvement or local search phase} & \textbf{Stopping criteria} & \textbf{Tackling the constraints} \\
\hline
\citep{breunig2017two} & LNS & Cheapest insertion algorithm with a labeling algorithm & \emph{Related node removal, random route removal, close satellite, open all satellites, remove single customer routes} & \emph{Dynamic programming labeling algorithm}& \emph{2-opt, 2-opt*, relocate, swap, swap 2-1} & Maximum number of iterations & Using repair operators \\
\hline
\citep{agardi2019two} & GA & Randomly generated and nearest neighbor algorithm& \emph{Order crossover (OX), cycle crossover (CX), partially matched crossover (PMX), 2-opt} & \emph{Insertion algorithms} & \emph{Elitism strategy, hill climbing algorithm} & Maximum number of iterations or time limit & Using repair operators \\
\hline
\citep{breunig2019electric} & LNS & Cheapest insertion algorithm with a labeling algorithm& \emph{Related node removal, random route removal, close satellite, open all satellites, remove single customer routes} & \emph{Reinsertion of customer visits, simplified cheapest insertion, dynamic programming labeling algorithm}& \emph{2-opt, 2-opt*, relocate, swap, swap 2-1,  systematic CSs relocations}& Maximum number of iterations and time limit & Using repair operators \\
\hline
\citep{affi2020general} & GVNS & Randomly-generated & \emph{Customer re-insertion} & - & \emph{Customer shift, swap, satellite closure, satellite open, satellite shift, satellite swap, satellite replace, CS replace, CS drop, CS addition} & When no improvement observed  & N/A\\ 
\hline
\citep{wang2021two} & VNS & Similar to \cite{schneider2014electric} for \nth{2} echelon& \emph{2-Opt*, inter-Relocate, inter-Exchange, inter-Swap, and inter-Shift, Satellite-Change} & - & \emph{intra-Relocate, intra-Exchange, intra-Swap, intra-Shift, Station-InRe} & Time limit & Penalizing infeasible solutions \\
\hline
\citep{zijlstra2021integrating} & LNS & Cheapest insertion algorithm & \emph{Related node removal, random route removal, close satellite, open all satellites, remove single customer routes} & \emph{Reinsertion of customer visits, simplified cheapest insertion, dynamic programming labeling algorithm} & - & Maximum number of iterations and time limit & Using repair operators \\
\hline
\citep{akbay2022variable} & VNS & \emph{Savings-based initial solution construction algorithm} & \emph{Random cyclic exchange, Random sequence relocation, Random customer removal, Random route removal, Close satellite} & \emph{Greedy customer insertion, Greedy customer insertion with noise, Best customer insertion, Greedy CS insertion} & \emph{exchange(1,1), shift(1,0), swap, relocation, 2-opt, and CS-reinsertion} & Time limit &  Using repair operators and penalizing the infeasible solutions \\
\hline
\citep{akbay2023application} & Adapt-CMSA & Randomly-generated & Probabilistic Clark and Wright savings algorithm and insertion algorithm  & \emph{CS-insertion}  &  \emph{relocation, swap, 2-opt, exchange (1,1), shift (1,0)} and set-covering-based MILP  &Time limit & Generating feasible solutions by MILP and using repair operators\\
\hline
\end{tabular}
\end{adjustbox}
\end{sidewaystable}

Also, Table \ref{meta-problems-table} highlights the benchmarks (dataset) and size of solved instances by each proposed metaheuristic. Several methods, such as those in \cite{breunig2017two,breunig2019electric,zijlstra2021integrating}, handle relatively large datasets with customer sizes ranging from 21 to 200. In contrast, methods like the GA \citep{agardi2019two} and VNS \citep{akbay2022variable} focus on smaller instance sizes, such as 48 to 96 and 5 to 100 customers, respectively. These differences reflect the trade-off between computational complexity and solution scalability inherent to different metaheuristic frameworks. Also, some studies, such as \citep{agardi2019two}, rely on randomly generated datasets, potentially limiting direct comparisons with other methods. However, most (recent) papers addressed the 2E-EVRP benchmarks generated by \cite{breunig2019electric}, enabling method verification and comparative studies. Additionally, only one metaheuristic is developed for each 2E-EVRPTW-BSS and 2E-EVRP-SPD variants, showing a promising gap for future research.

\begin{table}[ht!]
    \centering
    \caption{Problems and size of instances solved by the proposed metaheuristics in the literature of 2E-EVRP}
    \label{meta-problems-table}
    \begin{adjustbox}{max width=\textwidth}
    \begin{tabular}{|l|l|l|l|l|}
    \hline
       \textbf{Ref.}  &\textbf{Method} & \textbf{Problem} & \textbf{Dataset} & \textbf{Size of instances} \\
       \hline
\cite{breunig2017two} & LNS & 2E-EVRP & 2E-VRP \citep{hemmelmayr2012adaptive} & $|C|=\{100,200\}$ \\
\citep{agardi2019two}& GA & 2E-EVRP & Randomly-generated & $|C|=\{48,96\}$\\
 \citep{breunig2019electric} & LNS & 2E-EVRP &  2E-VRP \citep{perboli2011two,hemmelmayr2012adaptive,baldacci2013exact}  & $|C|=[21,200]$\\
\citep{affi2020general}& GVNS & 2E-EVRP & 2E-EVRP \citep{breunig2019electric} & $|C|=100$  \\
\citep{wang2021two}& VNS &2E-EVRPTW-BSS& 2E-EVRP \citep{breunig2019electric} & $|C|=[5,200]$\\
\citep{zijlstra2021integrating} & LNS & 2E-EVRP & 2E-EVRP \citep{breunig2019electric} & $|C|=[21,200]$ \\
\citep{akbay2022variable} & VNS & 2E-EVRP-TW & EVRP-TW \citep{schneider2014electric} & $|C|=[5,100]$\\
\citep{akbay2023application}& Adapt-CMSA &  2E-EVRP-SPD &  2E-EVRP \citep{breunig2019electric} & $|C|=[21,200]$ \\
\hline
\end{tabular}
\end{adjustbox}
\end{table}

\subsection{Matheuristic}\label{matheuristic-method}
Matheuristics are hybrid optimization approaches that combine the strengths of mathematical programming and heuristic techniques to solve complex optimization problems efficiently. They are particularly effective in large-scale, real-world challenges such as vehicle routing, scheduling, and supply chain management, where finding exact solutions is computationally prohibitive \citep{archetti2014survey,kramer2015matheuristic}. Matheuristics is crucial because it balances solution quality and computational efficiency, making it a critical tool for decision-makers seeking near-optimal solutions within reasonable timeframes. Additionally, their adaptability to various problem structures and constraints enhances their applicability across diverse domains, from logistics and manufacturing to energy and healthcare. Only one work has developed a metaheuristic for the 2E-EVRPs.

\cite{jie2019two} developed an integrated CG and adaptive LNS (CG-ALNS) for the 2E-EVRP-BSS. It combines the CG method for the first echelon with the ALNS for the second. This hybrid method effectively handles the complexity of the two-echelon structure by separately optimizing each echelon while maintaining their interdependencies. The process begins with the second echelon, where an ALNS approach is used to solve the multi-depot electric vehicle routing problem (MDEVRP). This step assigns customers to satellites and determines second-echelon routes by treating satellites as depots. ALNS employs various removal and insertion operations to explore the solution space and optimize customer-to-satellite assignments. Once the second-echelon solution is obtained, the first echelon is modeled as a split delivery electric vehicle routing problem (SDEVRP). The demand of each satellite, determined from the second echelon, is split across multiple large EVs. This echelon is optimized using the CG method, where a master problem determines the best routes and a subproblem identifies new routes with reduced costs. The subproblem is solved using a dynamic programming algorithm, which constructs and evaluates feasible routes for first-echelon EVs.

The CG-ALNS algorithm iterates between these two echelons, refining solutions through repeated optimization. A key feature of the approach is its ability to temporarily accommodate infeasible solutions, using penalized objective functions to guide the search toward feasibility. Additionally, a label-based BSS optimization algorithm dynamically adjusts the use of BSSs in both echelons to further minimize costs. To escape local optima and ensure diverse solution exploration., simulated annealing-based acceptance criteria are employed. The method integrates preprocessing techniques, such as arc elimination rules, to enhance computational efficiency by reducing the solution space. The CG-ALNS framework performs well for large-scale 2E-EVRP-BSS instances. Those instances were generated based on 2E-VRP instances \citep{crainic2010two,perboli2011two,hemmelmayr2012adaptive,baldacci2013exact}. They successfully solved the 2E-EVRP-BSS instances with customer sizes up to 200 while comparing the method's performance on 2E-VRP benchmarks. Its hybrid optimization structure effectively balances solution quality and computational efficiency.

\section{Study gaps \& future research directions}\label{gaps-future}
This section highlights significant literature trends and identifies gaps or under-researched areas. Based on the gaps and challenges identified, we suggest potential directions for future research.
These directions could involve new problem variants and extensions (Section \ref{problem-extensions-future}), practical applications with emerging technologies (Section \ref{real-apps}), and methodological improvements (Section \ref{novel-methods}), explained as follows.

\subsection{Problem extensions and new variants}\label{problem-extensions-future}
According to the previous sections and after a deep analysis of the papers on 2E-EVRP, several study gaps regarding problem extensions could be observed. First, no paper has tacked the delivery tardiness or late delivery in their model, although some papers have assumed time window-related constraints \citep{wang2019two,caggiani2021green,wang2021two,akbay2022variable}. Delivery tardiness, resulting in penalties or reduced customer satisfaction, is a critical challenge in two-echelon systems. Future research could incorporate on-time delivery with predictive analytics to mitigate tardiness while maintaining cost efficiency. Further, including stochastic travel times \citep{shao2017electric} or waiting times at CSs \citep{keskin2019electric} and demand variability in 2E-EVRP models could provide more robust solutions to minimize delays. Although EVs reduce emissions compared to ICVs, their lifecycle environmental impact, including electricity sourced from non-renewable resources, battery manufacturing, and disposal, remains a concern. Studies could explore models incorporating renewable energy sources for charging infrastructure, optimizing routes to minimize energy usage, and integrating carbon credits or penalties into objective functions \citep{li2020electric}. A comprehensive analysis of the trade-offs between environmental and operational costs in multi-echelon systems is also warranted.

High service levels, such as timely delivery and customer satisfaction, are pivotal in competitive logistics environments \citep{qin2019vehicle,moradi2025prize,boroujeni2025last}. Research could explore customer-oriented 2E-EVRP models that dynamically adjust service levels based on customer preferences, urgency, or loyalty tiers. Furthermore, hybrid optimization methods combining machine learning (ML) with mathematical modeling could be employed to maximize on-time deliveries and overall service quality. In addition, balancing trade-offs between conflicting objectives—such as minimizing costs, reducing environmental impact, and improving service quality—is a complex challenge in 2E-EVRPs. Future studies could focus on developing advanced multi-objective algorithms, such as decomposition-based approaches or evolutionary algorithms, to handle these trade-offs effectively. Additionally, sensitivity analyses on the weights assigned to different objectives can help decision-makers prioritize based on their strategic goals.

Moreover, based on the collected papers, we noticed no 2E-EVRP with multiple depots. Existing research often assumes a single central depot for the first echelon. However, real-world logistics systems frequently involve multiple depots to serve dispersed urban areas better \citep{wang2021solving,wang2023collaborative}. Expanding the 2E-EVRP to include multiple depot locations and their interactions could improve scalability and efficiency. Challenges such as depot selection, synchronization, and load balancing across depots warrant further investigation. Additionally, allowing split deliveries, where customer orders can be fulfilled in multiple trips, can improve resource utilization and reduce delivery costs \citep{mehlawat2019hybrid}. Future research could explore the computational complexity introduced by split deliveries and develop algorithms to ensure their practical applicability. Studies could also analyze the impact of split deliveries on customer satisfaction and operational efficiency.

Furthermore, 2E-EVRP should include heterogeneous fleets with varying capacities, costs, and energy consumption rates in both echelons, successfully integrated with EVRP \citep{park2020electric,wang2025heterogeneous} and 2E-VRP \citep{kancharla2019multi} in recent years. This extension would reflect real-world scenarios where multiple types of vehicles (e.g., small EVs, e-bikes, or hybrid vehicles) are deployed. Also, urban logistics are heavily influenced by time-varying traffic conditions, making static travel time assumptions unrealistic. Incorporating time-dependent travel times into 2E-EVRP models could significantly improve their applicability in real-world scenarios. Research could focus on developing efficient algorithms that account for peak and off-peak traffic patterns, leveraging historical and real-time data for more accurate routing decisions. Also, the current 2E-EVRPs could be extended to cover multi-period planning horizons, addressing inventory replenishment at satellites and dynamic route planning over multiple days. This detailed discussion of gaps and future research opportunities could address theoretical and practical challenges in sustainable logistics.

\subsection{Real-life applications}\label{real-apps}
Based on the study gas found in the literature on 2E-EVRPs, several directions for future research regarding real-life applications could be discussed. Aligning 2E-EVRP models with urban policies and regulations, such as low-emission zones or traffic restrictions, is a growing area of importance. Collaborative studies between logistics researchers and urban planners could propose policy recommendations and optimized logistics models tailored to urban constraints. Also, adapting 2E-EVRP models for emergency or disaster relief logistics could have significant societal benefits. These systems could deliver medical supplies, food, or other critical resources under time-sensitive and resource-constrained conditions. Research could explore the integration of resilience metrics and robust optimization to account for uncertainties in such scenarios. Also, in addition to small EVs in the second echelon, 2E-EVRPs could include delivery robots \citep{srinivas2022autonomous,moradi2024covering}, drones \citep{moradi2024urban} or autonomous vehicles \citep{moradi2023last} in the second echelon (while considering recharging stations for them on routes) for faster and more cost-efficient last-mile deliveries, especially in areas with challenging accessibility.

Moreover, the 2E-EVRPs could explicitly model urban congestion effects, incorporating regulations like congestion charges or restricted delivery hours. These constraints can guide the development of more practical and compliant solutions. For another research direction, the 2E-EVRP may be combined with inventory routing problems (IRP) \citep{moin2007inventory,malladi2018sustainability}, focusing on inventory levels at satellites and the central depot. This integration can optimize the flow of goods throughout the supply chain. Also, resilience metrics could be added to evaluate the system’s ability to recover from disruptions (e.g., traffic jams, vehicle breakdowns, or CS outages). Robust optimization methods could be employed to address such uncertainties. The 2E-EVRPs could also address the joint optimization of routing decisions and charging infrastructure deployment, considering strategic placement of CSs to balance setup costs and operational efficiency, similar to CS location and EVRP integration \citep{hung2022novel,wang2022electric}.

In addition, integrating parcel lockers \citep{tsai2024last} into the 2E-EVRP offers a promising avenue to enhance efficiency, sustainability, and customer satisfaction in last-mile logistics. Previously, \cite{enthoven2020two} integrated parcel lockers with 2E-VRP, and future work could be inspired by it. Also, combining parcel lockers (also known as pickup stations) into delivery routing problems has been reviewed under a general ``set covering routing problems" topic by \cite{moradi2024set}, which could be extended to 2E-EVRP with lockers. Parcel lockers can serve as intermediate delivery points in the second echelon, reducing the need for direct home deliveries, saving energy, and minimizing operational costs. Future extensions could focus on optimizing locker locations, incorporating dynamic locker assignments, considering capacity and time-window constraints, and exploring hybrid systems that combine locker deliveries for non-priority parcels with home delivery for high-priority items. Additionally, integrating locker placement with EV charging infrastructure, leveraging internet-of-things (IoT)-enabled smart lockers, and incorporating multi-objective optimization to balance cost, energy, and service level trade-offs are key areas for further research. These advancements can simplify last-mile delivery, improve scalability, and promote sustainable urban logistics.

Furthermore, incorporating same-day delivery services (similar to the prize-collecting concept in VRPs \citep{li2016two,trachanatzi2020firefly,moradi2024electric,moradi2025prize,boroujeni2025last}) into the 2E-EVRP addresses the growing demand for rapid and reliable last-mile logistics. Same-day delivery introduces stricter time constraints, requiring dynamic and time-sensitive routing strategies prioritizing urgency while maintaining efficiency. Future research could explore integrating real-time traffic data, dynamic customer order updates, and time-window constraints to ensure timely deliveries. Additionally, hybrid models could balance same-day and regular deliveries by segmenting customers based on service levels. At the same time, optimization techniques could minimize trade-offs between cost, energy consumption, and service speed. These extensions would enhance the 2E-EVRP's applicability in high-demand urban logistics, enabling companies to meet customer expectations and maintain competitive advantages.

\subsection{Novel solution approaches}\label{novel-methods}
Analysis of the 2E-EVRP papers shows a promising opportunity for future research on novel methodologies. In the past decade, there has been a surge in using artificial intelligence (AI) to improve solution quality, reduce time complexity, and solve combinatorial optimization problems, including VRPs and EVRPs. Building on this momentum, leveraging ML has the ability to greatly improve efficiency and scalability of solving 2E-EVRP. Predicting parcel demand accurately is crucial for balancing vehicle utilization and minimizing delivery costs. ML models, such as long short-term memory (LSTM), graph neural networks (GNN), transformers, gradient boosting machines (GBM), and pointer networks, have been used to analyze historical data, seasonal trends, and external factors like weather or economic activity~\citep{10639499}. These demand forecasts can seamlessly integrate into routing models like 2E-EVRPs, allowing for more informed and adaptive planning.

Also, clustering algorithms are widely combined with optimization algorithms to solve VRP~\citep{prodhon2014survey} and EVRP~\citep{ye2022electric,10139512} to some extent. Traditional methods like K-Means and Hierarchical Clustering balance workloads effectively but struggle in dynamic urban environments. Dynamic approaches, such as density-based spatial clustering of applications with noise (DBSCAN) and fuzzy clustering, adapt to real-time data, while multi-criteria clustering incorporates factors like time windows, parcel weights, and vehicle battery levels. These methods could integrate with routing algorithms like 2E-EVRPs to improve efficiency and determine optimal locations for satellites and routing decisions. However, challenges like real-time computational demands and seamless integration with routing models persist. The advent of pointer networks (Ptr-Net)~\citep{vinyals2015pointer}  introduced a learning-based approach to solving combinatorial optimization problems. Building on this foundation, reinforcement learning (RL) has been combined with GNNs, Ptr-Nets, and traditional heuristic and metaheuristic algorithms to effectively address problems such as VRP and its variants~\citep{MAZYAVKINA2021105400,JMLR:v24:21-0449,10379532}. These approaches enable more flexible and scalable solutions than traditional optimization methods, particularly for high-dimensional, dynamic, and real-time problem instances. Recently, RL has also been applied to EVRP~\citep{10.3389/fdata.2021.586481,9520134,BASSO2022102496,TANG2023121711}, demonstrating its potential in dynamic routing scenarios. By leveraging the adaptability of RL, these methods can accommodate evolving constraints such as charging requirements, real-time traffic, and energy consumption, providing a robust framework for optimizing delivery operations in 2E-EVRPs. 

Moreover, dynamic routing and real-time decision-making have become critical components of modern solutions to EVRP, particularly in addressing the uncertainties of urban delivery operations. These methods leverage real-time data, such as traffic conditions, customer demand, and vehicle battery levels, to dynamically adjust routes and depot assignments, enhancing operational flexibility and efficiency. Recent advancements have integrated real-time traffic updates and IoT-enabled sensors to optimize routing on the fly, reducing delays and energy consumption while maintaining customer satisfaction~\citep{9713755,10239234,DASTPAK2024104411,9481864}. The ideas in these papers could be extended to 2E-EVRPs as a promising future research area. Moreover, event-driven optimization frameworks, which trigger rerouting in response to disruptions like vehicle breakdowns or weather changes, have shown promise in improving system resilience. Dynamic approaches often incorporate machine learning algorithms, such as RL, to adapt to evolving constraints, enabling better handling of real-world complexities~\citep{9481864,BASSO2022102496,Ouertani2024} in routing problems like 2E-EVRPs. 

Simulation-based approaches have also emerged as valuable tools for studying and solving the EVRP. By simulating real-world scenarios, researchers can evaluate the performance of routing algorithms under various conditions, such as demand uncertainty, vehicle energy constraints, and traffic dynamics~\citep{ye2022electric}. Agent-based models are highly effective in capturing the interactions between vehicles, depots, and customers, providing insights into system performance and bottlenecks~\citep{10286625}. Discrete-event simulation is another widely used technique, enabling the analysis of time-sensitive decisions, such as scheduling and dynamic rerouting~\citep{KESKIN2021105060}. Recent studies have integrated simulation with optimization algorithms, such as metaheuristics~\citep{KESKIN2021105060}, to test and refine their solutions for EVRP. These simulation approaches not only support the development of robust algorithms but also provide a practical means to validate their applicability in real-world logistics networks. Recently, multi-method simulation frameworks~\citep{AFTABI2025125681} have proven effective in capturing the dynamic nature and stochasticity of the system across different levels, providing a comprehensive approach to analyzing complex interactions in EVRP. This idea could be effectively applied in 2E-EVRP to capture stochastic parameters and efficiently tackle the complexity of the problem. Future research could further explore and enhance these frameworks, integrating them with optimization and ML techniques to develop robust, scalable, and adaptable solutions for real-world applications.

\section{Conclusion}\label{conclusion}
Multi-echelon parcel delivery systems, especially those utilizing EVs, have become essential for managing urban logistics complexity while promoting sustainability. These systems, structured in two stages—larger vehicles transporting goods to satellites and smaller EVs managing last-mile deliveries—enable route optimization, reduced emissions, and enhanced service reliability. The 2E-EVRP, an extension of the traditional 2E-VRP, incorporates EV-specific challenges such as battery constraints and recharging station placements, addressing environmental impacts, urban congestion, and e-commerce-driven delivery demands. Despite its advantages in reducing costs, energy consumption, and greenhouse gas emissions, significant modeling challenges persist due to multi-echelon structures and EV-specific limitations. Unlike existing reviews, this paper systematically analyzes 2E-EVRP literature. The study addresses three key research questions: identifying optimization problems and variants, investigating mathematical models and solution techniques, and exploring research gaps and future directions on the 2E-EVRPs from the OR perspective. A systematic search strategy utilized targeted keywords and stringent inclusion/exclusion criteria, narrowing the selection to 12 relevant studies.

Based on collected papers, this work presents a classification scheme based on problem variants, objectives, constraints, and solution methods. This scheme facilitates structured analysis and highlights patterns, challenges, and gaps. This review identifies key gaps in the literature on 2E-EVRP and proposes several future research directions. Problem extensions could include delivery tardiness, environmental cost trade-offs, multi-objective optimization, multiple depots, split deliveries, heterogeneous fleets, and time-dependent travel conditions. Practical applications could align 2E-EVRP models with urban policies, integrate parcel lockers and same-day delivery services, adapt to emergency logistics, and include advanced technologies like autonomous vehicles, drones, and IoT-enabled systems. Methodological advancements suggest leveraging artificial intelligence, machine learning, reinforcement learning, and simulation-based approaches to address dynamic routing, real-time decision-making, and stochastic parameters. These directions aim to enhance 2E-EVRP’s applicability in sustainable and efficient urban logistics while addressing theoretical and practical challenges in future studies.

\section*{CRediT authorship contribution statement}
\textbf{Nima Moradi \& Niloufar Mirzavand Boroujeni:} Conceptualization, Data curation, Formal analysis, Investigation, Methodology, Software, Supervision, Validation, Visualization, Roles/Writing - original draft, and Writing - review \& editing. \textbf{Navid Aftabi:} Formal analysis, Resources, Software, Validation, Writing - review \& editing.
\textbf{Amin Aslani:} Formal analysis, Resources, Validation, Writing - review \& editing.

\section*{Acknowledgment}
N/A

\section*{Conflict of interest} 
The authors have no conflicts of interest to disclose.

\section*{Funding} 
The authors received no financial support for this paper's research, authorship, and publication.

\section*{Data Availability Statement}
Data, models, and codes are available upon request.

%% The Appendices part is started with the command \appendix;
%% appendix sections are then done as normal sections
%% \appendix

%% \section{}
%% \label{}

%% If you have bibdatabase file and want bibtex to generate the
%% bibitems, please use
%%
\bibliographystyle{elsarticle-harv} 
\bibliography{references.bib}

\begin{thebibliography}{96}
\expandafter\ifx\csname natexlab\endcsname\relax\def\natexlab#1{#1}\fi
\providecommand{\url}[1]{\texttt{#1}}
\providecommand{\href}[2]{#2}
\providecommand{\path}[1]{#1}
\providecommand{\DOIprefix}{doi:}
\providecommand{\ArXivprefix}{arXiv:}
\providecommand{\URLprefix}{URL: }
\providecommand{\Pubmedprefix}{pmid:}
\providecommand{\doi}[1]{\href{http://dx.doi.org/#1}{\path{#1}}}
\providecommand{\Pubmed}[1]{\href{pmid:#1}{\path{#1}}}
\providecommand{\bibinfo}[2]{#2}
\ifx\xfnm\relax \def\xfnm[#1]{\unskip,\space#1}\fi
%Type = Inproceedings
\bibitem[{Affi(2020)}]{affi2020general}
\bibinfo{author}{Affi, M.}, \bibinfo{year}{2020}.
\newblock \bibinfo{title}{General variable neighborhood search approach for solving the electric two-echelon vehicle routing problem}, in: \bibinfo{booktitle}{2020 International Multi-Conference on:“Organization of Knowledge and Advanced Technologies”(OCTA)}, \bibinfo{organization}{IEEE}. pp. \bibinfo{pages}{1--3}.
%Type = Article
\bibitem[{Aftabi et~al.(2025)Aftabi, Moradi, Mahroo and Kianfar}]{AFTABI2025125681}
\bibinfo{author}{Aftabi, N.}, \bibinfo{author}{Moradi, N.}, \bibinfo{author}{Mahroo, F.}, \bibinfo{author}{Kianfar, F.}, \bibinfo{year}{2025}.
\newblock \bibinfo{title}{Sd-abm-ism: An integrated system dynamics and agent-based modeling framework for information security management in complex information systems with multi-actor threat dynamics}.
\newblock \bibinfo{journal}{Expert Systems with Applications} \bibinfo{volume}{263}, \bibinfo{pages}{125681}.
\newblock \URLprefix \url{https://www.sciencedirect.com/science/article/pii/S095741742402548X}, \DOIprefix\doi{https://doi.org/10.1016/j.eswa.2024.125681}.
%Type = Article
\bibitem[{Ag{\'a}rdi et~al.(2019)Ag{\'a}rdi, Kov{\'a}cs and B{\'a}nyai}]{agardi2019two}
\bibinfo{author}{Ag{\'a}rdi, A.}, \bibinfo{author}{Kov{\'a}cs, L.}, \bibinfo{author}{B{\'a}nyai, T.}, \bibinfo{year}{2019}.
\newblock \bibinfo{title}{Two-echelon vehicle routing problem with recharge stations}.
\newblock \bibinfo{journal}{Transport and Telecommunication Journal} \bibinfo{volume}{20}, \bibinfo{pages}{305--317}.
%Type = Inproceedings
\bibitem[{Akbay et~al.(2023)Akbay, Kalayci and Blum}]{akbay2023application}
\bibinfo{author}{Akbay, M.A.}, \bibinfo{author}{Kalayci, C.B.}, \bibinfo{author}{Blum, C.}, \bibinfo{year}{2023}.
\newblock \bibinfo{title}{Application of adapt-cmsa to the two-echelon electric vehicle routing problem with simultaneous pickup and deliveries}, in: \bibinfo{booktitle}{European Conference on Evolutionary Computation in Combinatorial Optimization (Part of EvoStar)}, \bibinfo{organization}{Springer}. pp. \bibinfo{pages}{16--33}.
%Type = Article
\bibitem[{Akbay et~al.(2022)Akbay, Kalayci, Blum and Polat}]{akbay2022variable}
\bibinfo{author}{Akbay, M.A.}, \bibinfo{author}{Kalayci, C.B.}, \bibinfo{author}{Blum, C.}, \bibinfo{author}{Polat, O.}, \bibinfo{year}{2022}.
\newblock \bibinfo{title}{Variable neighborhood search for the two-echelon electric vehicle routing problem with time windows}.
\newblock \bibinfo{journal}{Applied Sciences} \bibinfo{volume}{12}, \bibinfo{pages}{1014}.
%Type = Inproceedings
\bibitem[{Almasan et~al.(2021)Almasan, Suárez-Varela, Wu, Xiao, Barlet-Ros and Cabellos-Aparicio}]{9481864}
\bibinfo{author}{Almasan, P.}, \bibinfo{author}{Suárez-Varela, J.}, \bibinfo{author}{Wu, B.}, \bibinfo{author}{Xiao, S.}, \bibinfo{author}{Barlet-Ros, P.}, \bibinfo{author}{Cabellos-Aparicio, A.}, \bibinfo{year}{2021}.
\newblock \bibinfo{title}{Towards real-time routing optimization with deep reinforcement learning: Open challenges}, in: \bibinfo{booktitle}{2021 IEEE 22nd International Conference on High Performance Switching and Routing (HPSR)}, pp. \bibinfo{pages}{1--6}.
\newblock \DOIprefix\doi{10.1109/HPSR52026.2021.9481864}.
%Type = Article
\bibitem[{Archetti and Speranza(2014)}]{archetti2014survey}
\bibinfo{author}{Archetti, C.}, \bibinfo{author}{Speranza, M.G.}, \bibinfo{year}{2014}.
\newblock \bibinfo{title}{A survey on matheuristics for routing problems}.
\newblock \bibinfo{journal}{EURO Journal on Computational Optimization} \bibinfo{volume}{2}, \bibinfo{pages}{223--246}.
%Type = Article
\bibitem[{Baldacci et~al.(2013)Baldacci, Mingozzi, Roberti and Calvo}]{baldacci2013exact}
\bibinfo{author}{Baldacci, R.}, \bibinfo{author}{Mingozzi, A.}, \bibinfo{author}{Roberti, R.}, \bibinfo{author}{Calvo, R.W.}, \bibinfo{year}{2013}.
\newblock \bibinfo{title}{An exact algorithm for the two-echelon capacitated vehicle routing problem}.
\newblock \bibinfo{journal}{Operations research} \bibinfo{volume}{61}, \bibinfo{pages}{298--314}.
%Type = Article
\bibitem[{Basso et~al.(2022)Basso, Kulcsár, Sanchez-Diaz and Qu}]{BASSO2022102496}
\bibinfo{author}{Basso, R.}, \bibinfo{author}{Kulcsár, B.}, \bibinfo{author}{Sanchez-Diaz, I.}, \bibinfo{author}{Qu, X.}, \bibinfo{year}{2022}.
\newblock \bibinfo{title}{Dynamic stochastic electric vehicle routing with safe reinforcement learning}.
\newblock \bibinfo{journal}{Transportation Research Part E: Logistics and Transportation Review} \bibinfo{volume}{157}, \bibinfo{pages}{102496}.
\newblock \URLprefix \url{https://www.sciencedirect.com/science/article/pii/S1366554521002581}, \DOIprefix\doi{https://doi.org/10.1016/j.tre.2021.102496}.
%Type = Article
\bibitem[{Benhassine and Shan(2023)}]{benhassine2023optimization}
\bibinfo{author}{Benhassine, A.}, \bibinfo{author}{Shan, B.}, \bibinfo{year}{2023}.
\newblock \bibinfo{title}{Optimization of cost and carbon emissions in a multi-echelon distribution network} .
%Type = Article
\bibitem[{Bogyrbayeva et~al.(2024)Bogyrbayeva, Meraliyev, Mustakhov and Dauletbayev}]{10379532}
\bibinfo{author}{Bogyrbayeva, A.}, \bibinfo{author}{Meraliyev, M.}, \bibinfo{author}{Mustakhov, T.}, \bibinfo{author}{Dauletbayev, B.}, \bibinfo{year}{2024}.
\newblock \bibinfo{title}{Machine learning to solve vehicle routing problems: A survey}.
\newblock \bibinfo{journal}{IEEE Transactions on Intelligent Transportation Systems} \bibinfo{volume}{25}, \bibinfo{pages}{4754--4772}.
\newblock \DOIprefix\doi{10.1109/TITS.2023.3334976}.
%Type = Article
\bibitem[{Boroujeni et~al.(2025)Boroujeni, Moradi, Jamalzadeh and Boroujeni}]{boroujeni2025last}
\bibinfo{author}{Boroujeni, N.M.}, \bibinfo{author}{Moradi, N.}, \bibinfo{author}{Jamalzadeh, S.}, \bibinfo{author}{Boroujeni, N.M.}, \bibinfo{year}{2025}.
\newblock \bibinfo{title}{Last-mile delivery optimization: Leveraging electric vehicles and parcel lockers for prime customer service}.
\newblock \bibinfo{journal}{Computers \& Industrial Engineering} \bibinfo{volume}{203}, \bibinfo{pages}{110991}.
%Type = Article
\bibitem[{Breunig et~al.(2019)Breunig, Baldacci, Hartl and Vidal}]{breunig2019electric}
\bibinfo{author}{Breunig, U.}, \bibinfo{author}{Baldacci, R.}, \bibinfo{author}{Hartl, R.F.}, \bibinfo{author}{Vidal, T.}, \bibinfo{year}{2019}.
\newblock \bibinfo{title}{The electric two-echelon vehicle routing problem}.
\newblock \bibinfo{journal}{Computers \& Operations Research} \bibinfo{volume}{103}, \bibinfo{pages}{198--210}.
%Type = Inproceedings
\bibitem[{Breunig et~al.(2017)Breunig, Hartl and Vidal}]{breunig2017two}
\bibinfo{author}{Breunig, U.}, \bibinfo{author}{Hartl, R.F.}, \bibinfo{author}{Vidal, T.}, \bibinfo{year}{2017}.
\newblock \bibinfo{title}{The two-echelon vehicle routing problem with electric vehicles}, in: \bibinfo{booktitle}{XLIX Simpósio Brasileiro de Pesquisa Operacional}.
%Type = Article
\bibitem[{Caggiani et~al.(2021)Caggiani, Colovic, Prencipe and Ottomanelli}]{caggiani2021green}
\bibinfo{author}{Caggiani, L.}, \bibinfo{author}{Colovic, A.}, \bibinfo{author}{Prencipe, L.P.}, \bibinfo{author}{Ottomanelli, M.}, \bibinfo{year}{2021}.
\newblock \bibinfo{title}{A green logistics solution for last-mile deliveries considering e-vans and e-cargo bikes}.
\newblock \bibinfo{journal}{Transportation Research Procedia} \bibinfo{volume}{52}, \bibinfo{pages}{75--82}.
%Type = Article
\bibitem[{Cappart et~al.(2023)Cappart, Ch{\'e}telat, Khalil, Lodi, Morris and Veli{\v{c}}kovi{\'c}}]{JMLR:v24:21-0449}
\bibinfo{author}{Cappart, Q.}, \bibinfo{author}{Ch{\'e}telat, D.}, \bibinfo{author}{Khalil, E.B.}, \bibinfo{author}{Lodi, A.}, \bibinfo{author}{Morris, C.}, \bibinfo{author}{Veli{\v{c}}kovi{\'c}, P.}, \bibinfo{year}{2023}.
\newblock \bibinfo{title}{Combinatorial optimization and reasoning with graph neural networks}.
\newblock \bibinfo{journal}{Journal of Machine Learning Research} \bibinfo{volume}{24}, \bibinfo{pages}{1--61}.
\newblock \URLprefix \url{http://jmlr.org/papers/v24/21-0449.html}.
%Type = Book
\bibitem[{Crainic et~al.(2008)Crainic, Mancini, Perboli, Tadei et~al.}]{crainic2008clustering}
\bibinfo{author}{Crainic, T.G.}, \bibinfo{author}{Mancini, S.}, \bibinfo{author}{Perboli, G.}, \bibinfo{author}{Tadei, R.}, et~al., \bibinfo{year}{2008}.
\newblock \bibinfo{title}{Clustering-based heuristics for the two-echelon vehicle routing problem}. volume~\bibinfo{volume}{46}.
\newblock \bibinfo{publisher}{CIRRELT Montr{\'e}al}.
%Type = Article
\bibitem[{Crainic et~al.(2010)Crainic, Perboli, Mancini and Tadei}]{crainic2010two}
\bibinfo{author}{Crainic, T.G.}, \bibinfo{author}{Perboli, G.}, \bibinfo{author}{Mancini, S.}, \bibinfo{author}{Tadei, R.}, \bibinfo{year}{2010}.
\newblock \bibinfo{title}{Two-echelon vehicle routing problem: a satellite location analysis}.
\newblock \bibinfo{journal}{Procedia-Social and Behavioral Sciences} \bibinfo{volume}{2}, \bibinfo{pages}{5944--5955}.
%Type = Article
\bibitem[{Crainic et~al.(2009)Crainic, Ricciardi and Storchi}]{crainic2009models}
\bibinfo{author}{Crainic, T.G.}, \bibinfo{author}{Ricciardi, N.}, \bibinfo{author}{Storchi, G.}, \bibinfo{year}{2009}.
\newblock \bibinfo{title}{Models for evaluating and planning city logistics systems}.
\newblock \bibinfo{journal}{Transportation science} \bibinfo{volume}{43}, \bibinfo{pages}{432--454}.
%Type = Article
\bibitem[{Cuda et~al.(2015)Cuda, Guastaroba and Speranza}]{cuda2015survey}
\bibinfo{author}{Cuda, R.}, \bibinfo{author}{Guastaroba, G.}, \bibinfo{author}{Speranza, M.G.}, \bibinfo{year}{2015}.
\newblock \bibinfo{title}{A survey on two-echelon routing problems}.
\newblock \bibinfo{journal}{Computers \& Operations Research} \bibinfo{volume}{55}, \bibinfo{pages}{185--199}.
%Type = Article
\bibitem[{Dastpak et~al.(2024)Dastpak, Errico, Jabali and Malucelli}]{DASTPAK2024104411}
\bibinfo{author}{Dastpak, M.}, \bibinfo{author}{Errico, F.}, \bibinfo{author}{Jabali, O.}, \bibinfo{author}{Malucelli, F.}, \bibinfo{year}{2024}.
\newblock \bibinfo{title}{Dynamic routing for the electric vehicle shortest path problem with charging station occupancy information}.
\newblock \bibinfo{journal}{Transportation Research Part C: Emerging Technologies} \bibinfo{volume}{158}, \bibinfo{pages}{104411}.
\newblock \URLprefix \url{https://www.sciencedirect.com/science/article/pii/S0968090X23004011}, \DOIprefix\doi{https://doi.org/10.1016/j.trc.2023.104411}.
%Type = Misc
\bibitem[{DHL.com, 2024()}]{dhl}
DHL.com, 2024, \bibinfo{year}{2024}.
\newblock \bibinfo{note}{Sustainability Trends in Logistics for 2024, https://www.dhl.com/discover/en-us/global-logistics-advice/sustainability-and-green-logistics/sustainability-trends-in-logistics?}
%Type = Article
\bibitem[{Dorokhova et~al.(2021)Dorokhova, Ballif and Wyrsch}]{10.3389/fdata.2021.586481}
\bibinfo{author}{Dorokhova, M.}, \bibinfo{author}{Ballif, C.}, \bibinfo{author}{Wyrsch, N.}, \bibinfo{year}{2021}.
\newblock \bibinfo{title}{Routing of electric vehicles with intermediary charging stations: A reinforcement learning approach}.
\newblock \bibinfo{journal}{Frontiers in Big Data} \bibinfo{volume}{4}.
\newblock \URLprefix \url{https://www.frontiersin.org/journals/big-data/articles/10.3389/fdata.2021.586481}, \DOIprefix\doi{10.3389/fdata.2021.586481}.
%Type = Article
\bibitem[{Duarte et~al.(2018)Duarte, Mladenovic, S{\'a}nchez-Oro and Todosijevi{\'c}}]{duarte2018variable}
\bibinfo{author}{Duarte, A.}, \bibinfo{author}{Mladenovic, N.}, \bibinfo{author}{S{\'a}nchez-Oro, J.}, \bibinfo{author}{Todosijevi{\'c}, R.}, \bibinfo{year}{2018}.
\newblock \bibinfo{title}{Variable neighborhood descent}.
\newblock \bibinfo{journal}{Handbook of heuristics} , \bibinfo{pages}{341--367}.
%Type = Article
\bibitem[{Durach et~al.(2017)Durach, Kembro and Wieland}]{durach2017new}
\bibinfo{author}{Durach, C.F.}, \bibinfo{author}{Kembro, J.}, \bibinfo{author}{Wieland, A.}, \bibinfo{year}{2017}.
\newblock \bibinfo{title}{A new paradigm for systematic literature reviews in supply chain management}.
\newblock \bibinfo{journal}{Journal of Supply Chain Management} \bibinfo{volume}{53}, \bibinfo{pages}{67--85}.
%Type = Article
\bibitem[{Enthoven et~al.(2020)Enthoven, Jargalsaikhan, Roodbergen, Uit~het Broek and Schrotenboer}]{enthoven2020two}
\bibinfo{author}{Enthoven, D.L.}, \bibinfo{author}{Jargalsaikhan, B.}, \bibinfo{author}{Roodbergen, K.J.}, \bibinfo{author}{Uit~het Broek, M.A.}, \bibinfo{author}{Schrotenboer, A.H.}, \bibinfo{year}{2020}.
\newblock \bibinfo{title}{The two-echelon vehicle routing problem with covering options: City logistics with cargo bikes and parcel lockers}.
\newblock \bibinfo{journal}{Computers \& Operations Research} \bibinfo{volume}{118}, \bibinfo{pages}{104919}.
%Type = Article
\bibitem[{Flood(1956)}]{flood1956traveling}
\bibinfo{author}{Flood, M.M.}, \bibinfo{year}{1956}.
\newblock \bibinfo{title}{The traveling-salesman problem}.
\newblock \bibinfo{journal}{Operations research} \bibinfo{volume}{4}, \bibinfo{pages}{61--75}.
%Type = Misc
\bibitem[{Ford.ca, 2024()}]{ford.ca}
Ford.ca, 2024, \bibinfo{year}{2024}.
\newblock \bibinfo{note}{Hardworking efficiency to power your work., https://www.ford.ca/commercial-trucks/e-transit/}.
%Type = Article
\bibitem[{Goeke(2019)}]{goeke2019vrptw}
\bibinfo{author}{Goeke, D.}, \bibinfo{year}{2019}.
\newblock \bibinfo{title}{E-vrptw instances}.
\newblock \bibinfo{journal}{Mendeley Data, v1} .
%Type = Article
\bibitem[{Gonzalez-Feliu(2011)}]{gonzalez2011two}
\bibinfo{author}{Gonzalez-Feliu, J.}, \bibinfo{year}{2011}.
\newblock \bibinfo{title}{Two-echelon freight transport optimisation: unifying concepts via a systematic review}.
\newblock \bibinfo{journal}{Working Papers on Operations Management} \bibinfo{volume}{2}, \bibinfo{pages}{18--30}.
%Type = Article
\bibitem[{Hemmelmayr et~al.(2012)Hemmelmayr, Cordeau and Crainic}]{hemmelmayr2012adaptive}
\bibinfo{author}{Hemmelmayr, V.C.}, \bibinfo{author}{Cordeau, J.F.}, \bibinfo{author}{Crainic, T.G.}, \bibinfo{year}{2012}.
\newblock \bibinfo{title}{An adaptive large neighborhood search heuristic for two-echelon vehicle routing problems arising in city logistics}.
\newblock \bibinfo{journal}{Computers \& operations research} \bibinfo{volume}{39}, \bibinfo{pages}{3215--3228}.
%Type = Article
\bibitem[{Holland(1992)}]{holland1992genetic}
\bibinfo{author}{Holland, J.H.}, \bibinfo{year}{1992}.
\newblock \bibinfo{title}{Genetic algorithms}.
\newblock \bibinfo{journal}{Scientific american} \bibinfo{volume}{267}, \bibinfo{pages}{66--73}.
%Type = Article
\bibitem[{Hung and Michailidis(2022)}]{hung2022novel}
\bibinfo{author}{Hung, Y.C.}, \bibinfo{author}{Michailidis, G.}, \bibinfo{year}{2022}.
\newblock \bibinfo{title}{A novel data-driven approach for solving the electric vehicle charging station location-routing problem}.
\newblock \bibinfo{journal}{IEEE Transactions on Intelligent Transportation Systems} \bibinfo{volume}{23}, \bibinfo{pages}{23858--23868}.
%Type = Misc
\bibitem[{IEA.org, 2021()}]{iea}
IEA.org, 2021, \bibinfo{year}{2021}.
\newblock \bibinfo{note}{Comparative life-cycle greenhouse gas emissions of a mid-size BEV and ICE vehicle, https://www.iea.org/data-and-statistics/charts/comparative-life-cycle-greenhouse-gas-emissions-of-a-mid-size-bev-and-ice-vehicle}.
%Type = Article
\bibitem[{Jie et~al.(2019)Jie, Yang, Zhang and Huang}]{jie2019two}
\bibinfo{author}{Jie, W.}, \bibinfo{author}{Yang, J.}, \bibinfo{author}{Zhang, M.}, \bibinfo{author}{Huang, Y.}, \bibinfo{year}{2019}.
\newblock \bibinfo{title}{The two-echelon capacitated electric vehicle routing problem with battery swapping stations: Formulation and efficient methodology}.
\newblock \bibinfo{journal}{European Journal of Operational Research} \bibinfo{volume}{272}, \bibinfo{pages}{879--904}.
%Type = Article
\bibitem[{Kancharla and Ramadurai(2019)}]{kancharla2019multi}
\bibinfo{author}{Kancharla, S.R.}, \bibinfo{author}{Ramadurai, G.}, \bibinfo{year}{2019}.
\newblock \bibinfo{title}{Multi-depot two-echelon fuel minimizing routing problem with heterogeneous fleets: Model and heuristic}.
\newblock \bibinfo{journal}{Networks and Spatial Economics} \bibinfo{volume}{19}, \bibinfo{pages}{969--1005}.
%Type = Article
\bibitem[{Kayvanfar et~al.(2018)Kayvanfar, S.~Sajadieh, Moattar~Husseini and Karimi}]{kayvanfar2018analysis}
\bibinfo{author}{Kayvanfar, V.}, \bibinfo{author}{S.~Sajadieh, M.}, \bibinfo{author}{Moattar~Husseini, S.}, \bibinfo{author}{Karimi, B.}, \bibinfo{year}{2018}.
\newblock \bibinfo{title}{Analysis of a multi-echelon supply chain problem using revised multi-choice goal programming approach}.
\newblock \bibinfo{journal}{Kybernetes} \bibinfo{volume}{47}, \bibinfo{pages}{118--141}.
%Type = Article
\bibitem[{Keskin and {\c{C}}atay(2016)}]{keskin2016partial}
\bibinfo{author}{Keskin, M.}, \bibinfo{author}{{\c{C}}atay, B.}, \bibinfo{year}{2016}.
\newblock \bibinfo{title}{Partial recharge strategies for the electric vehicle routing problem with time windows}.
\newblock \bibinfo{journal}{Transportation research part C: emerging technologies} \bibinfo{volume}{65}, \bibinfo{pages}{111--127}.
%Type = Article
\bibitem[{Keskin et~al.(2019)Keskin, Laporte and {\c{C}}atay}]{keskin2019electric}
\bibinfo{author}{Keskin, M.}, \bibinfo{author}{Laporte, G.}, \bibinfo{author}{{\c{C}}atay, B.}, \bibinfo{year}{2019}.
\newblock \bibinfo{title}{Electric vehicle routing problem with time-dependent waiting times at recharging stations}.
\newblock \bibinfo{journal}{Computers \& Operations Research} \bibinfo{volume}{107}, \bibinfo{pages}{77--94}.
%Type = Article
\bibitem[{Keskin et~al.(2021)Keskin, Çatay and Laporte}]{KESKIN2021105060}
\bibinfo{author}{Keskin, M.}, \bibinfo{author}{Çatay, B.}, \bibinfo{author}{Laporte, G.}, \bibinfo{year}{2021}.
\newblock \bibinfo{title}{A simulation-based heuristic for the electric vehicle routing problem with time windows and stochastic waiting times at recharging stations}.
\newblock \bibinfo{journal}{Computers \& Operations Research} \bibinfo{volume}{125}, \bibinfo{pages}{105060}.
\newblock \URLprefix \url{https://www.sciencedirect.com/science/article/pii/S0305054820301775}, \DOIprefix\doi{https://doi.org/10.1016/j.cor.2020.105060}.
%Type = Article
\bibitem[{Ko{\c{c}} et~al.(2020)Ko{\c{c}}, Laporte and T{\"u}kenmez}]{kocc2020review}
\bibinfo{author}{Ko{\c{c}}, {\c{C}}.}, \bibinfo{author}{Laporte, G.}, \bibinfo{author}{T{\"u}kenmez, {\.I}.}, \bibinfo{year}{2020}.
\newblock \bibinfo{title}{A review of vehicle routing with simultaneous pickup and delivery}.
\newblock \bibinfo{journal}{Computers \& Operations Research} \bibinfo{volume}{122}, \bibinfo{pages}{104987}.
%Type = Article
\bibitem[{Kramer et~al.(2015)Kramer, Subramanian, Vidal and Luc{\'\i}dio~dos Anjos}]{kramer2015matheuristic}
\bibinfo{author}{Kramer, R.}, \bibinfo{author}{Subramanian, A.}, \bibinfo{author}{Vidal, T.}, \bibinfo{author}{Luc{\'\i}dio~dos Anjos, F.C.}, \bibinfo{year}{2015}.
\newblock \bibinfo{title}{A matheuristic approach for the pollution-routing problem}.
\newblock \bibinfo{journal}{European Journal of Operational Research} \bibinfo{volume}{243}, \bibinfo{pages}{523--539}.
%Type = Article
\bibitem[{Li et~al.(2021)Li, Chen, Wang and Bai}]{li2021ground}
\bibinfo{author}{Li, H.}, \bibinfo{author}{Chen, J.}, \bibinfo{author}{Wang, F.}, \bibinfo{author}{Bai, M.}, \bibinfo{year}{2021}.
\newblock \bibinfo{title}{Ground-vehicle and unmanned-aerial-vehicle routing problems from two-echelon scheme perspective: A review}.
\newblock \bibinfo{journal}{European Journal of Operational Research} \bibinfo{volume}{294}, \bibinfo{pages}{1078--1095}.
%Type = Article
\bibitem[{Li et~al.(2016)Li, Zhang, Lv and Chang}]{li2016two}
\bibinfo{author}{Li, H.}, \bibinfo{author}{Zhang, L.}, \bibinfo{author}{Lv, T.}, \bibinfo{author}{Chang, X.}, \bibinfo{year}{2016}.
\newblock \bibinfo{title}{The two-echelon time-constrained vehicle routing problem in linehaul-delivery systems}.
\newblock \bibinfo{journal}{Transportation Research Part B: Methodological} \bibinfo{volume}{94}, \bibinfo{pages}{169--188}.
%Type = Article
\bibitem[{Li et~al.(2023)Li, Liu and Wang}]{10239234}
\bibinfo{author}{Li, J.}, \bibinfo{author}{Liu, R.}, \bibinfo{author}{Wang, R.}, \bibinfo{year}{2023}.
\newblock \bibinfo{title}{Elastic strategy-based adaptive genetic algorithm for solving dynamic vehicle routing problem with time windows}.
\newblock \bibinfo{journal}{IEEE Transactions on Intelligent Transportation Systems} \bibinfo{volume}{24}, \bibinfo{pages}{13930--13947}.
\newblock \DOIprefix\doi{10.1109/TITS.2023.3308593}.
%Type = Article
\bibitem[{Li et~al.(2020)Li, Wang and He}]{li2020electric}
\bibinfo{author}{Li, J.}, \bibinfo{author}{Wang, F.}, \bibinfo{author}{He, Y.}, \bibinfo{year}{2020}.
\newblock \bibinfo{title}{Electric vehicle routing problem with battery swapping considering energy consumption and carbon emissions}.
\newblock \bibinfo{journal}{Sustainability} \bibinfo{volume}{12}, \bibinfo{pages}{10537}.
%Type = Article
\bibitem[{Lin et~al.(2022)Lin, Ghaddar and Nathwani}]{9520134}
\bibinfo{author}{Lin, B.}, \bibinfo{author}{Ghaddar, B.}, \bibinfo{author}{Nathwani, J.}, \bibinfo{year}{2022}.
\newblock \bibinfo{title}{Deep reinforcement learning for the electric vehicle routing problem with time windows}.
\newblock \bibinfo{journal}{IEEE Transactions on Intelligent Transportation Systems} \bibinfo{volume}{23}, \bibinfo{pages}{11528--11538}.
\newblock \DOIprefix\doi{10.1109/TITS.2021.3105232}.
%Type = Article
\bibitem[{Malladi and Sowlati(2018)}]{malladi2018sustainability}
\bibinfo{author}{Malladi, K.T.}, \bibinfo{author}{Sowlati, T.}, \bibinfo{year}{2018}.
\newblock \bibinfo{title}{Sustainability aspects in inventory routing problem: A review of new trends in the literature}.
\newblock \bibinfo{journal}{Journal of Cleaner Production} \bibinfo{volume}{197}, \bibinfo{pages}{804--814}.
%Type = Inproceedings
\bibitem[{Mazar et~al.(2023)Mazar, Saouabe, Sahnoun and Mourtaji}]{10286625}
\bibinfo{author}{Mazar, M.}, \bibinfo{author}{Saouabe, A.}, \bibinfo{author}{Sahnoun, M.}, \bibinfo{author}{Mourtaji, I.}, \bibinfo{year}{2023}.
\newblock \bibinfo{title}{Electric vehicle route simulation: A preliminary approach}, in: \bibinfo{booktitle}{2023 International Conference on Decision Aid Sciences and Applications (DASA)}, pp. \bibinfo{pages}{641--645}.
\newblock \DOIprefix\doi{10.1109/DASA59624.2023.10286625}.
%Type = Article
\bibitem[{Mazyavkina et~al.(2021)Mazyavkina, Sviridov, Ivanov and Burnaev}]{MAZYAVKINA2021105400}
\bibinfo{author}{Mazyavkina, N.}, \bibinfo{author}{Sviridov, S.}, \bibinfo{author}{Ivanov, S.}, \bibinfo{author}{Burnaev, E.}, \bibinfo{year}{2021}.
\newblock \bibinfo{title}{Reinforcement learning for combinatorial optimization: A survey}.
\newblock \bibinfo{journal}{Computers \& Operations Research} \bibinfo{volume}{134}, \bibinfo{pages}{105400}.
\newblock \URLprefix \url{https://www.sciencedirect.com/science/article/pii/S0305054821001660}, \DOIprefix\doi{https://doi.org/10.1016/j.cor.2021.105400}.
%Type = Article
\bibitem[{Mehlawat et~al.(2019)Mehlawat, Gupta, Khaitan and Pedrycz}]{mehlawat2019hybrid}
\bibinfo{author}{Mehlawat, M.K.}, \bibinfo{author}{Gupta, P.}, \bibinfo{author}{Khaitan, A.}, \bibinfo{author}{Pedrycz, W.}, \bibinfo{year}{2019}.
\newblock \bibinfo{title}{A hybrid intelligent approach to integrated fuzzy multiple depot capacitated green vehicle routing problem with split delivery and vehicle selection}.
\newblock \bibinfo{journal}{IEEE Transactions on Fuzzy Systems} \bibinfo{volume}{28}, \bibinfo{pages}{1155--1166}.
%Type = Article
\bibitem[{Mladenovi{\'c} and Hansen(1997)}]{mladenovic1997variable}
\bibinfo{author}{Mladenovi{\'c}, N.}, \bibinfo{author}{Hansen, P.}, \bibinfo{year}{1997}.
\newblock \bibinfo{title}{Variable neighborhood search}.
\newblock \bibinfo{journal}{Computers \& operations research} \bibinfo{volume}{24}, \bibinfo{pages}{1097--1100}.
%Type = Article
\bibitem[{Mladenovi{\'c} et~al.(2012)Mladenovi{\'c}, Uro{\v{s}}evi{\'c}, Ili{\'c} et~al.}]{mladenovic2012general}
\bibinfo{author}{Mladenovi{\'c}, N.}, \bibinfo{author}{Uro{\v{s}}evi{\'c}, D.}, \bibinfo{author}{Ili{\'c}, A.}, et~al., \bibinfo{year}{2012}.
\newblock \bibinfo{title}{A general variable neighborhood search for the one-commodity pickup-and-delivery travelling salesman problem}.
\newblock \bibinfo{journal}{European Journal of Operational Research} \bibinfo{volume}{220}, \bibinfo{pages}{270--285}.
%Type = Article
\bibitem[{Moin and Salhi(2007)}]{moin2007inventory}
\bibinfo{author}{Moin, N.H.}, \bibinfo{author}{Salhi, S.}, \bibinfo{year}{2007}.
\newblock \bibinfo{title}{Inventory routing problems: a logistical overview}.
\newblock \bibinfo{journal}{Journal of the Operational Research Society} \bibinfo{volume}{58}, \bibinfo{pages}{1185--1194}.
%Type = Article
\bibitem[{Moradi and Boroujeni(2025)}]{moradi2025prize}
\bibinfo{author}{Moradi, N.}, \bibinfo{author}{Boroujeni, N.M.}, \bibinfo{year}{2025}.
\newblock \bibinfo{title}{Prize-collecting electric vehicle routing model for parcel delivery problem}.
\newblock \bibinfo{journal}{Expert Systems with Applications} \bibinfo{volume}{259}, \bibinfo{pages}{125183}.
%Type = Inproceedings
\bibitem[{Moradi et~al.(2024a)Moradi, Kayvanfar and Baldacci}]{moradi2024electric}
\bibinfo{author}{Moradi, N.}, \bibinfo{author}{Kayvanfar, V.}, \bibinfo{author}{Baldacci, R.}, \bibinfo{year}{2024}a.
\newblock \bibinfo{title}{Electric-vehicle routing problem with time windows and energy minimization: green logistics with same-day delivery approaches}, in: \bibinfo{booktitle}{2024 International Conference on Electrical, Computer and Energy Technologies (ICECET}, \bibinfo{organization}{IEEE}. pp. \bibinfo{pages}{1--6}.
%Type = Inproceedings
\bibitem[{Moradi et~al.(2024b)Moradi, Mafakheri and Wang}]{moradi2024covering}
\bibinfo{author}{Moradi, N.}, \bibinfo{author}{Mafakheri, F.}, \bibinfo{author}{Wang, C.}, \bibinfo{year}{2024}b.
\newblock \bibinfo{title}{Covering routing problem with robots and parcel lockers: A sustainable last-mile delivery approach}, in: \bibinfo{booktitle}{IISE Annual Conference. Proceedings}, \bibinfo{organization}{Institute of Industrial and Systems Engineers (IISE)}. pp. \bibinfo{pages}{1--6}.
%Type = Article
\bibitem[{Moradi et~al.(2024c)Moradi, Mafakheri and Wang}]{moradi2024set}
\bibinfo{author}{Moradi, N.}, \bibinfo{author}{Mafakheri, F.}, \bibinfo{author}{Wang, C.}, \bibinfo{year}{2024}c.
\newblock \bibinfo{title}{Set covering routing problems: A review and classification scheme}.
\newblock \bibinfo{journal}{Computers \& Industrial Engineering} , \bibinfo{pages}{110730}.
%Type = Article
\bibitem[{Moradi et~al.(2023)Moradi, Sadati and {\c{C}}atay}]{moradi2023last}
\bibinfo{author}{Moradi, N.}, \bibinfo{author}{Sadati, {\.I}.}, \bibinfo{author}{{\c{C}}atay, B.}, \bibinfo{year}{2023}.
\newblock \bibinfo{title}{Last mile delivery routing problem using autonomous electric vehicles}.
\newblock \bibinfo{journal}{Computers \& Industrial Engineering} \bibinfo{volume}{184}, \bibinfo{pages}{109552}.
%Type = Article
\bibitem[{Moradi et~al.(2024d)Moradi, Wang and Mafakheri}]{moradi2024urban}
\bibinfo{author}{Moradi, N.}, \bibinfo{author}{Wang, C.}, \bibinfo{author}{Mafakheri, F.}, \bibinfo{year}{2024}d.
\newblock \bibinfo{title}{Urban air mobility for last-mile transportation: A review}.
\newblock \bibinfo{journal}{Vehicles} \bibinfo{volume}{6}, \bibinfo{pages}{1383--1414}.
%Type = Article
\bibitem[{Nielsen et~al.(2024)Nielsen, Dahanayaka, Perera, Thibbotuwawa and Kilic}]{nielsen2024systematic}
\bibinfo{author}{Nielsen, P.}, \bibinfo{author}{Dahanayaka, M.}, \bibinfo{author}{Perera, H.N.}, \bibinfo{author}{Thibbotuwawa, A.}, \bibinfo{author}{Kilic, D.K.}, \bibinfo{year}{2024}.
\newblock \bibinfo{title}{A systematic review of vehicle routing problems and models in multi-echelon distribution networks}.
\newblock \bibinfo{journal}{Supply Chain Analytics} , \bibinfo{pages}{100072}.
%Type = Inbook
\bibitem[{Ouertani et~al.(2024)Ouertani, Ben-Romdhane and Krichen}]{Ouertani2024}
\bibinfo{author}{Ouertani, N.}, \bibinfo{author}{Ben-Romdhane, H.}, \bibinfo{author}{Krichen, S.}, \bibinfo{year}{2024}.
\newblock \bibinfo{title}{The Dynamic Vehicle Routing Problem: A Comprehensive Survey}. \bibinfo{publisher}{Springer Nature Switzerland}, \bibinfo{address}{Cham}.
\newblock pp. \bibinfo{pages}{1--36}.
%Type = Article
\bibitem[{Park and Jin(2020)}]{park2020electric}
\bibinfo{author}{Park, H.}, \bibinfo{author}{Jin, S.}, \bibinfo{year}{2020}.
\newblock \bibinfo{title}{Electric vehicle routing problem with heterogeneous vehicles and partial charge}.
\newblock \bibinfo{journal}{International Journal of Industrial Engineering and Management} \bibinfo{volume}{11}, \bibinfo{pages}{215--225}.
%Type = Article
\bibitem[{Perboli et~al.(2011)Perboli, Tadei and Vigo}]{perboli2011two}
\bibinfo{author}{Perboli, G.}, \bibinfo{author}{Tadei, R.}, \bibinfo{author}{Vigo, D.}, \bibinfo{year}{2011}.
\newblock \bibinfo{title}{The two-echelon capacitated vehicle routing problem: Models and math-based heuristics}.
\newblock \bibinfo{journal}{Transportation Science} \bibinfo{volume}{45}, \bibinfo{pages}{364--380}.
%Type = Inproceedings
\bibitem[{Phu-Ang(2023)}]{10139512}
\bibinfo{author}{Phu-Ang, A.}, \bibinfo{year}{2023}.
\newblock \bibinfo{title}{Solving the capacitated electric vehicle (ev) routing problem by the differential evolutionary algorithm with adaptive k-means}, in: \bibinfo{booktitle}{2023 Joint International Conference on Digital Arts, Media and Technology with ECTI Northern Section Conference on Electrical, Electronics, Computer and Telecommunications Engineering (ECTI DAMT \& NCON)}, pp. \bibinfo{pages}{225--228}.
\newblock \DOIprefix\doi{10.1109/ECTIDAMTNCON57770.2023.10139512}.
%Type = Article
\bibitem[{Prodhon and Prins(2014)}]{prodhon2014survey}
\bibinfo{author}{Prodhon, C.}, \bibinfo{author}{Prins, C.}, \bibinfo{year}{2014}.
\newblock \bibinfo{title}{A survey of recent research on location-routing problems}.
\newblock \bibinfo{journal}{European journal of operational research} \bibinfo{volume}{238}, \bibinfo{pages}{1--17}.
%Type = Article
\bibitem[{Prodhon and Prins(2016)}]{prodhon2016metaheuristics}
\bibinfo{author}{Prodhon, C.}, \bibinfo{author}{Prins, C.}, \bibinfo{year}{2016}.
\newblock \bibinfo{title}{Metaheuristics for vehicle routing problems}.
\newblock \bibinfo{journal}{Metaheuristics} , \bibinfo{pages}{407--437}.
%Type = Misc
\bibitem[{Purolator.com, 2021()}]{e-truck-puro}
Purolator.com, 2021, \bibinfo{year}{2021}.
\newblock \bibinfo{note}{Purolator hits the road as first national courier to deploy fully electric delivery vehicles, https://www.purolator.com/en/articles/purolator-hits-road-first-national-courier-deploy-fully-electric-delivery-vehicles}.
%Type = Article
\bibitem[{Qin et~al.(2019)Qin, Tao and Li}]{qin2019vehicle}
\bibinfo{author}{Qin, G.}, \bibinfo{author}{Tao, F.}, \bibinfo{author}{Li, L.}, \bibinfo{year}{2019}.
\newblock \bibinfo{title}{A vehicle routing optimization problem for cold chain logistics considering customer satisfaction and carbon emissions}.
\newblock \bibinfo{journal}{International journal of environmental research and public health} \bibinfo{volume}{16}, \bibinfo{pages}{576}.
%Type = Article
\bibitem[{Reddy and Narayana(2022)}]{reddy2022meta}
\bibinfo{author}{Reddy, A.K.V.K.}, \bibinfo{author}{Narayana, K.V.L.}, \bibinfo{year}{2022}.
\newblock \bibinfo{title}{Meta-heuristics optimization in electric vehicles-an extensive review}.
\newblock \bibinfo{journal}{Renewable and Sustainable Energy Reviews} \bibinfo{volume}{160}, \bibinfo{pages}{112285}.
%Type = Misc
\bibitem[{Ryder.com, 2024()}]{ryder}
Ryder.com, 2024, \bibinfo{year}{2024}.
\newblock \bibinfo{note}{Prioritizing Sustainability in Last Mile Delivery, https://www.ryder.com/en-us/insights/blogs/last-mile/sustainability-last-mile-delivery?}
%Type = Article
\bibitem[{Schneider et~al.(2014)Schneider, Stenger and Goeke}]{schneider2014electric}
\bibinfo{author}{Schneider, M.}, \bibinfo{author}{Stenger, A.}, \bibinfo{author}{Goeke, D.}, \bibinfo{year}{2014}.
\newblock \bibinfo{title}{The electric vehicle-routing problem with time windows and recharging stations}.
\newblock \bibinfo{journal}{Transportation science} \bibinfo{volume}{48}, \bibinfo{pages}{500--520}.
%Type = Article
\bibitem[{Shao et~al.(2017)Shao, Guan, Ran, He and Bi}]{shao2017electric}
\bibinfo{author}{Shao, S.}, \bibinfo{author}{Guan, W.}, \bibinfo{author}{Ran, B.}, \bibinfo{author}{He, Z.}, \bibinfo{author}{Bi, J.}, \bibinfo{year}{2017}.
\newblock \bibinfo{title}{Electric vehicle routing problem with charging time and variable travel time}.
\newblock \bibinfo{journal}{Mathematical Problems in Engineering} \bibinfo{volume}{2017}, \bibinfo{pages}{5098183}.
%Type = Inproceedings
\bibitem[{Shaw(1998)}]{shaw1998using}
\bibinfo{author}{Shaw, P.}, \bibinfo{year}{1998}.
\newblock \bibinfo{title}{Using constraint programming and local search methods to solve vehicle routing problems}, in: \bibinfo{booktitle}{International conference on principles and practice of constraint programming}, \bibinfo{organization}{Springer}. pp. \bibinfo{pages}{417--431}.
%Type = Article
\bibitem[{Siragusa et~al.(2022)Siragusa, Tumino, Mangiaracina and Perego}]{siragusa2022electric}
\bibinfo{author}{Siragusa, C.}, \bibinfo{author}{Tumino, A.}, \bibinfo{author}{Mangiaracina, R.}, \bibinfo{author}{Perego, A.}, \bibinfo{year}{2022}.
\newblock \bibinfo{title}{Electric vehicles performing last-mile delivery in b2c e-commerce: An economic and environmental assessment}.
\newblock \bibinfo{journal}{International Journal of Sustainable Transportation} \bibinfo{volume}{16}, \bibinfo{pages}{22--33}.
%Type = Article
\bibitem[{Sluijk et~al.(2023)Sluijk, Florio, Kinable, Dellaert and Van~Woensel}]{sluijk2023two}
\bibinfo{author}{Sluijk, N.}, \bibinfo{author}{Florio, A.M.}, \bibinfo{author}{Kinable, J.}, \bibinfo{author}{Dellaert, N.}, \bibinfo{author}{Van~Woensel, T.}, \bibinfo{year}{2023}.
\newblock \bibinfo{title}{Two-echelon vehicle routing problems: A literature review}.
\newblock \bibinfo{journal}{European Journal of Operational Research} \bibinfo{volume}{304}, \bibinfo{pages}{865--886}.
%Type = Article
\bibitem[{Srinivas et~al.(2022)Srinivas, Ramachandiran and Rajendran}]{srinivas2022autonomous}
\bibinfo{author}{Srinivas, S.}, \bibinfo{author}{Ramachandiran, S.}, \bibinfo{author}{Rajendran, S.}, \bibinfo{year}{2022}.
\newblock \bibinfo{title}{Autonomous robot-driven deliveries: A review of recent developments and future directions}.
\newblock \bibinfo{journal}{Transportation research part E: logistics and transportation review} \bibinfo{volume}{165}, \bibinfo{pages}{102834}.
%Type = Article
\bibitem[{Tang et~al.(2023)Tang, Zhuang, Li, Liu, Song and Yin}]{TANG2023121711}
\bibinfo{author}{Tang, M.}, \bibinfo{author}{Zhuang, W.}, \bibinfo{author}{Li, B.}, \bibinfo{author}{Liu, H.}, \bibinfo{author}{Song, Z.}, \bibinfo{author}{Yin, G.}, \bibinfo{year}{2023}.
\newblock \bibinfo{title}{Energy-optimal routing for electric vehicles using deep reinforcement learning with transformer}.
\newblock \bibinfo{journal}{Applied Energy} \bibinfo{volume}{350}, \bibinfo{pages}{121711}.
\newblock \URLprefix \url{https://www.sciencedirect.com/science/article/pii/S0306261923010759}, \DOIprefix\doi{https://doi.org/10.1016/j.apenergy.2023.121711}.
%Type = Article
\bibitem[{Trachanatzi et~al.(2020)Trachanatzi, Rigakis, Marinaki and Marinakis}]{trachanatzi2020firefly}
\bibinfo{author}{Trachanatzi, D.}, \bibinfo{author}{Rigakis, M.}, \bibinfo{author}{Marinaki, M.}, \bibinfo{author}{Marinakis, Y.}, \bibinfo{year}{2020}.
\newblock \bibinfo{title}{A firefly algorithm for the environmental prize-collecting vehicle routing problem}.
\newblock \bibinfo{journal}{Swarm and Evolutionary Computation} \bibinfo{volume}{57}, \bibinfo{pages}{100712}.
%Type = Article
\bibitem[{Tsai et~al.(2024)Tsai, Ngo and Che}]{tsai2024last}
\bibinfo{author}{Tsai, J.F.}, \bibinfo{author}{Ngo, H.N.}, \bibinfo{author}{Che, Z.H.}, \bibinfo{year}{2024}.
\newblock \bibinfo{title}{Last-mile delivery during covid-19: A systematic review of parcel locker adoption and consumer experience}.
\newblock \bibinfo{journal}{Acta Psychologica} \bibinfo{volume}{249}, \bibinfo{pages}{104462}.
%Type = Article
\bibitem[{Vinyals et~al.(2015)Vinyals, Fortunato and Jaitly}]{vinyals2015pointer}
\bibinfo{author}{Vinyals, O.}, \bibinfo{author}{Fortunato, M.}, \bibinfo{author}{Jaitly, N.}, \bibinfo{year}{2015}.
\newblock \bibinfo{title}{Pointer networks}.
\newblock \bibinfo{journal}{Advances in neural information processing systems} \bibinfo{volume}{28}.
%Type = Article
\bibitem[{Wang et~al.(2021)Wang, Guo and Zuo}]{wang2021solving}
\bibinfo{author}{Wang, C.}, \bibinfo{author}{Guo, C.}, \bibinfo{author}{Zuo, X.}, \bibinfo{year}{2021}.
\newblock \bibinfo{title}{Solving multi-depot electric vehicle scheduling problem by column generation and genetic algorithm}.
\newblock \bibinfo{journal}{Applied Soft Computing} \bibinfo{volume}{112}, \bibinfo{pages}{107774}.
%Type = Article
\bibitem[{Wang and Zhou(2021)}]{wang2021two}
\bibinfo{author}{Wang, D.}, \bibinfo{author}{Zhou, H.}, \bibinfo{year}{2021}.
\newblock \bibinfo{title}{A two-echelon electric vehicle routing problem with time windows and battery swapping stations}.
\newblock \bibinfo{journal}{Applied Sciences} \bibinfo{volume}{11}, \bibinfo{pages}{10779}.
%Type = Inproceedings
\bibitem[{Wang et~al.(2019)Wang, Zhou and Feng}]{wang2019two}
\bibinfo{author}{Wang, D.}, \bibinfo{author}{Zhou, H.}, \bibinfo{author}{Feng, R.}, \bibinfo{year}{2019}.
\newblock \bibinfo{title}{A two-echelon vehicle routing problem involving electric vehicles with time windows}, in: \bibinfo{booktitle}{Journal of Physics: Conference Series}, \bibinfo{organization}{IOP Publishing}. p. \bibinfo{pages}{012071}.
%Type = Article
\bibitem[{Wang et~al.(2017)Wang, Shao and Zhou}]{wang2017matheuristic}
\bibinfo{author}{Wang, K.}, \bibinfo{author}{Shao, Y.}, \bibinfo{author}{Zhou, W.}, \bibinfo{year}{2017}.
\newblock \bibinfo{title}{Matheuristic for a two-echelon capacitated vehicle routing problem with environmental considerations in city logistics service}.
\newblock \bibinfo{journal}{Transportation Research Part D: Transport and Environment} \bibinfo{volume}{57}, \bibinfo{pages}{262--276}.
%Type = Article
\bibitem[{Wang et~al.(2025)Wang, Adulyasak, Cordeau and He}]{wang2025heterogeneous}
\bibinfo{author}{Wang, W.}, \bibinfo{author}{Adulyasak, Y.}, \bibinfo{author}{Cordeau, J.F.}, \bibinfo{author}{He, G.}, \bibinfo{year}{2025}.
\newblock \bibinfo{title}{The heterogeneous-fleet electric vehicle routing problem with nonlinear charging functions}.
\newblock \bibinfo{journal}{Transportation Research Part C: Emerging Technologies} \bibinfo{volume}{170}, \bibinfo{pages}{104932}.
%Type = Article
\bibitem[{Wang et~al.(2023)Wang, Zhou, Sun, Fan, Wang and Wang}]{wang2023collaborative}
\bibinfo{author}{Wang, Y.}, \bibinfo{author}{Zhou, J.}, \bibinfo{author}{Sun, Y.}, \bibinfo{author}{Fan, J.}, \bibinfo{author}{Wang, Z.}, \bibinfo{author}{Wang, H.}, \bibinfo{year}{2023}.
\newblock \bibinfo{title}{Collaborative multidepot electric vehicle routing problem with time windows and shared charging stations}.
\newblock \bibinfo{journal}{Expert Systems with Applications} \bibinfo{volume}{219}, \bibinfo{pages}{119654}.
%Type = Article
\bibitem[{Wang et~al.(2022)Wang, Zhou, Sun, Wang, Zhe and Wang}]{wang2022electric}
\bibinfo{author}{Wang, Y.}, \bibinfo{author}{Zhou, J.}, \bibinfo{author}{Sun, Y.}, \bibinfo{author}{Wang, X.}, \bibinfo{author}{Zhe, J.}, \bibinfo{author}{Wang, H.}, \bibinfo{year}{2022}.
\newblock \bibinfo{title}{Electric vehicle charging station location-routing problem with time windows and resource sharing}.
\newblock \bibinfo{journal}{Sustainability} \bibinfo{volume}{14}, \bibinfo{pages}{11681}.
%Type = Article
\bibitem[{Wen et~al.(2024)Wen, Lin, Wu, Mao, Cai, Hou, Guo, Liang, Jin, Zhao, Zimmermann, Ye and Wan}]{10639499}
\bibinfo{author}{Wen, H.}, \bibinfo{author}{Lin, Y.}, \bibinfo{author}{Wu, L.}, \bibinfo{author}{Mao, X.}, \bibinfo{author}{Cai, T.}, \bibinfo{author}{Hou, Y.}, \bibinfo{author}{Guo, S.}, \bibinfo{author}{Liang, Y.}, \bibinfo{author}{Jin, G.}, \bibinfo{author}{Zhao, Y.}, \bibinfo{author}{Zimmermann, R.}, \bibinfo{author}{Ye, J.}, \bibinfo{author}{Wan, H.}, \bibinfo{year}{2024}.
\newblock \bibinfo{title}{A survey on service route and time prediction in instant delivery: Taxonomy, progress, and prospects}.
\newblock \bibinfo{journal}{IEEE Transactions on Knowledge and Data Engineering} \bibinfo{volume}{36}, \bibinfo{pages}{7516--7535}.
\newblock \DOIprefix\doi{10.1109/TKDE.2024.3441309}.
%Type = Article
\bibitem[{Wu and Zhang(2023)}]{wu2023branch}
\bibinfo{author}{Wu, Z.}, \bibinfo{author}{Zhang, J.}, \bibinfo{year}{2023}.
\newblock \bibinfo{title}{A branch-and-price algorithm for two-echelon electric vehicle routing problem}.
\newblock \bibinfo{journal}{Complex \& Intelligent Systems} \bibinfo{volume}{9}, \bibinfo{pages}{2475--2490}.
%Type = Article
\bibitem[{Yao et~al.(2022)Yao, Chen and Yang}]{9713755}
\bibinfo{author}{Yao, C.}, \bibinfo{author}{Chen, S.}, \bibinfo{author}{Yang, Z.}, \bibinfo{year}{2022}.
\newblock \bibinfo{title}{Online distributed routing problem of electric vehicles}.
\newblock \bibinfo{journal}{IEEE Transactions on Intelligent Transportation Systems} \bibinfo{volume}{23}, \bibinfo{pages}{16330--16341}.
\newblock \DOIprefix\doi{10.1109/TITS.2022.3149942}.
%Type = Article
\bibitem[{Ye et~al.(2022)Ye, He and Chen}]{ye2022electric}
\bibinfo{author}{Ye, C.}, \bibinfo{author}{He, W.}, \bibinfo{author}{Chen, H.}, \bibinfo{year}{2022}.
\newblock \bibinfo{title}{Electric vehicle routing models and solution algorithms in logistics distribution: A systematic review}.
\newblock \bibinfo{journal}{Environmental Science and Pollution Research} \bibinfo{volume}{29}, \bibinfo{pages}{57067--57090}.
\newblock \DOIprefix\doi{https://doi.org/10.1007/s11356-022-21559-2}.
%Type = Article
\bibitem[{Yu and Puchinger(2024)}]{yu2024collaborative}
\bibinfo{author}{Yu, S.}, \bibinfo{author}{Puchinger, J.}, \bibinfo{year}{2024}.
\newblock \bibinfo{title}{Collaborative truck--robot deliveries: challenges, models, and methods}.
\newblock \bibinfo{journal}{Annals of Operations Research} , \bibinfo{pages}{1--38}.
%Type = Article
\bibitem[{Zhou et~al.(2024)Zhou, Zhang, Yuan, Wang, Zhou and Bell}]{zhou2024pickup}
\bibinfo{author}{Zhou, S.}, \bibinfo{author}{Zhang, D.}, \bibinfo{author}{Yuan, W.}, \bibinfo{author}{Wang, Z.}, \bibinfo{author}{Zhou, L.}, \bibinfo{author}{Bell, M.G.}, \bibinfo{year}{2024}.
\newblock \bibinfo{title}{Pickup and delivery problem with electric vehicles and time windows considering queues}.
\newblock \bibinfo{journal}{Transportation Research Part C: Emerging Technologies} \bibinfo{volume}{167}, \bibinfo{pages}{104829}.
%Type = Article
\bibitem[{Zhou et~al.(2021)Zhou, Huang, Shi, Wang and Huang}]{zhou2021electric}
\bibinfo{author}{Zhou, Y.}, \bibinfo{author}{Huang, J.}, \bibinfo{author}{Shi, J.}, \bibinfo{author}{Wang, R.}, \bibinfo{author}{Huang, K.}, \bibinfo{year}{2021}.
\newblock \bibinfo{title}{The electric vehicle routing problem with partial recharge and vehicle recycling}.
\newblock \bibinfo{journal}{Complex \& Intelligent Systems} \bibinfo{volume}{7}, \bibinfo{pages}{1445--1458}.
%Type = Article
\bibitem[{Zijlstra et~al.(2021)Zijlstra, Galindo~Pecin and Spliet}]{zijlstra2021integrating}
\bibinfo{author}{Zijlstra, T.}, \bibinfo{author}{Galindo~Pecin, D.}, \bibinfo{author}{Spliet, R.}, \bibinfo{year}{2021}.
\newblock \bibinfo{title}{Integrating electric vehicles in a the two-echelon vehicle routing problem} .

\end{thebibliography}

%% else use the following coding to input the bibitems directly in the
%% TeX file.

\end{document}